\documentclass[reqno,a4paper]{amsart}

\usepackage{mathrsfs,amsmath,amssymb,graphicx}
\usepackage{stmaryrd}
\usepackage[unicode]{hyperref}
\usepackage{latexsym}
\usepackage{layout}
\usepackage[english]{babel}
\usepackage{color}
\usepackage{fancyvrb}
\usepackage{color}
\usepackage{amsmath,amscd}
\usepackage{txfonts}

\usepackage{multicol} 

\usepackage{subcaption}

\usepackage{transparent}
\usepackage{xifthen}

\usepackage{tikz}
\usetikzlibrary{calc}
\usetikzlibrary{matrix,arrows,decorations}

\makeatletter
\@namedef{subjclassname@2020}{\textup{2020} Mathematics Subject Classification}
\makeatother

\theoremstyle{plain}
\newtheorem*{theorem*}{Theorem} 
\newtheorem{theorem}{Theorem}[section]
\newtheorem{lemma}[theorem]{Lemma}

\newtheorem{corollary}[theorem]{Corollary}

\theoremstyle{definition}
\newtheorem{definition}{Definition}[section]
 
\theoremstyle{remark}
\newtheorem{remark}{Remark}[section]
\newtheorem{question}{Question}[section]
\newtheorem{problem}{Problem}[section]

\newcommand{\HK}{{\rm HK}}
\newcommand{\HKAtau}{{\rm HK}^{\mathcal{A},\tau}}

\newcommand{\Sbb} {\mathbb S}
\newcommand{\sphere} {\Sbb^3}
\newcommand{\disk}{{\mathcal D}}
\newcommand{\dsystem}{{\mathbf D}}
\newcommand{\id}{\rm id}
 
\newcommand{\bCompl}[1]{\partial E(#1)}
\newcommand{\Compl}[1]{E(#1)}
\newcommand{\ComplHKA}{\Compl{\HK_A}}
\newcommand{\ComplHK}{\Compl \HK}
\newcommand{\mirror}{\mathfrak{m}}

\newcommand{\pair}{(\sphere,\HK)}
\newcommand{\pairA}{(\sphere,\HK_A)}
\newcommand{\pairM}{(\sphere,\HK_M)}
\newcommand{\pairKtau}{(\sphere,\rnbhd{K\cup\tau})}
\newcommand{\pairAtau}{(\sphere,\HK^{\mathcal{A},\tau})}
\newcommand{\pairAtauA}{(\sphere,\HK^{\mathcal{A},\tau}_{A})}
\newcommand{\pairandA}{(\sphere,\HK,A)}
\newcommand{\pairandmA}{(\sphere,\HK,\mirror A)}

\newcommand{\Kpq}{\mathcal{K}_{p,q}}
\newcommand{\pairKpq}{(\sphere,\Kpq)}
\newcommand{\ComplKpq}{\Compl{\Kpq}}

\newcommand{\Tmn}{\mathcal{T}_{\mu,\nu}}
\newcommand{\pairTmn}{(\sphere,\Tmn)}
\newcommand{\ComplTmn}{\Compl{\Tmn}}
\newcommand{\ComplTmnA}{\Compl{\mathcal{T}_{\mu,\nu, A}}} 
\newcommand{\pairTthreethree}{(\sphere,\mathcal{T}_{3,3})} 

\newcommand{\Imn}{\mathcal{I}_{\mu,\nu}}
\newcommand{\pairImn}{(\sphere,\Imn)}
\newcommand{\ComplImn}{\Compl{\Imn}}
\newcommand{\ComplImnA}{\Compl{\mathcal{I}_{\mu,\nu, A}}} 
 
\newcommand{\Tmu}{\mathcal{T}_{\mu,2-\mu}}
\newcommand{\pairTmu}{(\sphere,\Tmu)}
\newcommand{\ComplTmu}{\Compl{\Tmu}}
\newcommand{\ComplTmuA}{\Compl{\mathcal{T}_{\mu,2-\mu, A}}}

\newcommand{\Umn}{\mathcal{U}_{\mu,\nu}}
\newcommand{\pairUmn}{(\sphere,\Umn)}
\newcommand{\ComplUmn}{\Compl{\Umn}}

\newcommand{\Tp}{\mathcal{T}_{p,p}} 
\newcommand{\pairTp}{(\sphere,\Tp)}

\newcommand{\pairHKt}{(\sphere,\HK^t)}
\newcommand{\pairHKtA}{(\sphere,\HK^t_A)}
\newcommand{\HKt}{\HK^t}
\newcommand{\ComplHKt}{\Compl{\HKt}}
\newcommand{\ComplHKtA}{\Compl{\HKt_A}}

\newcommand{\HKc}{\HK^c}
\newcommand{\pairHKc}{(\sphere,\HKc)}
\newcommand{\pairHKcA}{(\sphere,\HKc_A)}
\newcommand{\pairHKcAM}{(\sphere,\HKc_{A,M})} 
\newcommand{\ComplHKc}{\Compl{\HKc}}
\newcommand{\ComplHKcA}{\Compl{\HKc_A}}
 
\newcommand{\lk}[2]{{\ell \mathit k}(#1,#2)} 
\newcommand{\IN}[2]{{\mathcal I}([#1],[#2])} 
\newcommand{\In}{{\mathcal I}}

\newcommand{\cl}[1]{{\rm cl}(#1)}
 
\newcommand{\cout}[1]   {}

\newcommand{\op}[1]{\operatorname{#1}}

\newcommand{\rnbhd}[1]{\mathfrak N(#1)} 
\newcommand{\opennbhd}[2]{\mathring{\mathfrak N}(#1; #2)}
\newcommand{\openrnbhd}[1]{\mathring{\mathfrak N}(#1)}
\newcommand{\rnbhdKpq}{\rnbhd{\Kpq}}

\newcommand{\Sym}[2][\sphere]{\mathcal MCG(#1, #2)}
\newcommand{\pSym}[2][\sphere]{\mathcal MCG_+(#1, #2)}
\newcommand{\gSym}[2][\sphere]{\mathcal MCG_{(+)}(#1, #2)}
\newcommand{\MCG}[1]{\mathcal MCG(#1)}
\newcommand{\pMCG}[1]{\mathcal MCG_+(#1)}
\newcommand{\Aut}[1]{\mathcal Homeo(#1)}
\newcommand{\pAut}[1]{\mathcal Homeo_+(#1)}

\newcommand{\TSG}[1]{{\rm TSG}(\sphere, #1)}
 
\newcommand{\rel}{{\rm rel\,}} 
\newcommand{\absp}{\vert p \vert}

\definecolor{mygray}{rgb}{0.92,0.92,0.92}

\numberwithin{equation}{section}
\numberwithin{figure}{section}

\title[Cylindrical handlebody-knots]{Rigidity and symmetry of cylindrical handlebody-knots}
 
\author{Yi-Sheng Wang}
\address{Institute of Mathematics, Academia Sinica, Taipei City 106, Taiwan}
\email{yisheng@gate.sinica.edu.tw}

\date{\today}

\begin{document}

\thanks{The work was supported by
National Center for Theoretical Sciences, Academia Sinica, 
and MoST (grant no. 110-2115-M-001-004-MY3), Taiwan}

\begin{abstract} 
A recent result of Funayoshi-Koda shows
that a handlebody-knot of genus two 
has a finite symmetry group if and only if
it is hyperbolic---the exterior admits a hyperbolic structure 
with totally geodesic boundary---or irreducible, atoroidal, 
cylindrical---the exterior contains no essential disks or tori
but contains an essential annulus.
Based on the Koda-Ozawa
classification theorem, essential annuli 
in an irreducible, atoroidal handlebody-knots of genus two
are classified into four classes:
type $2$, type $3$-$2$, type $3$-$3$ and type $4$-$1$.
We show that under mild conditions
most genus two cylindrical handlebody-knot exteriors
contain no essential disks or tori, 
and when a type $3$-$3$ annulus exists, 
it is often unique up to isotopy; 
a classification result for symmetry groups of such cylindrical handlebody-knots
is also obtained.
\end{abstract}

\maketitle
 
\section{Introduction}\label{sec:intro}
By Thurston's hyperbolization theorem, 
knots are classified into four categories:
trivial, torus, satellite and hyperbolic knots,
based on the existence of
essential surfaces with non-negative Euler characteristics 
in knot exteriors. 
In particular, a torus knot is characterized by 
the existence of an essential annulus
and absence of any essential disk or torus
in its exterior. Furthermore, essential annuli
in a torus knot exterior are all isotopic as shown in Tsau \cite{Tsa:94}. 
The uniqueness of essential annulus implies 
a classical result of Schreier \cite{Sch:29}
which states that torus knots are chiral and their symmetry groups are all 
isomorphic to $\mathbb{Z}_2$
as explained in Section \ref{subsec:torus_knot}.

A genus $g$ handlebody-knot $\pair$ is a
genus $g$ handlebody $\HK$ embedded in an oriented $3$-sphere $\sphere$;
the study of genus one handlebody-knots is equivalent to classical knot theory. The present work is concerned with genus two
handlbody-knots, abbreviated to handlebody-knots
hereafter, unless otherwise specified. 
As with the case of knots, by Thurston's hyperbolization theorem,
together with 
the equivariant torus theorem by Holzmann \cite{Hol:91} 
and the fixed point theorem by Tollefson \cite{Tol:81},
handlebody-knots are classified into four classes:
\begin{multicols}{2}
\begin{itemize}
\item reducible 
\item irreducible, toroidal 
\item irreducible, atoroidal, cylindrical  
\item hyperbolic  
\end{itemize}
\end{multicols}
A reducible handlebody-knot is a handlebody-knot
whose exterior $\Compl\HK:=\overline{\sphere-\HK}$
contains an essential disk---this should be contrasted 
with the trivial knot,
whereas a cylindrical (resp.\ toroidal) handlebody-knot is 
characterized by the existence of an
essential annulus (resp.\ torus) in its exterior. 
Particularly, irreducible, \emph{toroidal} handlebody-knots
correspond to satellite knots, and
irreducible, atoroidal, \emph{cylindrical} handlebody-knots to torus knots. A handlebody-knot is hyperbolic if its exterior
admits a hyperbolic metric with totally geodesic boundary.
In the present paper, we study \emph{cylindrical} handlebody-knots 
and their symmetries. 

\begin{figure}[h]
\centering
\def\svgwidth{.35\columnwidth}
\begingroup%
  \makeatletter%
  \providecommand\color[2][]{%
    \errmessage{(Inkscape) Color is used for the text in Inkscape, but the package 'color.sty' is not loaded}%
    \renewcommand\color[2][]{}%
  }%
  \providecommand\transparent[1]{%
    \errmessage{(Inkscape) Transparency is used (non-zero) for the text in Inkscape, but the package 'transparent.sty' is not loaded}%
    \renewcommand\transparent[1]{}%
  }%
  \providecommand\rotatebox[2]{#2}%
  \newcommand*\fsize{\dimexpr\f@size pt\relax}%
  \newcommand*\lineheight[1]{\fontsize{\fsize}{#1\fsize}\selectfont}%
  \ifx\svgwidth\undefined%
    \setlength{\unitlength}{1133.85826772bp}%
    \ifx\svgscale\undefined%
      \relax%
    \else%
      \setlength{\unitlength}{\unitlength * \real{\svgscale}}%
    \fi%
  \else%
    \setlength{\unitlength}{\svgwidth}%
  \fi%
  \global\let\svgwidth\undefined%
  \global\let\svgscale\undefined%
  \makeatother%
  \begin{picture}(1,0.75)%
    \lineheight{1}%
    \setlength\tabcolsep{0pt}%
    \put(0,0){\includegraphics[width=\unitlength,page=1]{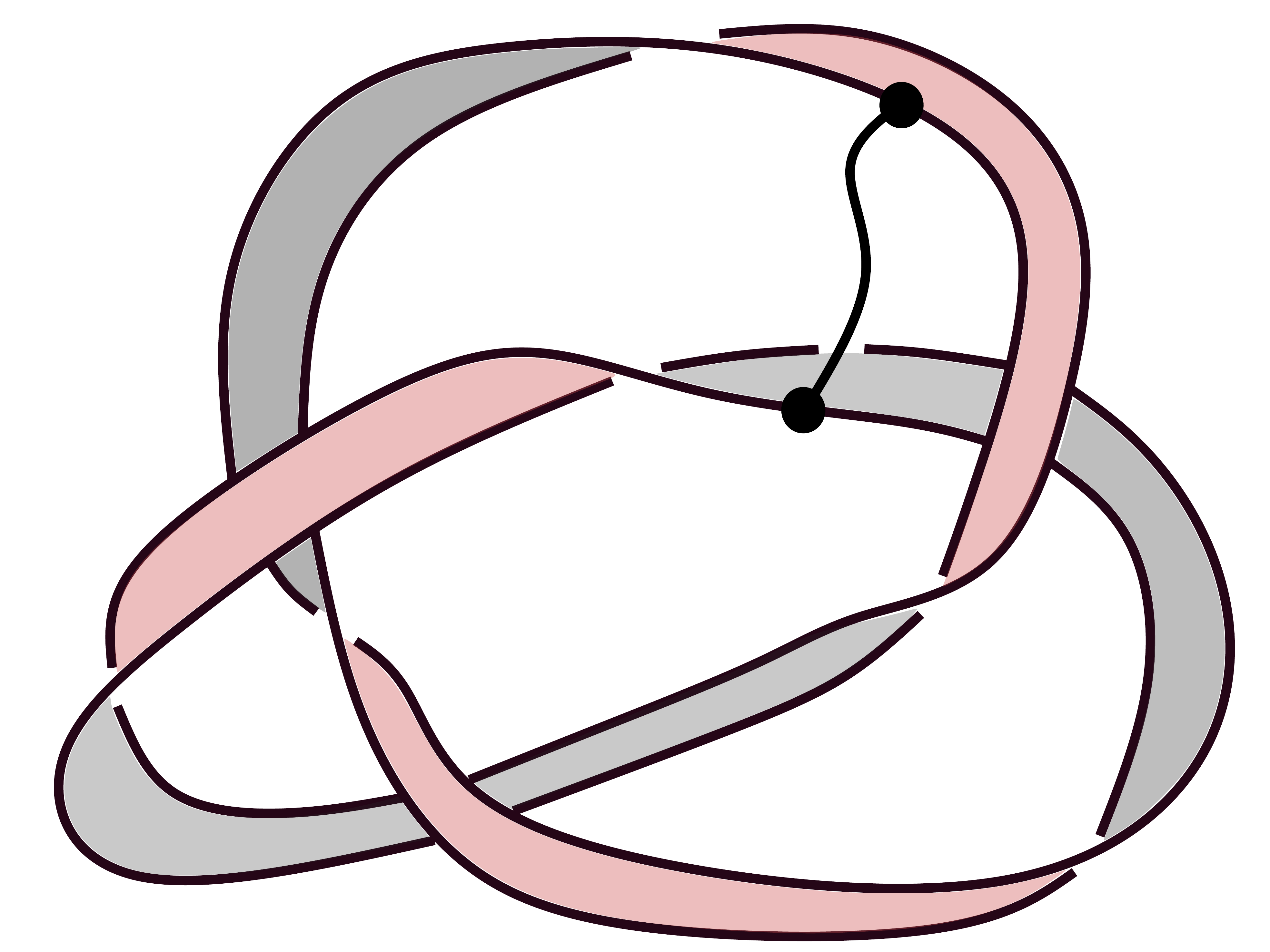}}%
    \put(0.26048897,0.37648033){\color[rgb]{0,0,0}\makebox(0,0)[lt]{\lineheight{1.25}\smash{\begin{tabular}[t]{l}{\footnotesize $A$}\end{tabular}}}}%
    \put(0.90414185,0.16837715){\color[rgb]{0,0,0}\makebox(0,0)[lt]{\lineheight{1.25}\smash{\begin{tabular}[t]{l}{\footnotesize $A$}\end{tabular}}}}%
  \end{picture}%
\endgroup%

\caption{Reducible $\pair$ with $\pairA$ trivial.} 
\label{fig:intro:hk_reducible_hkA_trivial}
\end{figure}

Contrary to the case of knots, 
the existence of essential annulus (resp.\
essential torus) does not entail the 
non-existence of essential disks in a handlebody-knot exterior, 
and there are many
reducible, cylindrical handlebody-knots (e.g.\ the handlebody-knot in Fig.\ \ref{fig:intro:hk_reducible_hkA_trivial}). 
On the other hand, akin to the case of knots,
irreducible, \emph{toroidal, cylindrical} handlebody-knots,
the analogue of cable knots,
abound. This raises the following question.
\begin{question}\label{ques:irre_atoro}
When is a cylindrical handlebody-knot $\pair$ 
irreducible and atoroidal?
\end{question}

\begin{figure}[h]
\begin{subfigure}{0.45\textwidth}
\centering
\def\svgwidth{.7\columnwidth}
\begingroup%
  \makeatletter%
  \providecommand\color[2][]{%
    \errmessage{(Inkscape) Color is used for the text in Inkscape, but the package 'color.sty' is not loaded}%
    \renewcommand\color[2][]{}%
  }%
  \providecommand\transparent[1]{%
    \errmessage{(Inkscape) Transparency is used (non-zero) for the text in Inkscape, but the package 'transparent.sty' is not loaded}%
    \renewcommand\transparent[1]{}%
  }%
  \providecommand\rotatebox[2]{#2}%
  \newcommand*\fsize{\dimexpr\f@size pt\relax}%
  \newcommand*\lineheight[1]{\fontsize{\fsize}{#1\fsize}\selectfont}%
  \ifx\svgwidth\undefined%
    \setlength{\unitlength}{1133.85826772bp}%
    \ifx\svgscale\undefined%
      \relax%
    \else%
      \setlength{\unitlength}{\unitlength * \real{\svgscale}}%
    \fi%
  \else%
    \setlength{\unitlength}{\svgwidth}%
  \fi%
  \global\let\svgwidth\undefined%
  \global\let\svgscale\undefined%
  \makeatother%
  \begin{picture}(1,0.75)%
    \lineheight{1}%
    \setlength\tabcolsep{0pt}%
    \put(0,0){\includegraphics[width=\unitlength,page=1]{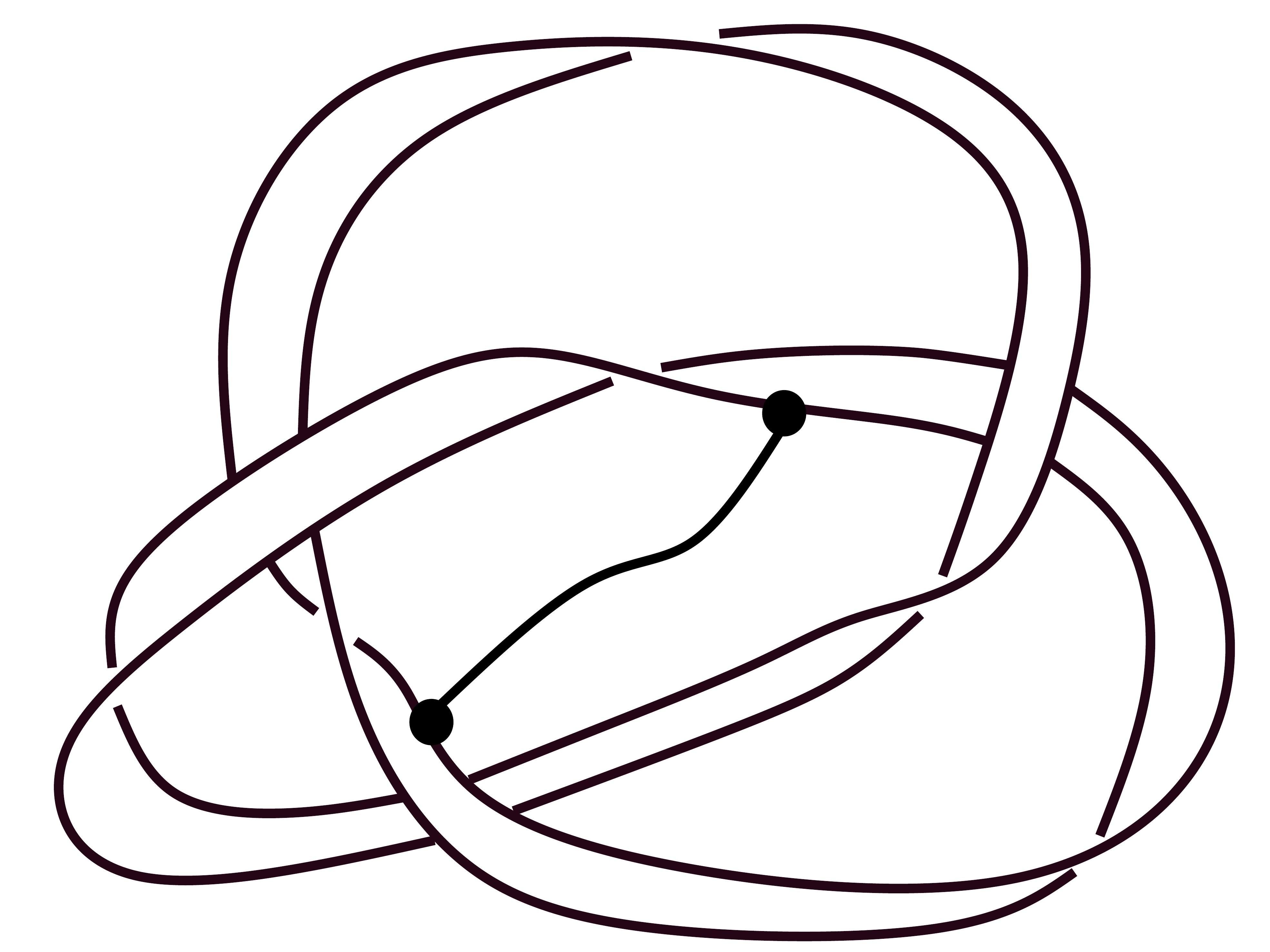}}%
    \put(0.00786382,0.70282885){\color[rgb]{0,0,0}\makebox(0,0)[lt]{\lineheight{1.25}\smash{\begin{tabular}[t]{l}{\footnotesize $\pairHKt$}\end{tabular}}}}%
  \end{picture}%
\endgroup%

\caption{$\partial A$ is a torus link.}
\label{fig:intro:example_two_annuli_torus_link}
\end{subfigure}
\begin{subfigure}{0.5\textwidth}
\centering
\def\svgwidth{.7 \columnwidth}
\begingroup%
  \makeatletter%
  \providecommand\color[2][]{%
    \errmessage{(Inkscape) Color is used for the text in Inkscape, but the package 'color.sty' is not loaded}%
    \renewcommand\color[2][]{}%
  }%
  \providecommand\transparent[1]{%
    \errmessage{(Inkscape) Transparency is used (non-zero) for the text in Inkscape, but the package 'transparent.sty' is not loaded}%
    \renewcommand\transparent[1]{}%
  }%
  \providecommand\rotatebox[2]{#2}%
  \newcommand*\fsize{\dimexpr\f@size pt\relax}%
  \newcommand*\lineheight[1]{\fontsize{\fsize}{#1\fsize}\selectfont}%
  \ifx\svgwidth\undefined%
    \setlength{\unitlength}{1133.85826772bp}%
    \ifx\svgscale\undefined%
      \relax%
    \else%
      \setlength{\unitlength}{\unitlength * \real{\svgscale}}%
    \fi%
  \else%
    \setlength{\unitlength}{\svgwidth}%
  \fi%
  \global\let\svgwidth\undefined%
  \global\let\svgscale\undefined%
  \makeatother%
  \begin{picture}(1,0.775)%
    \lineheight{1}%
    \setlength\tabcolsep{0pt}%
    \put(0,0){\includegraphics[width=\unitlength,page=1]{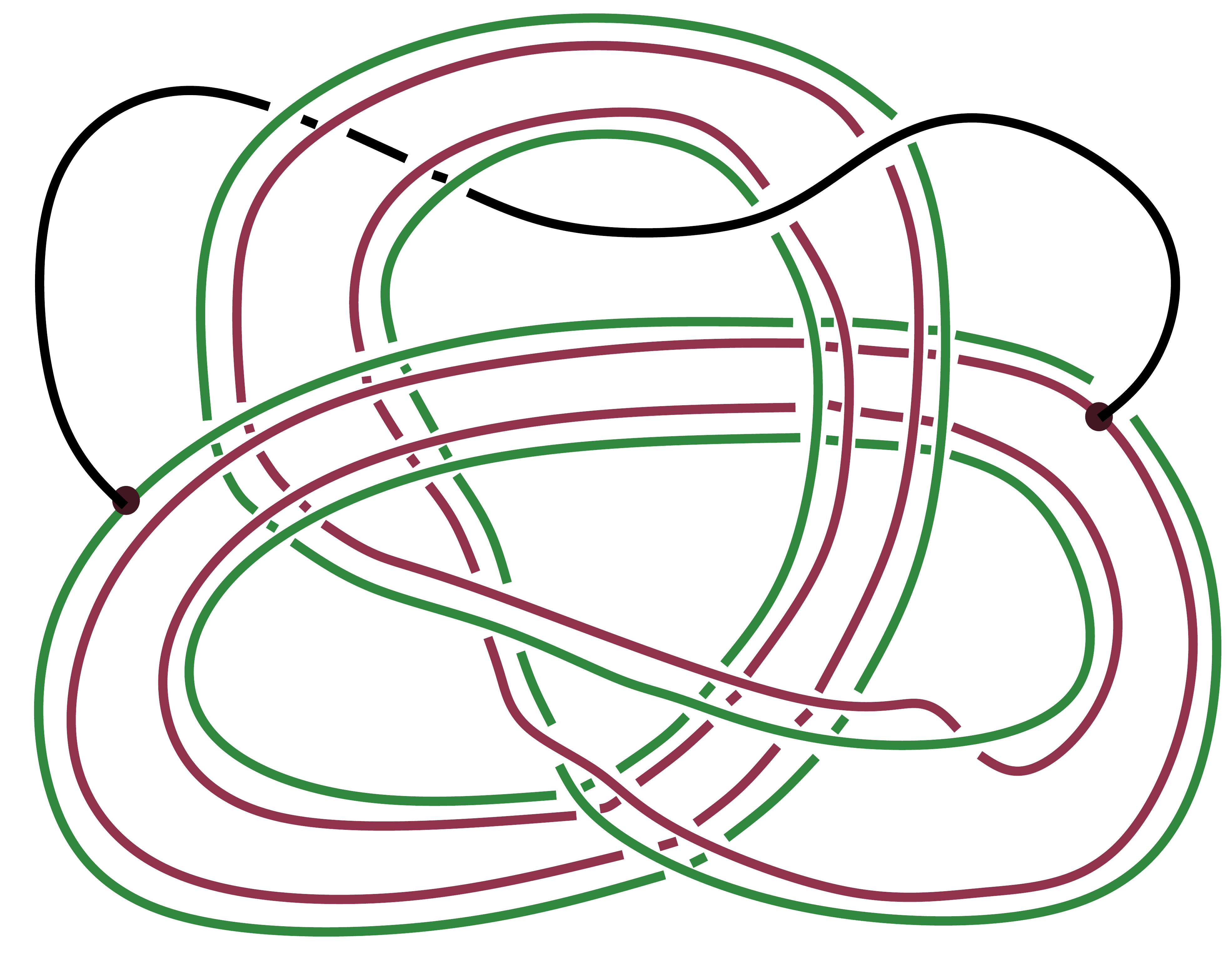}}%
    \put(0.67750946,0.73659466){\color[rgb]{0,0,0}\makebox(0,0)[lt]{\lineheight{1.25}\smash{\begin{tabular}[t]{l}{\footnotesize $\pairHKc$}\end{tabular}}}}%
  \end{picture}%
\endgroup%

\caption{$\partial A$ is a cable link.}
\label{fig:intro:example_two_annuli_cable_link}
\end{subfigure}
\caption{Handlebody-knot exteriors that admit non-isotopic type $3$-$3$ annuli.}
\label{fig:intro:nonuniqueness}
\end{figure}

In contrast to a torus knot,
an irreducible, atoroidal, cylindrical handlebody-knots
may admit non-isotopic essential annuli in its exterior, for instance. handlebody-knots $4_1,6_{10}$ 
in the Ishii-Kishimoto-Moriuchi-Suzuki handlebody-knot table 
\cite{IshKisMorSuz:12}, and handlebody-knots in Fig.\ \ref{fig:intro:nonuniqueness} and Section \ref{subsec:example_non_uniqueness}; 
as we should see later, in some cases irreducible,
atoroidal, cylindrical handlebody-knots behaves more like torus links.
Nonetheless, the existence of an essential
annulus does often impose strong constraints
on the existence of other non-isotopic essential annuli
as observed by Funayoshi-Koda \cite{FunKod:20}.
This leads to the next question. 
\begin{question}\label{ques:uniqueness}
When does 
an irreducible, atoroidal, cylindrical handlebody-knot $\pair$ 
admit a unique essential annulus in its exterior, up to isotopy?
\end{question}   
The uniqueness of essential annuli provides 
rigidity that allows us to compute the (positive) symmetry group
of $\pair$, and in many instances, it reduces 
the computation to studying spatial graph symmetries \cite{Sim:86}, \cite{ChoKod:13}, \cite{Kod:15}.
The symmetry group $\Sym\HK$ of $\pair$ 
is defined as the group of components 
\[
 \pi_0\big(\Aut{\sphere,\HK}\big)
\]
of the topological group of self-homeomorphisms of $\sphere$
preserving $\HK$ setwise, whereas the positive symmetry group
$\pSym\HK$ is the subgroup of $\Sym\HK$ given by
the components of topological subgroup 
$\pAut{\sphere,\HK}$ of orientation-preserving homeomorphisms
in $\Aut{\sphere,\HK}$. Note that when $\pair$ is trivial---namely,
$\Compl\HK$ is a handlebody, 
$\Sym\HK$ is the genus two Goeritz group \cite{Goe:33}. 

It follows from Funayoshi-Koda \cite{FunKod:20} that the symmetry group of a handlebody-knot is finite if 
and only if it is hyperbolic or 
irreducible, atoroidal, cylindrical (see also \cite[Remark $2.1$]{Wan:21}). 
Some examples with a trivial symmetry group 
are computed by Koda \cite{Kod:15} and the author \cite{Wan:21}, but apart from them, little is known about
the structure of these finite groups, in contrast to finite symmetry groups of knots, which are cyclic or dihedral \cite{Kaw:96}.  
This leads to the following classification problem, which, together with Questions \ref{ques:irre_atoro}
and \ref{ques:uniqueness}, is what the present study and its sequels seek to address. 
\begin{problem}\label{prob:symmetry}
Classify the structures of symmetry groups of 
irreducible, atoroidal, cylindrical handlebody-knots 
whose exterior contain a unique essential annulus, up to isotopy.
\end{problem}

Based on the classification theorem by Koda-Ozawa 
\cite{KodOzaGor:15}, essential annuli in a handlebo\-dy-knot exterior
are classified into seven types; 
as observed by Funayoshi-Koda \cite[Lemma $3.2$]{FunKod:20} only four types among the seven,
that is, types $2$, $3$-$2$, $3$-$3$ and $4$-$1$,
can exist in an irreducible, atoroidal
handlebody-knot exterior. 
On the other hand, 
the existence of essential annuli
of one of these four types 
does \emph{not} imply the irreducibility and atoroidality in general. 

These four types of essential annuli 
are characterized by their boundary 
in relation to the handlebody $\HK$. 
A type $2$ annulus $A$ has exactly one component of 
$\partial A$ bounds a disk in $\HK$, while
a type $3$-$2$ (resp.\ type $3$-$3$) annulus 
has parallel (resp.\ non-parallel) 
boundary components $\partial A$ in $\partial \HK$ that bound no disk in $\HK$, and there exists an essential disk $\disk\subset \HK$ 
disjoint from $A$.
A type $4$-$1$ annulus also has parallel $\partial A$ in $\partial \HK$, but 
no essential disk $\disk\subset \HK$ disjoint from $A$ exists. 

Let $\HK_A$ be the union of $\HK$ and a regular neighborhood
$\rnbhd{A}$ of $A$ in $\Compl\HK$. The case where $\HK_A$
is a handlebody is of particular interest, 
as the handlebody-knot
$\pairA$ is often ``simpler" than $\pair$;
in many cases, $\pairA$ is in fact trivial, and
we call $A$ 
an \emph{unknotting annulus}, following Koda \cite{Kod:15},
in such a situation.

Since the boundary components of 
a type $3$-$2$ or type $4$-$1$ annulus are parallel in $\partial\HK$,
$\HK_A$ is a handlebody only if 
$A$ is of type $2$ or type $3$-$3$.
Conversely, if $A$ is of type $2$, then $\HK_A$ 
is always a handlebody. 
The situation with type $3$-$3$ annuli is slightly more involved.
Given a handlebody-knot $\pair$ and 
a type $3$-$3$ annulus $A\subset\Compl\HK$,
by the definition,
there exists an essential disk $\disk_A\subset\HK$ 
disjoint from $A$. 
The disk $\disk_A$
is necessarily separating and unique, up to isotopy, in $\HK$ 
(see Section \ref{subsec:typethreethree}). It separates $\HK$
into two solid tori $W_1,W_2$; each meets $A$ at a component of
$\partial A$. Let $l_i:=\partial A\cap W_i$, $i=1,2$.
Then we define the \emph{slope pair} of $A$
to be the unordered pair $\{r_1,r_2\}$ of rational numbers 
with $r_i$ being the slope of $l_i\subset W_i$, $i=1,2$.

In Section \ref{sec:prelim}, we show that the slope pair 
$\{r_1,r_2\}$ of $A$ 
is either of the form $\{\frac{p}{q},\frac{q}{p}\}$, $pq\neq 0$,
or of the form $\{\frac{p}{q},pq\}$, $q\neq 0$, where $p,q\in \mathbb{Z}$,
and $\HK_A$ is a handlebody if and only 
if it is the latter.
The present paper is concerned primarily with the case where 
$A$ is of type $3$-$3$ with a slope pair $(p,p)$, namely, 
$q=1$, and $p\neq 0$. For the sake of simplicity, in this case
we say $A$ has a non-trivial boundary slope of $p$.

Question \ref{ques:irre_atoro} is investigated in Section \ref{sec:irre_atoro}, where we
deal with irreducibility and atoroidality separately,
and consider general annuli as well as M\"obius bands. 
To state the result, 
we note that $\rnbhd{A}\cap\partial\HK_A$
consists of two annuli $A_+,A_-$ 
parallel to $A$ in $\ComplHK$,
and denote by $l_+,l_-$ essential loops of $A_+,A_-$, respectively.
Then as a corollary of 
Lemmas \ref{lm:irreducibility:boundary_irreducible_ComplHKA},
\ref{lm:irreducibility:trivial_HKA}, \ref{lm:atoroidality:trivial_HKA} 
and \ref{lm:atoroidality_general},
we have the following application to type $3$-$3$ annuli.

\begin{theorem}\label{teo:intro:irre_atoro}
Let $A$ be a type $3$-$3$ annulus with a slope pair $\{\frac{p}{q},pq\}$;
if in addition $q=\pm 1$, we assume that $\partial A\subset\sphere$ is not a $(2m,2n)$-torus link with $\vert m\vert, n >1$ and $mn=\pm p$.
Suppose one of the following holds:
\begin{enumerate}
\item $\pairA$ is irreducible and atoroidal; 
\item $\pairA$ is trivial, $q=\pm 1$, and
neither of $l_+,l_-$ represents 
the $\vert p\vert$-th power of some primitive element of $\pi_1(\ComplHKA)$, up to conjugation.
\end{enumerate} 
Then $\pair$ is irreducible and atoroidal.
\end{theorem}
 
We remark that 
the condition of $\partial A$ being
not a $(2m,2n)$-torus link, 
$\vert m\vert, n >1$ is used to guarantee 
\emph{atoroidality} when $\pairA$ is irreducible,
while it is there to ensure \emph{irreducibility} 
when $\pairA$ is trivial.
The handlebody-knot $\pair$ and type $3$-$3$ annulus $A$
in Fig.\ \ref{fig:intro:hk_reducible_hkA_trivial} 
is a reducible handlebody-knot with $\pairA$ trivial.
Note that $\partial A\subset \sphere$ is a $(6,4)$-torus link
with $3\cdot 2$ being the boundary slope of $A$.
The handlebody-knot $\pairTmn$ 
in Fig.\ \ref{fig:example_Tmn} with $\mu=\pm 1$  
is another such example; in this case, 
one of $l_+,l_-$ represents the $\vert p\vert$-th power of 
some primitive element of $\pi_1(\ComplHKA)$, up to conjugation, where $p$ is the boundary slope of $A$. 

In Section \ref{subsec:examples_irreatoro}, 
we construct several families of 
handlebody-knots whose exterior admit a type $3$-$3$
annulus, and use the criteria developed in Section \ref{subsec:criteria_irre_atoro} to examine their 
irreducibility and atoroidality.
Other methods for detecting irreducibility of a handlebody-knot
are developed by Ishii-Kishimoto \cite{IshKis:11} using quandle invariant, and by Bellettini-Paolini-Wang \cite{BelPaoWan:20} via homomorphisms on fundamental group,
and by Okazaki \cite{Oka:20} using Alexander polynomial.

The uniqueness problem (Question \ref{ques:uniqueness}) 
is studied in Section \ref{sec:uniqueness} 
where attention
is restricted to essential annuli with a non-trivial boundary slope of $p$.
In general, it is not difficult to construct a
handlebody-knot whose exterior contains two 
non-isotopic type $3$-$3$ annuli.
One way is to start with  
a $(2m,2n)$-torus link with $\vert m\vert, n >1$ or a $(2m,2n)$-cable link with $n>1$,
and then choose an arc that connects the two components
but does not intersect the two non-isotopic annuli in the link exterior.
For instance, the exterior of each handlebody-knot in Fig.\ \ref{fig:intro:nonuniqueness} 
contains two non-isotopic type $3$-$3$ annuli.
To verify the resulting handlebody-knot
is irreducible and atoroidal, however, is often a harder undertaking.
The irreducibility and atoroidality of the handlebody-knots in Fig.\ \ref{fig:intro:nonuniqueness} are detected by a strengthening
of Theorem \ref{teo:intro:irre_atoro}
as explained in Section \ref{subsubsec:irre_atoro_of_non_uniqueness_examples}. 

On the other hand, the condition of $\partial A$ being not 
a $(2m,2n)$-torus link, $\vert m\vert, n >1$, 
or a $(2m,2n)$-cable link, $n>1$, with $mn=p$
turns out to suffice to entail 
the uniqueness of most type $3$-$3$ 
essential annuli with a non-trivial boundary slope of $p$. 
Note that the condition is equivalent to saying that
$l_1,l_2$ are not $(m,n)$-torus or -cable knots in $\sphere$ with $mn=p$. 
Here by an $(m,n)$-torus 
knot, we understand a non-trivial torus knot, that is,
with our convention, $\vert m\vert,n>1$; similarly, an $(m,n)$-cable
knot always means a non-trivial cabling, namely $n>1$.
We denote by $\dagger$, $\dagger\dagger$ 
the following conditions, respectively:
\begin{align}
\label{cond:torus_cable}\tag{$\dagger$}
&\text{$l_1,l_2$ are not $(m,n)$-torus or $(m,n)$-cable knots
in $\sphere$ with 
$mn=p$};\\
\label{cond:torus}\tag{$\dagger\dagger$}
&\text{$l_1,l_2$ are not $(m,n)$-torus knots in $\sphere$ with $
mn=p$}.
\end{align}
As will be made clear in the proofs of 
Lemmas \ref{lm:separating_disjoint}, \ref{lm:trivial_intersection_not_primitive}, 
and \ref{lm:parallel}, 
when $A$ is unknotting, only
the latter, weaker condition \eqref{cond:torus}
is required owing to the classification of tunnel number one 
non-simple knots and links by Morimoto-Sakuma \cite{MorSak:91}, Eudave-Mu\~{n}oz-Uchida\cite{EudUch:96}, respectively. 
Recall that a simple loop $l$ in the boundary of a handlebody $V$ 
is \emph{primitive} if there exists a meridian disk $D\subset V$
such that $l\cap D$ is a point. The following summarizes
Theorems \ref{teo:uniqueness_HKA_irreducible} and
\ref{teo:uniqueness_HKA_trivial_not_primitive}.

\begin{theorem}\label{teo:intro:uniqueness_1}
Let $\pair$ be an irreducible, atoroidal handlebody-knots
and $A\subset \ComplHK$ a type $3$-$3$ annulus with a non-trivial boundary slope of $p$.
Suppose $A$ satisfies the condition \eqref{cond:torus_cable}
or the condition \eqref{cond:torus} if $A$ is unknotting, 
and one of the following holds:
\begin{enumerate}
\item $\vert p\vert =1$; 
\item $\vert p\vert >1$, none of $l_+,l_-$
represents the $\vert p\vert$-th power of some element in $\pi_1(\ComplHKA)$,
up to conjugation, 
and if $A$ is unknotting, at least one of $l_+,l_-\subset \ComplHKA$ is not primitive. 
\end{enumerate}
Then, up to isotopy, $A$ is the unique type $3$-$3$
annulus in $\ComplHK$. 
\end{theorem}
Section \ref{subsec:example_non_uniqueness} gives
an irreducible, atoroidal handlebody-knot that
fails the condition of $l_+,l_-$ not representing  
the $\vert p\vert$-th powers of some elements in $\pi_1(\ComplHKA)$
and has two non-isotopic type $3$-$3$ annuli in its exterior. 
On the other hand, the assumption 
that one of $l_+,l_-\subset\ComplHKA$ is not primitive can be dropped
when $\vert p\vert$ is odd and greater than $1$.

\begin{theorem}[\textbf{Theorem \ref{teo:uniqueness_HKA_trivial_odd_p}}]\label{teo:intro_uniqueness_2}
Let $\pair$ be an irreducible, atoroidal handlebody-knot  
and $A\subset \ComplHK$ a type $3$-$3$ annulus with a non-trivial boundary slope of $\absp>1$ and $p$ is odd. 
Suppose $A$ is unknotting and satisfies the condition \eqref{cond:torus}, and neither of $l_+,l_-$ represents the $\vert p\vert$-th multiple of some generator 
of $H_1(\ComplHKA)$. 

Then, up to isotopy, $A$ is the unique type $3$-$3$
annulus in $\ComplHK$.
\end{theorem}

Lastly, Problem \ref{prob:symmetry} is addressed 
in Section \ref{sec:symmetry}, 
where we obtain a classification theorem
for symmetry groups 
of irreducible, atoroidal 
handlebody-knots whose exteriors contain a unique type $3$-$3$ annulus $A$ 
with a non-trivial boundary slope of $p$.

To state the classification result, we 
associate an order pair $(p_1,p_2)$ to
the annulus $A$, 
called the slope type of $A$. The pair 
$(p_1,p_2)$ are 
characterized by the properties: 
$[l_+]=(p_1,p_2), [l_-]=(p_1-1,p_2+1)$ in terms of 
a basis of $H_1(\ComplHKA)$ 
induced by meridian disks of $\HK_A$ disjoint from
$\disk_A$, $p_1>p_2$,
and either $0<p_1\leq p$ or $p<p_1\leq 0$,
depending on the sign of $p$.
We show in Section \ref{subsec:interinsic_disks}
the slope type is well-defined and depends only on $A$ and $\pair$.

\begin{theorem}\label{teo:intro:symmetry}
Let $\pair$ be an irreducible, atoroidal handlebody-knot
whose exterior $\ComplHK$ contains a unique type $3$-$3$
annulus $A$ with a non-trivial slope of $p$.
Then $\Sym\HK=\pSym\HK\leq\mathbb{Z}_2\times \mathbb{Z}_2$.
If in addition the slope type of $A$
is not $(\frac{p+1}{2},\frac{p-1}{2})$, 
then $\Sym\HK=\pSym\HK\leq\mathbb{Z}_2$.
\end{theorem}

Sections \ref{subsec:example_full} 
and \ref{subsec:examples_not_isomorphism}
compute the symmetry groups of 
several families of handlebody-knots, showing that 
$\mathbb{Z}_2\times \mathbb{Z}_2$, $\mathbb{Z}_2$ 
in Theorem \ref{teo:intro:symmetry} are optimal upper bounds, 
and the inequalities $\leq$'s there are in general not an isomorphism. 
As an application of Theorem \ref{teo:intro:symmetry}, 
we obtain the symmetry group
of $(\sphere,5_2)$ and $(\sphere,6_4)$ in the 
Ishii-Kishimoto-Moriuchi-Suzuki handlebody-knot table \cite{IshKisMorSuz:12}:
\[
\Sym{5_2}\simeq \pSym{5_2}\simeq \mathbb{Z}_2\times \mathbb{Z}_2,\quad\Sym{6_4}\simeq \pSym{6_4}\simeq \mathbb{Z}_2.\] 
Note that  
the symmetry groups of $(\sphere, 5_1)$ and $(\sphere, 6_1)$ 
are computed in Koda
\cite{Kod:15} using results from Motto \cite{Mott:90}, Lee-Lee \cite{LeeLee:12}, and they are $\Sym{5_1}=\Sym{6_1}=1$.   
The symmetry groups of the rest handlebody-knots up to six crossings 
seem to remain unknown. 
In fact, there are still four handlebody-knots in the table
whose chiraliry is yet to be determined 
as summarized in Ishii-Iwakiri-Jang-Oshiro \cite[Table $2$]{IshIwaJanOsh:13} (see also \cite{Mott:90}, \cite{IshIwa:12}, \cite{LeeLee:12}, \cite{IshKisOza:15}).

The paper is organized as follows: 
Section \ref{sec:prelim} fixes the notation
and summarizes relevant known results. Section \ref{sec:disks}
discusses natural basis of $H_1(\ComplHKA)$ and of $H_1(\bCompl{\HK_A})$ associated to a type $3$-$3$ annulus $A$ of a non-trivial slope, and
examines the existence and 
non-existence of various types of disks in $\ComplHKA$ when $A$ is unknotting. Results in 
Section \ref{sec:disks} are crucial 
for our investigation on Questions \ref{ques:irre_atoro} and \ref{ques:uniqueness} and Problem \ref{prob:symmetry} 
in Sections \ref{sec:irre_atoro}, \ref{sec:uniqueness} and \ref{sec:symmetry}, respectively. We include many examples and counterexamples along the way, in hope that 
they can provide a more comprehensive picture of
the topic and pave a way toward a complete classification
for symmetry groups of irreducible, atoroidal, cylindrical handlebody-knots.


\section{Preliminaries}\label{sec:prelim}

Throughout the paper we work in the piecewise linear category.
Given a subpolyhedron $X$ of $M$, 
$\mathring{X}$ denotes the interior of $X$, and
$\mathfrak{N}(X)$ a regular neighborhood of $X$ in $M$.
The exterior $\Compl X$ of $X$ in $M$ is the complement of $\opennbhd{X}{M}$ if $X$ has codimension greater than zero, and  
is the closure of $M-X$ otherwise. 
Submanifolds of a manifold $M$ are 
proper and in general position except in some obvious cases
where submanifolds are in $\partial M$.  
A surface in a three-manifold is essential if 
it is incompressible, $\partial$-incompressible, and non-boundary parallel. We denote by $(\sphere,X)$ an embedding 
of $X$ in the oriented $3$-sphere $\sphere$. 
When $X$ is a handlebody, an essential disk in $X$ is called a \emph{meridian} disk.

\subsection{Mapping class group} 
Given subpolyhedra $X_1,\cdots,X_n$ of a manifold $M$, we denote by 
\begin{equation}\label{eq:space_homeo}
\Aut{M,X_1,\dots,X_n} 
\end{equation}  
the space of self-homeomorphisms of $M$ preserving $X_i$, $i=1,\dots, n$, 
setwise, and by 
\begin{equation}\label{eq:mapping_class_group}
\MCG{M,X_1,\dots,X_n}:=\pi_0(\Aut{M,X_1,\dots ,X_n})
\end{equation}
the corresponding mapping class group.
The ``+'' subscript is added when only orientation-preserving
homeomorphisms are considered, for instance, the subspace
\[ 
\pAut{M,X_1,\dots, X_n}
\] 
of \eqref{eq:space_homeo}, and
the subgroup $\pMCG{M, X_1,\dots,X_n}$
of \eqref{eq:mapping_class_group}. 

A meridian system $\mathbf{D}$ of a handlebody ${\rm H}$ of genus $2$ 
is a triplet $\{\disk_1,\disk_2,\disk_3\}$ of disjoint, 
non-parallel, meridian disks in ${\rm H}$. 
The exterior of $\cup\dsystem:=\cup_{i=1}^{3} \disk_i$ in ${\rm H}$
consists of two $3$-balls, and $\dsystem$  
determines a trivalent spine of ${\rm H}$. 
In particular, given a handlebody-knot $\pair$ and 
a meridian system $\mathbf{D}$ of $\HK$, 
then the induced spine is either a spatial 
$\theta$-curve or handcuff graph.
Given a spatial graph $(\sphere,\Gamma)$, 
$\TSG\Gamma$ denotes the topological symmetry 
group defined in \cite{Sim:86}, which is the image of $\Sym \Gamma$
in $\MCG\Gamma$. For instance, if $\Gamma$
is a handcuff graph, then $\TSG\Gamma$ is a subgroup of the dihedral
group $\mathsf{D}_4\simeq \MCG\Gamma$.

The next two lemmas follow from the Alexander trick and \cite{Ham:66}, \cite[Section $2$]{Hat:76}, 
\cite[Theorem $1$]{Hat:99}
(see also \cite[Section $2$]{ChoKod:13},\cite[Section $2$]{Kod:15}).

\begin{lemma}\label{lm:symmetry_groups_hk_gamma}
Given a handlebody-knot $\pair$, let $\dsystem$ be a meridian system of $\HK$, and $\Gamma$ 
the induced spatial graph. 
Then 
\begin{itemize}
\item the natural homomorphisms
\[\Sym{\HK,\cup\dsystem}\rightarrow \Sym\HK, \quad
\pSym{\HK,\cup\dsystem}\rightarrow \pSym\HK\]
are injective;
\item the natural homomorphism given by the Alexander trick
\[\Sym{\HK,\cup\dsystem}\rightarrow \Sym\Gamma, \quad
\pSym{\HK,\cup\dsystem}\rightarrow\pSym\Gamma\]
are isomorphisms.
\end{itemize}
\end{lemma}

\begin{lemma}\label{lm:A_preserving_homo}
Given a handlebody-knot $\pair$ and an essential annulus $A$ in $\Compl\HK$, the natural homomorphisms
\[\Sym{\HK,A}\rightarrow \Sym\HK, \quad
\pSym{\HK,A}\rightarrow \pSym\HK\]
are injective.
\end{lemma}  
 
As a direct consequence of 
\cite[Theorems $2.5$ and $3.2$]{ChoKod:13}, we have the following Lemma.
\begin{lemma}\label{lm:spine_atoro_irre_hk}
If $\Gamma$ is a handcuff spine of an irreducible atoroidal handlebody-knot, then $\Sym\Gamma\simeq \TSG\Gamma<\mathsf{D}_4$.
\end{lemma}

\cout{
Recall that \cite[Theorem $3.3$]{KodOzaGor:15} classifies 
annuli in a handlebody-knot exterior
$\Compl\HK$ into four types based on the topology of their boundary
in relation to $\partial \HK$.  
In particular, an annulus $A$ is of type $3$
if no component of $\partial A$ bounds a disk in 
$V$ and there exists a disk $D$ either in $\Compl\HK$
or in $\HK$ disjoint from $A$. 
If $D\subset \Compl\HK$, $A$ is called a type $3$-$1$ annulus.
If $D\subset \HK$, and $\partial A$
is parallel in $\partial \HK$, 
then $A$ is called a type $3$-$2$ annulus.
If $D\subset \HK$ and $\partial A$ is not parallel
in $\partial \HK$, $A$ is called a type $3$-$3$ annulus.
} 
 
\subsection{Torus knot symmetry}\label{subsec:torus_knot}
We detour here to show that the symmetry group of 
a torus knot can be computed via 
the uniqueness of essential annuli in its exterior.
The underlying idea reappears
in the proof of Theorem \ref{teo:intro:symmetry};
on the other hand, proving the uniqueness of essential
annuli in a handlebody-knot exterior often requires more effort.

Let $\pairKpq$ be the torus knot
given by $(\frac{1}{\sqrt{2}}e^{2\pi i p t},\frac{1}{\sqrt{2}} e^{2\pi i q t})\subset\sphere\subset\mathbb{C}^2$,
$t\in \mathbb{R}$, where the integers $p,q$ are relatively prime with 
$\vert p\vert, q>1$.

\begin{theorem}[\cite{Sch:29}]
$\Sym{\Kpq}\simeq \pSym{\Kpq}\simeq \mathbb{Z}_2$.
\end{theorem} 
\begin{proof}
Let $A$ be an essential annulus in $\ComplKpq$,
and $l$ an essential loop in $A$.
$A$ cuts $\ComplKpq$
into two solid tori $V,W$.
It may be assumed that $l$
has a slope of $\frac{p}{q}$ (resp.\ $\frac{q}{p}$)
in $V$ (resp.\ $W$).
Orient $A$, $l$ and 
the cores $l_v, l_w$ 
of $V,W$,
respectively.
By the uniqueness of the annulus $A$,
we have the following isomorphisms:
\[\Sym\Kpq
\simeq \Sym{\rnbhdKpq,A},\quad
\pSym\Kpq
\simeq \pSym{\rnbhdKpq,A}
.\]

\centerline{\textbf{Claim: no $f\in\Aut{\sphere,\rnbhdKpq,A}$
swaps $V,W$}.}
If such an $f$ exists, then $f(l)$ is isotopic to $\pm l$ in $A$ and 
$f(l_v)$ isotopic to $\pm l_w$ in $W$.
In particular, we have
\[p=\lk{l}{l_v}=\pm \lk{l}{l_w}=\pm q,\] 
contradicting the assumption that $p,q$ are coprime and not equal to $\pm 1$.
In particular, the homomorphism   
\[r:
\pSym{\rnbhd{\Kpq},A}\rightarrow \MCG{\mathring{A}}
\]
given by restriction has its image in $\pMCG{\mathring{A}}$.  

\smallskip
\centerline{\textbf{Claim: $r$ is injective.}}
Suppose $f\vert_{\mathring{A}}$ is isotopic to the identity.
Then $f$ can be isotoped in $\Aut{\sphere,\rnbhdKpq,A}$ 
so that it restricts to the identity on $A$.
Let $D$ be an oriented meridian disk of $V$, and observe that 
$f(D)$ is also a meridian disk of $V$
since $f$ does not swap $V,W$, and therefore $f(\partial D)$ has an infinite slope in $V$.

Let $B_v\subset V,B_w\subset W$ be
the annuli cut off  
from $\partial\rnbhdKpq$
by $\partial A$, and 
$e_v,e_w$ two essential arcs in $B_v,B_w$, respectively.
Suppose $f(e_v)$ is not isotopic to $e_v$ in 
\[\pAut{B_v,\rel B_v}.\]
Then $f\vert_{B_v}$ is isotopic in $\pAut{B_v,\rel B_v}$ 
to the homeomorphism given by performing
Dehn twist along an essential loop
of $B_v$ $k\neq 0$ times.
This implies that $f(\partial D)$ has 
a slope of $\frac{kp+1}{kq}\in \mathbb{Q}$ in $V$, a contradiction.
$f(e_v)$ being isotopic to $e_v$ in $B_v$, relative to $\partial B_v$,
implies that $f$ can be isotoped in 
\[\pAut{\sphere,\rnbhdKpq,A}\]
so that $f\vert_{\partial V}=\id$. Applying the
same argument to $e_w$, 
we may assume $f\vert_{\partial W}=\id$.
Then, applying the Alexander trick twice, one can further isotope $f$
in $\pAut{\sphere,\rnbhdKpq,A}$ such that
$f\vert_V=\id, f\vert_W=\id$, and hence $f$
is isotopic to the identity in $\pAut{\sphere,\rnbhdKpq,A}$.
This proves the claim.

Now, consider the complex conjugation
\begin{align*}
J:\mathbb{C}^2&\rightarrow \mathbb{C}^2\\
(z_1,z_2)&\mapsto (\bar{z}_1,\bar{z}_2),
\end{align*}
and observe that $J$ 
induces a non-trivial element $g_J$ in $\Sym{\rnbhd{\Kpq},A}$ 
since $g_J$ is sent to the generator of $\pMCG{\mathring{A}}$,
where the essential annulus $A$
and the regular neighborhood $\rnbhd{\Kpq}$ are identified with 
\begin{multline*}
A_\epsilon:=\bigg\{
(\frac{1}{\sqrt{2}} e^{2\pi i p (t+s+\frac{1}{2q})},\frac{1}{\sqrt{2}} e^{2\pi i q  t})\mid 
t\in \mathbb{R}, \vert s \vert \leq \epsilon 
\bigg\} \quad\text{and}\\
\mathfrak{N}_\epsilon^{\delta_0,\delta_1}:=\bigg\{
( \sqrt{1-r} e^{2\pi i p   (t+u)}, \sqrt{r} e^{2\pi i q  t})\mid 
t\in \mathbb{R}, \vert u \vert \leq \frac{1}{2pq}-\epsilon,  \delta_0 \leq r\leq \delta_1 
\bigg\},
 \text{ respectively,}
\end{multline*}
for some $0< \epsilon <\frac{1}{2pq}$ and $0<\delta_0<\frac{1}{2}<\delta_1<1$.
This proves $\pSym\Kpq\simeq \mathbb{Z}_2$.
 
Suppose $\Sym\Kpq$ is not isomorphic to $\pSym\Kpq$.
Then there exists an orientation-reversing homeomorphism 
$f\in\Aut{\sphere,\rnbhdKpq,A}$
such that $f\vert_A$ reverses the orientation of $A$
but does not swap the two components of $\partial A$.
This implies that $f(l)$ is isotopic to $-l$ in $A$,
and since components of $\partial A$ are not swapped, 
$f(\partial D)$ is isotopic to $\partial D$ as an oriented loop,
and therefore
$f(l_v)$ is isotopic to $-l_v$. This implies 
\[p=\lk{l}{l_v}=-\lk{-l}{-l_v}=-p\]  
contradicting that $p\neq 0$. This completes the proof.
\end{proof}

\subsection{Annuli in a handlebody}
We collect here some facts about
annuli in a handlebody $V$ of genus $2$.
Let $A$ be 
an annulus in $V$, and note first the following lemma, a corollary of \cite[Lemma $9$]{BonOta:83}, \cite[Lemma $2.4$]{HayShi:01}.
\begin{lemma}\label{lm:incompressible_annulus}
If $A$ is incompressible in $V$, then $A$ is $\partial$-compressible. 
\end{lemma}

\begin{corollary}\label{cor:annuli_in_handlebody}
If $A$ is incompressible and separating, 
then $A$ divides $V$ into a solid torus $W$ 
and a handlebody $U$ of genus $2$.  
\end{corollary}
\begin{proof}
Suppose $A$ is boundary-parallel, and 
let $W$ be the solid torus through which
$A$ is parallel to an annulus in $\partial V$. Then 
$U:=V-\mathring{W}$, homeomorphic to
$V$, is a handlebody of genus $2$.
Suppose $A$ is non-boundary-parallel. 
Since $A$ is incompressible, there exists 
a $\partial$-compressing disk $D$ of $A$ by Lemma \ref{lm:incompressible_annulus}.
The boundary 
of a regular neighborhood $\rnbhd{D\cup A}$
of $D\cup A$ consists of a disk $D'$ and 
an annulus $A'$ parallel to $A$.
As $A$ is non-boundary-parallel and separating,
$D'$ is a separating essential disk in $V$,
which cuts $V$ into two solid tori
$V_1$, $V_2$, one of which, say $V_1$, contains $A$.
On the other hand, every incompressible annulus in a solid torus
is boundary-parallel, so $A$ 
cuts $V_1$ into two solid tori $W$, $U_1$,
through one of which, say $U_1$, $A$ is parallel
to an annulus in $\partial V_1$. $A$ being non-boundary-parallel 
also implies $D'\subset \partial U_1\cap \partial V_1$.
Thus, the closures of the two components of $V-A$ 
are $W$ and $U:=V_2\cup U_1$, respectively;    
the latter is a genus $2$ handlebody. 
\end{proof}

The next two corollaries follow readily from the proof of Corollary \ref{cor:annuli_in_handlebody}.
\begin{corollary}\label{cor:separating_disk}
If $A$ is incompressible, separating, and non-boundary parallel, 
then there exists a separating essential disk 
disjoint from $A$.    
\end{corollary}

\begin{corollary}\label{cor:incompressible_non_b_parallel_primitive_element}
Suppose $A$ is incompressible, separating and non-boundary-parallel,
and $W\subset V$ is the solid torus cut off by $A$.
Then the image of a generator of $\pi_1(W)$ under  
the homomorphism $\pi_1(W)\rightarrow \pi_1(V)$
induced by the inclusion
is a primitive element, up to conjugation.
\end{corollary}

\subsection{Type $3$-$3$ annulus}\label{subsec:typethreethree}
Given a handlebody-knot $\pair$ and an annulus $A\subset\ComplHK$,
we denote by $l_1,l_2$ the components of $\partial A$. 
\begin{definition}\label{def:type_three_three}
An annulus $A\subset\Compl\HK$
is of type $3$-$3$ if $l_1,l_2$ are not parallel in $\partial\HK$,
and do not bound disks in $\HK$,
and there exists a meridian disk $\disk\subset\HK$
disjoint from $A$.
\end{definition}
Henceforth $A$ denotes a type $3$-$3$ annulus unless otherwise specified.
\begin{lemma}\label{lm:separating_disks}
If $\disk\subset \HK$ is a meridian disk disjoint from $l_1\cup l_2$,
then $\disk$ is separating, and each component of $\HK-\openrnbhd{\disk}$
meets $A$. Furthermore, any two such disks are isotopic in $\HK-\partial A$.
\end{lemma}
\begin{proof}
Suppose $\disk$ is non-separating.
Then $W:=\HK-\openrnbhd{\disk}$ is a solid torus.
Since $l_1,l_2$ do not bound disks in $\HK$.
$l_1,l_2$ are parallel, essential, non-meridional loops
in $W$. $l_1\cup l_2$ cuts $\partial W$ into two annuli,
each of which meets $\rnbhd{\disk}$
since $l_1,l_2$ are not parallel in $\partial \HK$.
This, however, contradicts that $A\cap \HK=\emptyset$
since $A$ necessarily separates $\Compl W$ into two components.
Therefore $\HK-\openrnbhd{\disk}$ consists of two slid tori,
and both
meet $\partial A$ since $l_1,l_2$ are not parallel in $\partial \HK$. 
The second assertion follows from \cite[Lemma $2.3$]{LeeLee:12},
given the properties of $l_1,l_2$ just proved.
\end{proof}

\begin{remark}
By Lemma \ref{lm:separating_disks},
it is not difficult to see that Definition \ref{def:type_three_three}
is equivalent to the definition in \cite[Section $3$]{KodOzaGor:15}
(see also \cite[Theorem $3.3$]{KodOzaGor:15}).
\end{remark}

Let $\disk_A\subset \HK$ be a meridian disk disjoint from $A$,
and denote by $W_1,W_2$
the solid torus components of $\HK-\openrnbhd{\disk_A}$
with $l_i\subset W_i$, $i=1,2$,
and by $\disk_i$ a meridian disk of $W_i$, $i=1,2$.
By Lemma \ref{lm:separating_disks}, 
the meridian system 
$\dsystem_A:=\{\disk_A,\disk_1,\disk_2\}$
is determined, up to isotopy, by $A$, 
and $\dsystem_A$ induces a spatial handcuff graph
$\Gamma_A$.  
If in addition $A$ is unique, up to isotopy, 
then by Lemma \ref{lm:A_preserving_homo},
\[\gSym{\HK,A}\rightarrow \gSym\HK\]
is an isomorphism; moreover, since 
every $f\in\Aut{\sphere,\HK}$ can be isotoped
to one that preserves $\cup\dsystem_A$,    
\[\gSym{\HK,\cup\dsystem_A}\rightarrow \gSym\HK\]
is also an isomorphism by Lemma \ref{lm:symmetry_groups_hk_gamma}. 
As a result, 
we have the following corollary of Lemma \ref{lm:spine_atoro_irre_hk}

\begin{corollary}\label{cor:isom_mapping_class_groups}
If $\pair$ is irreducible, atoroidal, and $A\subset\ComplHK$ is unique, up to isotopy, then
\[\gSym{\HK,A}\simeq \gSym\HK\simeq \gSym{\HK,\cup\dsystem_A}\leq\mathsf{D}_4.\]  
\end{corollary}

A finer upper bound than $\mathsf{D}_4$
is given in Section \ref{sec:symmetry}.

\begin{definition}[\textbf{Slope Pair}]
The slope pair of $A$ is an unordered pair $\{r_1,r_2\}$ 
with $r_i$ the slope of $l_i\subset W_i$ ,
$i=1,2$.
\end{definition}
By Lemma \ref{lm:separating_disjoint},
the slope pair of $A$ is independent of the choice of $\disk_A$.
 
\begin{lemma}\label{lm:classification_slope}
If $\{r_1, r_2\}$ is the slope pair of $A$,
then either $\{r_1, r_2\}=\{\frac{p}{q},\frac{q}{p}\}$ with $pq\neq 0$,  
or  $\{r_1, r_2\}=\{\frac{p}{q},pq\}$ with $q\neq 0$,
where $p,q\in \mathbb{Z}$.
\end{lemma}
\begin{proof}
If both $r_1,r_2$ are not integers,
then $M:=W_1\cup \rnbhd{A}\cup W_2$ is a Seifert fiber space
in $\sphere$.
In particular, $M$ is the exterior of a $(p,q)$-torus knot,
and $W_1\cup W_2$ is a regular neighborhood of a Hopf link.
Hence, we have $\{r_1,r_2\}=\{\frac{p}{q},\frac{q}{p}\}$.

Suppose one of $r_1,r_2$ is integral, say $r_2$,
and let $r_1=\frac{p}{q}$. 
Then $W_1\cup\rnbhd{A}\cup W_2$ is a solid torus,
and $l_2$ has a slope of $\frac{p}{q}$ in 
the solid torus $W_1':=W_1\cup\rnbhd{A}$.
On the other hand, since $l_2\subset W_2$ has an integral slope of $r_2$, 
$r_2$ can be computed by the linking number 
of $l_2$ and any essential loop $\alpha$ 
in the annulus $\partial W_2\cap \rnbhd{A}$ 
disjoint from $l_2$. Now, $\alpha\subset W_1'$ 
also has a slope of $\frac{p}{q}$, so $r_2=pq$.
\end{proof} 

Let $\HK_A:=\rnbhd{A}\cup\HK$.
Then the following can be derived from the preceding proof.
\begin{corollary}
$\HK_A$ is a handlebody if and only if the slope pair of $A$
is $\{\frac{p}{q},pq\}$, $q\neq 0$, $p,q\in\mathbb{Z}$.
\end{corollary}

\begin{definition}[\textbf{Boundary Slope}]
$A$ is said to have a boundary slope of $p$ if $q=1$, namely, the slope pair being of the form $(p,p)$.
\end{definition}
The paper focuses primarily on the case where
$A$ has a \emph{non-trivial} boundary slope, that is, $q=1,p\neq 0$.


\section{Disks}\label{sec:disks}

Throughout the section, $A\subset\Compl\HK$ is a type $3$-$3$ 
annulus with a non-trivial boundary slope of $p$.
Unless otherwise specified, $A$ is assumed to be oriented, and
the components $l_1,l_2$ of $\partial A$ are oriented
so that $\partial A=l_1\cup -l_2$. 
Let $A_+,A_-$ be the components of $\rnbhd{A}\cap\partial\HK_A$
with the normal of $A$ in $\rnbhd{A}$ pointing toward $A_+$; 
$A_\pm$ are annuli parallel to $A$ in $\ComplHK$.
Let $l_\pm$ be essential loops in $A_\pm$, respectively, and 
orient $l_\pm$ so that they represent the same homology class
as $l_1, l_2$ in $H_1(\rnbhd{A})$.
Note that by the definition of
type $3$-$3$ annulus,
$l_\pm$ are non-separating, and hence essential,
loops in $\partial \HK_A$.

\subsection{Intrinsic disks and basis}\label{subsec:interinsic_disks}
Recall that, by Lemma \ref{lm:separating_disks} 
there is a unique meridian disk 
$\disk_A\subset\HK$ separating $l_1,l_2$; $\disk_A$ induces 
a non-separating disk $D_A \subset \HK_A$. 
Denote by $W$ the complement
$\HK_A-\openrnbhd{D_A}$, and by $D_A^\pm\subset \partial W$
the disk components of $\partial W\cap \rnbhd{D_A}$.

A meridian disk $D\subset \HK_A$ associated to $D_A$ 
is a non-separating disk disjoint from and non-parallel to $D_A$.
Particularly, $D$ can be viewed as a meridian disk of $W$.
Orient $D$ such that
$\IN{\partial D}{l_+}=1=\IN{\partial D}{l_-}$,
where 
\[
\In:H_1(\partial \HK_A)
\times H_1(\partial \HK_A)\rightarrow\mathbb{Z}
\]
is the intersection form with
the orientation of $\partial \HK_A$
given by the induced orientation of $\HK_A\subset\sphere$. 
A disk system   
associated to $D_A$ is a pair   
$\{D_1,D_2\}$ of disjoint, non-parallel 
meridian disks associated to $D_A$.
We remark that, by the definition, $D_1,D_2$ separate 
$D_A^\pm$ in $\partial W$, 
and induce a 
basis $\{[\partial D_1],[\partial D_2]\}$
of $H_1(\ComplHKA)$.

\begin{figure}[b]
\begin{subfigure}{0.3\textwidth}
\centering
\def\svgwidth{.93\columnwidth}
\begingroup%
  \makeatletter%
  \providecommand\color[2][]{%
    \errmessage{(Inkscape) Color is used for the text in Inkscape, but the package 'color.sty' is not loaded}%
    \renewcommand\color[2][]{}%
  }%
  \providecommand\transparent[1]{%
    \errmessage{(Inkscape) Transparency is used (non-zero) for the text in Inkscape, but the package 'transparent.sty' is not loaded}%
    \renewcommand\transparent[1]{}%
  }%
  \providecommand\rotatebox[2]{#2}%
  \newcommand*\fsize{\dimexpr\f@size pt\relax}%
  \newcommand*\lineheight[1]{\fontsize{\fsize}{#1\fsize}\selectfont}%
  \ifx\svgwidth\undefined%
    \setlength{\unitlength}{850.39370079bp}%
    \ifx\svgscale\undefined%
      \relax%
    \else%
      \setlength{\unitlength}{\unitlength * \real{\svgscale}}%
    \fi%
  \else%
    \setlength{\unitlength}{\svgwidth}%
  \fi%
  \global\let\svgwidth\undefined%
  \global\let\svgscale\undefined%
  \makeatother%
  \begin{picture}(1,1)%
    \lineheight{1}%
    \setlength\tabcolsep{0pt}%
    \put(0,0){\includegraphics[width=\unitlength,page=1]{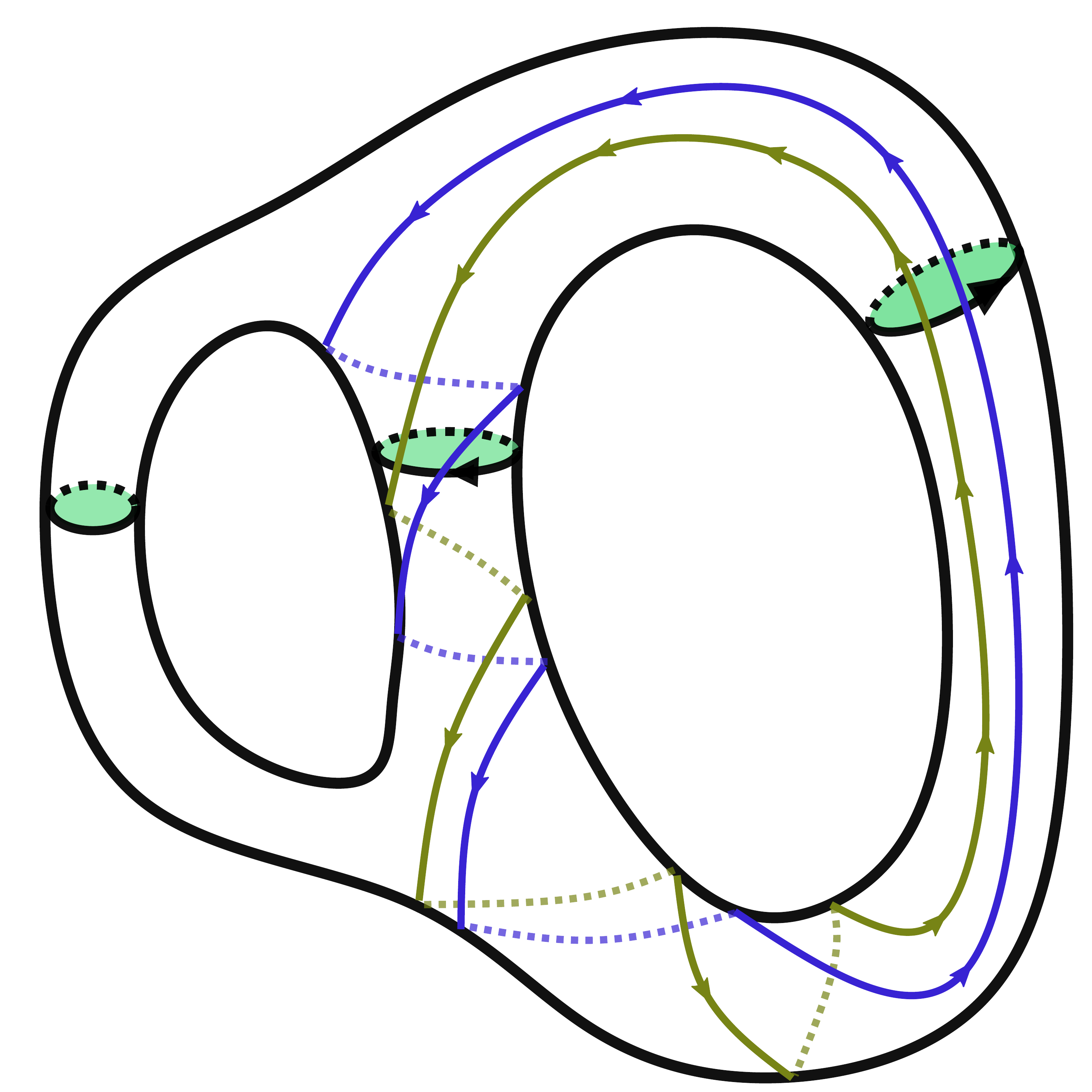}}%
    \put(0.14091874,0.53671444){\color[rgb]{0,0,0}\makebox(0,0)[lt]{\lineheight{1.25}\smash{\begin{tabular}[t]{l}{\footnotesize $D_A$}\end{tabular}}}}%
    \put(0.49041409,0.58601897){\color[rgb]{0,0,0}\makebox(0,0)[lt]{\lineheight{1.25}\smash{\begin{tabular}[t]{l}{\footnotesize $D_1$}\end{tabular}}}}%
    \put(0.70302314,0.64192749){\color[rgb]{0,0,0}\makebox(0,0)[lt]{\lineheight{1.25}\smash{\begin{tabular}[t]{l}{\footnotesize $D_2$}\end{tabular}}}}%
    \put(0.26965698,0.76263903){\color[rgb]{0,0,0}\makebox(0,0)[lt]{\lineheight{1.25}\smash{\begin{tabular}[t]{l}{\footnotesize $l_+$}\end{tabular}}}}%
    \put(0.47431548,0.7614435){\color[rgb]{0,0,0}\makebox(0,0)[lt]{\lineheight{1.25}\smash{\begin{tabular}[t]{l}{\footnotesize $l_-$}\end{tabular}}}}%
  \end{picture}%
\endgroup%

\caption{Disk system.}
\label{fig:disk_system}
\end{subfigure}
\begin{subfigure}{0.33\textwidth}
\centering
\def\svgwidth{.93\columnwidth}
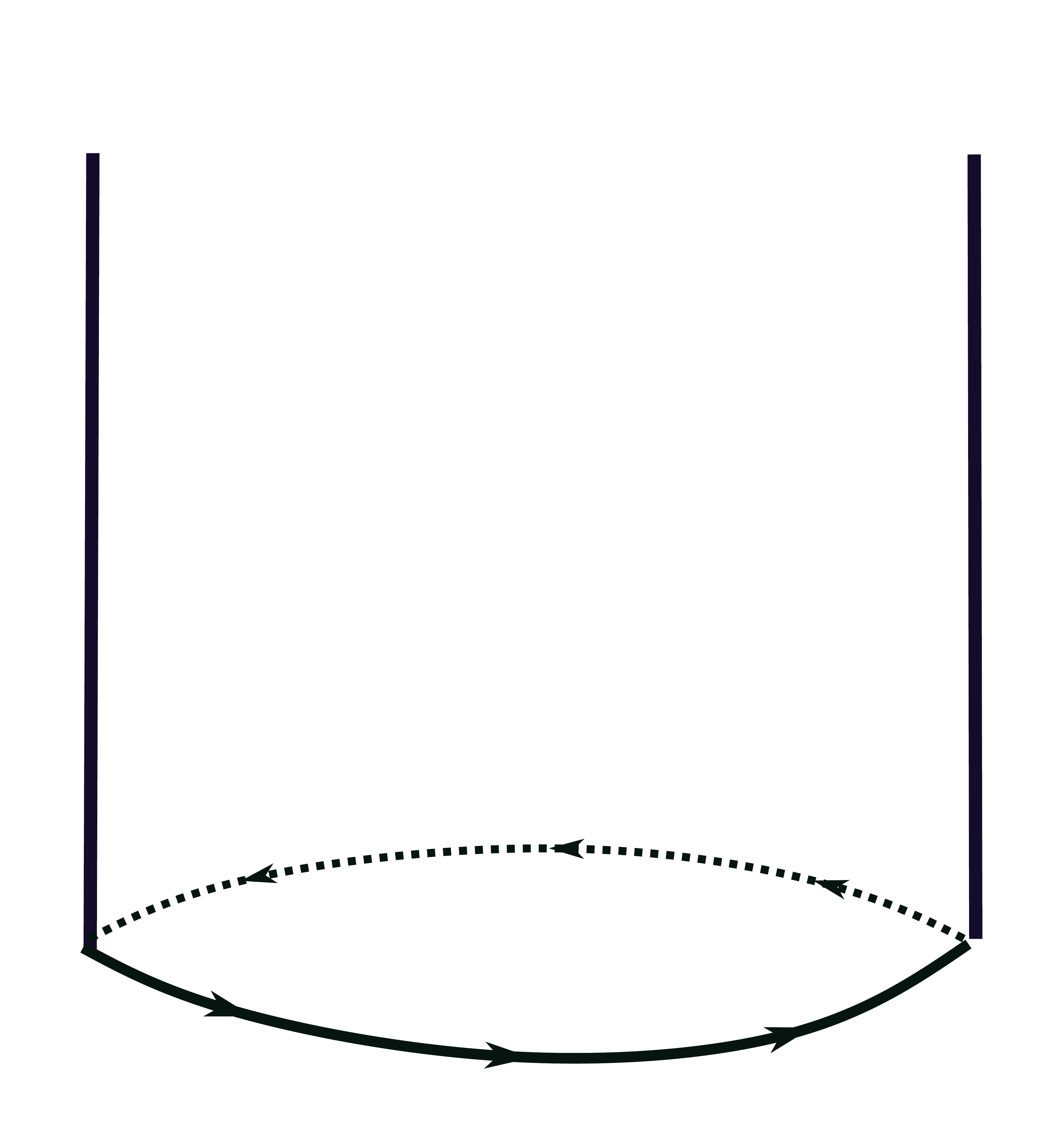
\caption{A disk $D$ with $[\partial D]=x'$.}
\label{fig:meridian_disk_n}
\end{subfigure}
\begin{subfigure}{0.33\textwidth}
\centering
\def\svgwidth{1.0\columnwidth}
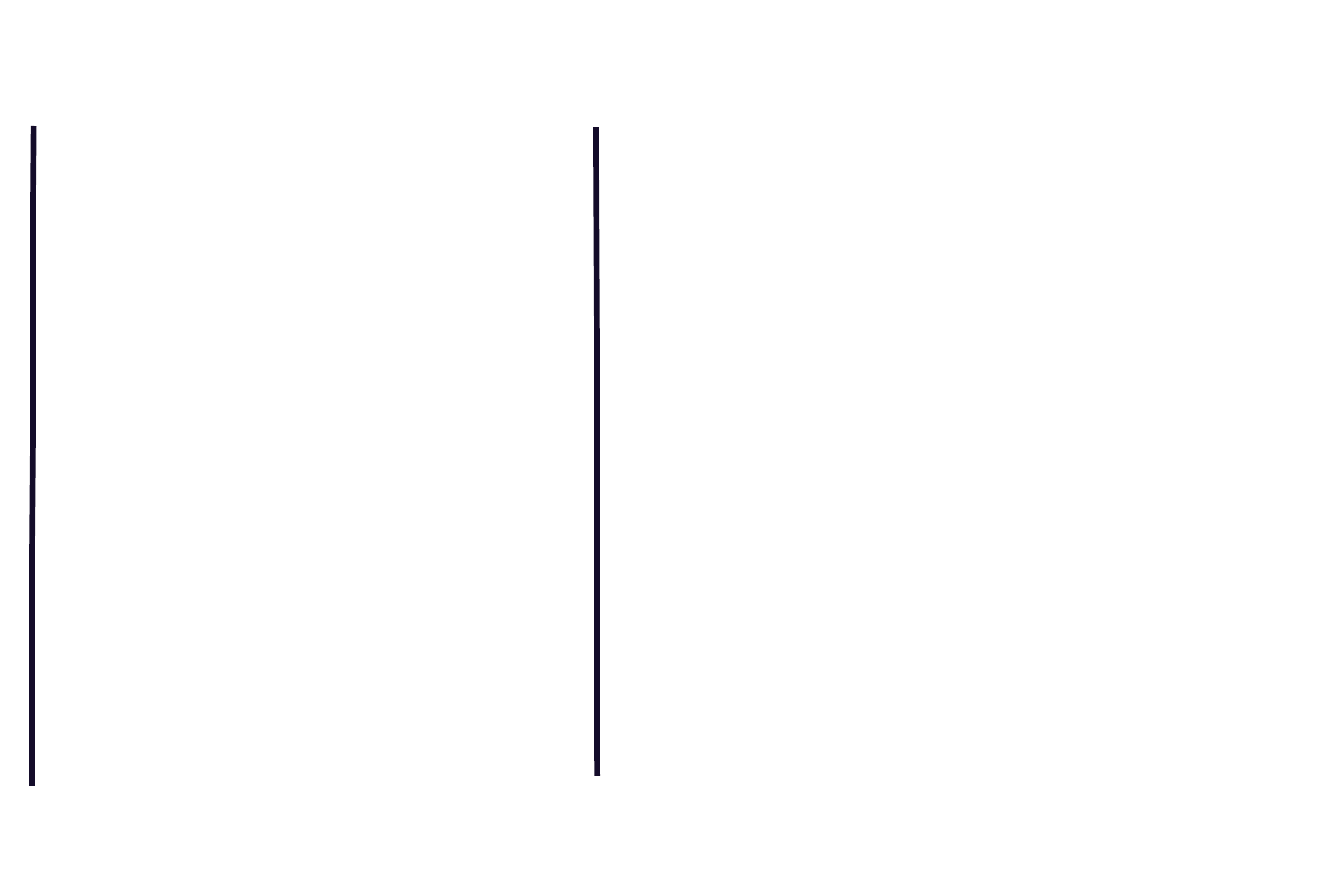
\caption{Disk $D_2'$.}
\label{fig:Dtwo_prime}
\end{subfigure}
\caption{Disk system and meridional basis.}
\end{figure} 

\begin{definition}\label{def:meridional_basis}
A basis $\{a_1,a_2\}$ of $H_1(\ComplHKA)$
is \emph{meridional} if 
it is induced by a disk system 
associated to $D_A$.
\end{definition}

The existence of a disk system  
is easy to check, It implies
the existence of a meridional basis of $H_1(\ComplHKA)$.

\begin{lemma}\label{lm:l_pm_via_meridional_basis}
Let $\{a_1,a_2\}$ be a meridional basis  
of $H_1(\ComplHKA)$.
Then $[l_+]-[l_-]\in H_1(\ComplHKA)$
is $(1-1)$ or $(-1,1)$ in terms of $\{a_1,a_2\}$. 
In addition, if   
$[l_+]=(q_1,q_2)$,
then $q_1+q_2=p$.
\end{lemma}
\begin{proof}
Note first that $l_\pm$
are parallel loops in $\partial W$
with each component of $\partial W-(l_+\cup l_-)$
containing one of 
$D_A^+,D_A^-$. 
In particular, one can orient $D_A$ so that 
$[l_+]=[l_-]+ [\partial D_A]\in H_1(\ComplHKA)$.
The first assertion then follows from the fact that
$[\partial D_A]$ is either $a_1-a_2$ or $a_2-a_1$.
 
Now, consider the homomorphism 
\[H_1(\ComplHKA)\rightarrow H_1(\Compl{W}),\]
and observe that it
sends both $a_1,a_2$ to a generator $a$ of $H_1(\Compl{W})$. 
Since $A$ has a boundary slope of $p$,  
$[l_+]=p a\in H_1(\Compl{W})$, and thus the second assertion.   
\end{proof}

For the next two lemmas, we fix a 
meridional basis $\{a_1,a_2\}$ of 
$H_1(\ComplHKA)$ induced by a disk system  
$\{D_1,D_2\}$.

\begin{lemma}\label{lm:associated_meridian_disk}
If $D$ is a meridian disk associated to $D_A$, then
\[[\partial D]=a_1+n(a_2-a_1),\quad\text{for some $n\in\mathbb{Z}$}.\]  
Conversely, if $a_1'=a_1+n(a_2-a_1)$, for some $n\in\mathbb{Z}$, 
then there exists 
a meridian disk $D$ associated to $D_A$ such that
$a_1'=[\partial D]$.
\end{lemma}   
\begin{proof}
Isotope $D$ in $\HK_A$ 
so that $m:=\# D\cap (D_1\cup D_2)$ is minimized.
We prove by induction on $m$.
If $m=0$, then $D$ is parallel to either $D_1$ or $D_2$ in $\HK_A$,
and thus $n=0$ or $1$.

Suppose the statement holds for any meridian disk associated to
$D_A$ with $0\leq m<k$, and
$D$ is a meridian disk associated to $D_A$
with $m=k$. 
Let $D^\alpha\subset D$ 
be an outermost disk cut off by an outermost arc 
$\alpha\subset D\cap (D_1\cup D_2)$ in $D$. It may be 
assumed that $\alpha\subset D_1$. 
$\alpha$ cuts off a disk $D_1^\alpha$ from $D_1$
such that $\tilde{D}^\alpha:=D^\alpha\cup D_1^\alpha\subset W$ 
is inessential in $W$.
By the minimality, $\partial \tilde{D}^\alpha$
bounds a disk $E^\alpha\subset \partial W$
that contains one or both of $D_A^\pm$.
Orient $D_A$ so that 
$[\partial D_A]$ is $a_2-a_1$. 

On the other hand, the union $(D-\mathring{D^\alpha}) \cup D_1^\alpha$
induces a disk $D'$ with $\# D'\cap (D_1\cup D_2)<k$.
Therefore by induction, $[\partial D']=a_1+n'(a_2-a_1)$, for some $n'\in \mathbb{Z}$.
The assertion then follows by observing that 
$[\partial D]=[\partial D']+t(a_2-a_1)$, where $t=0$ or $\pm 1$,
depending on $E^\alpha$ contains both or only one of $D_A^\pm$,  
respectively.   

To see the second assertion, note that, for 
$n=0,1$, one can take $D=D_1$, $D_2$, respectively.
Suppose $n\geq 2$. Then the meridian disk $D$ can be 
constructed as in Fig.\ \ref{fig:meridian_disk_n},
where the boundary of $D$ is depicted. 
It is not difficult to check that $[\partial D]=a_1'$.
\end{proof}
\begin{lemma}\label{lm:disk_system}
If $\{D_1',D_2'\}$ is a disk system associated to $D_A$, 
then
\begin{align*}
[\partial D_1']&=a_1+n(a_2-a_1)\\
[\partial D_2']&=a_1+(n\pm 1)(a_2-a_1),
\end{align*} 
for some $n\in\mathbb{Z}$.
Conversely, if 
\begin{align*}
a_1'&=a_1+n(a_2-a_1)\\
a_2'&=a_1+(n\pm 1)(a_2-a_1),\quad\text{for some $n\in\mathbb{Z}$,}
\end{align*}
then the basis $\{a_1',a_2'\}$ is meridional.
\end{lemma}
\begin{proof}
By Lemma \ref{lm:associated_meridian_disk}, 
$[\partial D_1']=a_1+n(a_2-a_1)$ for some $n\in \mathbb{Z}$.
Since $\partial D_1',\partial D_2'$ cobound
a cylinder containing exactly one of $D_A^\pm$.
$[\partial D_2']=[\partial D_1']\pm (a_1-a_2)$.

For the opposite direction, by Lemma \ref{lm:associated_meridian_disk},
there is a meridian disk $D_1'$ associated to $D_A$
with $[\partial D_1']=a_1'$. 
Properly choosing a partition of $\partial W-\partial D_1'$
that separates $\{D_A^+,D_A^-\}$ as in Fig.\ \ref{fig:Dtwo_prime},
we obtain another meridian disk $D_2'$ associated 
to $D_A$ with $\{a_1',a_2'\}$ induced by $\{D_1',D_2'\}$.   
\end{proof}  
\begin{definition}\label{def:normalized_basis}
A normalized basis $\{a_1,a_2\}$ of $H_1(\ComplHKA)$
is a meridional basis in terms of which
$[l_+]=(p_1,p_2)$ with either $0<p_1\leq p$ or $p<p_2\leq 0$
and $[l_+]-[l_-]=(1,-1)$.
\end{definition}
\begin{lemma}\label{lm:existence_uniqueness_normalized_basis}
There exists a normalized basis $\{a_1,a_2\}$ of $H_1(\ComplHKA)$.
Furthermore, given two normalized basis 
$\{a_1,a_2\}$, $\{a_1',a_2'\}$ of $H_1(\ComplHKA)$,
if $[l_+]=(p_1,p_2)$ and $[l_+]=(p_1',p_2')$ in terms of $\{a_1,a_2\}$ and $\{a_1',a_2'\}$, respectively,
then $p_i=p_i'$, $i=1,2$.
\end{lemma}
\begin{proof}
Start with a meridional basis $\{a_1',a_2'\}$ 
of $H_1(\ComplHKA)$, in terms of which 
$[l_+]=(q_1,q_2)$ and $[l_-]=(q_1-1,q_2+1)$.
%
There exists an $n\in\mathbb{Z}$
such that either $0<np+q_1\leq p$ or $p<np+q_1\leq 0$ holds.

By the second assertion of 
Lemma \ref{lm:disk_system}, there exists a 
meridional basis $\{a_1,a_2\}$ 
such that
\[
\{a_1, a_2\}=\{a_1', a_2'\}
\begin{bmatrix}
1-n&-n\\
n& n+1 
\end{bmatrix}, 
\]
In particular, if in terms of $\{a_1,a_2\}$, 
$[l_+]=(p_1,p_2)$, then $[l_-]=(p_1-1,p_2+1)$ and
\begin{align}\label{eq:coefficient_change}
p_1 &=np+q_1\\
p_2 &=-np+q_2\nonumber  
\end{align}   
because
\[
\begin{bmatrix}
p_1 \\
p_2
\end{bmatrix}
=
\begin{bmatrix}
n+1&n\\
-n&1-n
\end{bmatrix} 
\begin{bmatrix}
q_1\\
q_2
\end{bmatrix}. 
\] 
Therefore either $0<p_1\leq p$ or $p<p_1\leq 0$, and 
$\{a_1, a_2\}$ is a normalized basis. 

For the second assertion, we note that, by the first statement of 
Lemma \ref{lm:disk_system}, either
\begin{equation}\label{eq:lm:uniqueness_two_transition_matrices}
\{a_1,a_2\}=\{a_1',a_2'\}
\begin{bmatrix}
1-n&-n\\
n& n+1 
\end{bmatrix} 
\quad
\text{or}
\quad 
\{a_1,a_2\}=\{a_1',a_2'\}
\begin{bmatrix}
1-n&-n+2\\
n& n-1 
\end{bmatrix}.
\end{equation} 
The latter implies $a_1-a_2=a_2'-a_1'$ 
contradicting the definition of a normalized basis.
On the other hand, since $0< p_1,p_1'\leq p$ 
or $p<p_1,p_1'\leq 0$, 
the only possible $n$ in the first matrix of \eqref{eq:lm:uniqueness_two_transition_matrices}  
is $0$, and hence the assertion.
\end{proof}
\begin{corollary}\label{cor:l_pm}
$[l_+]\neq \pm [l_-]\in H_1(\Compl{\HK_A})$, 
and neither is trivial in $H_1(\Compl{\HK_A})$, and
$[l_+],[l_-]$ form a basis of $H_1(\Compl{\HK_A})$
if and only if $\vert p\vert= 1$.
\end{corollary}
\begin{proof}
Let $[l_+]=(p_1,p_2)$, $[l_-]=(p_1-1,p_2+1)$
in terms of a meridional basis. 
Then the first and second assertions 
follow from the fact that $p_1,p_2$ are integers
and $p_1+p_2=p\neq 0$, respectively.
The fact that the determinant  
\[
\begin{vmatrix}
p_1& p_2\\
p_1-1& p_2+1
\end{vmatrix}
=p
\] 
implies the third.
\end{proof}
Since changing the orientation of $l_+,l_-$
changes the orientation of a disk system,
the uniquenes part of Lemma \ref{lm:existence_uniqueness_normalized_basis} 
gives us the following invariant of the triplet $\pairandA$.
 
\begin{definition}
The slope invariant of $\pairandA$
is the pair $(p_1,p_2)$ of integers such that 
$[l_+]=(p_1,p_2)$ in terms of a normalized basis
of $H_1(\ComplHKA)$.
\end{definition}
If the orientation of $A$ is reversed, denoted by $\mirror A$, 
then $l_+,l_-$ are swapped and their orientation reversed.
In particular, $(a_1,a_2)$
is a normalized basis of $H_1(\ComplHKA)$
if and only if $(-a_2,-a_1)$ is a normalized basis
of $H_1(\Compl{\HK_{\mirror A}})$; thus
the slope invariant $\pairandmA$
is $(p_2+1,p_1-1)$.
\begin{corollary}\label{cor:reversing_orientation_A}
Let $f$ be a homeomorphism in $\Aut{\sphere,\HK,A}$
such that $f\vert_A$ is orientation-reversing.
Then the slope invariant of $\pairandA$
is $(\frac{p+1}{2},\frac{p-1}{2})$.
In particular, no such
homeomorphisms exist when $p$ is even.
\end{corollary}
\begin{proof}
Observe first that if  
$(a_1,a_2)$ is a normalized basis of $H_1(\ComplHKA)$,
then 
\[\big(f_\ast(a_1),f_\ast(a_2)\big)\]
is a normalized basis of $H_1(\Compl{\HK_{mA}})$.
Secondly, if 
$(p_1,p_2)$ is the slope invariant of $\pairandA$,
then the slope invariant of $\pairandmA$ is $(p_2+1,p_1-1)$.
%
As a result, we have
$(p_1,p_2)=(p_2+1,p_1-1)$, which together with 
$p_1+p_2=p$ implies $(p_1,p_2)=(\frac{p+1}{2},\frac{p-1}{2})$.
\end{proof}

\begin{definition}
An orientation of $A$ is called a \emph{preferred} orientation
if the slope invariant $\pairandA$
is $(p_1,p_2)$ with $p_1>p_2$.
\end{definition}

The existence of a preferred orientation 
can be easily verified since  
$p_1>p_2$ if and only if $p_1-1 \geq p_2+1$
when $(p_1,p_2)\neq (\frac{p+1}{2},\frac{p-1}{2})$.
In the event $(p_1,p_2)= (\frac{p+1}{2},\frac{p-1}{2})$,
both $A,\mirror A$ are preferred orientations.
The observation above allows us to define 
the slope type of an unoriented $A$.
\begin{definition}
The slope type of an unoriented $A\subset\Compl\HK$  
is the slope invariant of $\pairandA$
with $A$ endowed with a preferred orientation. 
\end{definition}
The slope type gives a finer classification
of type $3$-$3$ annuli with a non-trivial slope,
and is used to examine handlebody-knot families 
in Section \ref{subsec:examples_uniqueness}.
 
%
%
%
%
%
\begin{definition}\label{def:good_excellent_basis}
A \emph{good} basis  
of $H_1(\partial\HK_A)$
is a basis $\{a_1,b_1,a_2,b_2\}$ 
such that $\{b_1,b_2\}$ generates the kernel of
\[j:H_1(\partial\HK_A)\rightarrow H_1(\ComplHKA),\]
and $\{a_1,a_2\}$ viewed as elements of $H_1(\ComplHKA)$
is a meridional basis, and the intersection form $\In$ on $H_1(\partial \HK_A)$
is realized by the matrix
\begin{equation}\label{eq:intersection_form_matrix}
\begin{bmatrix}
0&1&0&0\\
-1&0&0&0\\
0&0&0&1\\
0&0&-1&0
\end{bmatrix},
\end{equation}
in terms of $\{a_1,b_1,a_2,b_2\}$.

A good basis $\{a_1,b_1,a_2,b_2\}$ 
is excellent if $\{a_1,a_2\}\subset H_1(\ComplHKA)$
is normalized.
\end{definition}
\begin{lemma}\label{lm:existence_good_basis}
Given a disk system $\{D_1,D_2\}$,
there exists a good basis $\{a_1,b_1,a_2,b_2\}$
of $H_1(\partial \HK_A)$
with $a_i=[\partial D_i]$, $i=1,2$.
\end{lemma}
\begin{proof}
Let $a_i$, $i=1,2$, be elements in $H_1(\partial \HK_A)$
represented by $\partial D_i$, and
choose two disjoint oriented loops $l_i$, $i=1,2$,
that satisfy the following properties: $D_i\cap l_j$ is a point when $i=j$
and empty otherwise, and $\IN{\partial D_i}{l_i}=1$, $i=1,2$.
Suppose $[l_i]=t_{i1}a_1+t_{i2}a_2\in H_1(\ComplHKA)$.
Then define $b_i:=[l_i]-t_{i1}a_1-t_{i2}a_2\in H_1(\partial\HK_A)$.
In particular, $\{b_1,b_2\}$ form a basis of $\op{Ker}(j)$.

Since 
$b_1,b_2$ are in the image of $H_2(\ComplHKA,\partial \HK_A)\xrightarrow{\partial} H_1(\partial\HK_A)$, $\mathcal{I}(b_1,b_2)=0$. 
On the other hand, by the construction of $l_i$, $i=1,2$,
we have $\mathcal{I}(a_i,b_j)$ is $1$ when $i=j$ and $0$ otherwise.
These, together with $\mathcal{I}(a_1,a_2)=0$, 
implies that 
the intersection form $\mathcal{I}$ is realized by
the matrix \eqref{eq:intersection_form_matrix}
in terms of the basis $\{a_1,b_1,a_2,b_2\}$. 
\end{proof}
As a result of Lemmas \ref{lm:l_pm_via_meridional_basis}, \ref{lm:existence_uniqueness_normalized_basis}
and \ref{lm:existence_good_basis}, we have the following corollaries. 
\begin{corollary}\label{def:excellent_basis}
There exists an excellent basis of $H_1(\partial \HK_A)$.
\end{corollary}

%
\begin{corollary}\label{cor:l_pm_via_good_basis}
Let $\{a_1,b_1,a_2,b_2\}$ be a good basis  
of $H_1(\partial \HK_A)$.
Then in terms of the basis, 
$[l_+]=(q_1,1,q_2,1)\in H_1(\partial \HK_A)$ 
with $q_1+q_2=p$, and
$[l_+]-[l_-]=(1,0,-1,0)$
or $(-1,0,1,0)$. 
\end{corollary}
\begin{proof}
Since $\IN{\partial D_i}{l_+}=1$ (resp.\ $\IN{\partial D_i}{l_-}=1$), $i=1,2$,
the coefficients of $b_1,b_2$ in $[l_+]$ (resp.\ $[l_-]$) are $1$.
The rest follows readily from Lemma \ref{lm:l_pm_via_meridional_basis}. 
\end{proof}
\subsection{Extrinsic disks}
Through the subsection, we assume $\pairA$ is trivial---namely, 
$\ComplHKA$ is a handlebody.
Given an oriented disk $D\subset\Compl{\HK_A}$, denote by 
$I_D$ the pair 
\[\big(\IN{l_+}{\partial D}, \IN{l_-}{\partial D}\big)\]
of intersection numbers, where $\In$
is the intersection form on $H_1(\partial \HK_A)$.
\begin{lemma}\label{lm:trivial_intersection}
$I_D=(0,0)$ if and only if $D$ is a separating disk.
\end{lemma}
\begin{proof}
The ``if'' part is clear.
For the ``only if'' part, observe that  
$[\partial D]=(0,d_1,0,d_2)$, for some $d_1,d_2\in\mathbb{Z}$, 
in terms of a good basis of $H_1(\partial \HK_A)$. 
$I_D=(0,0)$ then implies  
\[
\begin{cases}
 p_1d_1+p_2d_2=0\\
(p_1-1)d_1+(p_2+1)d_2=0.
\end{cases}
\]
Since $p_1(p_2+1)-(p_1-1)p_2 =p_1+p_2=p\neq 0$, we have 
$d_1=d_2=0$, and therefore $\partial D\subset \partial\HK_A$ is separating. Thus, $D\subset \ComplHKA$ is a separating disk.
\end{proof}
\begin{corollary}\label{cor:non_separating_disjoint}
There exists no non-separating disk in $\Compl{\HK_A}$
disjoint from $l_+\cup l_-$.
\end{corollary}
\begin{lemma}\label{lm:plus_minus_one_parallel}
Let $D_1,D_2$ be two disjoint oriented disks in $\ComplHKA$.
Suppose $I_{D_1}=I_{D_2}=(1,-1)$ or $(-1,1)$.
Then $D_1, D_2$ are parallel.
\end{lemma}
\begin{proof}
Since $I_{D_i}=(\pm 1,\mp 1),i=1,2$, the disks
$D_1,D_2$ are essential and non-separating.
If they are not parallel,
then $\ComplHKA-\openrnbhd{D_1\cup D_2}$
is a $3$-ball, and hence, 
there exist loops $\alpha_i$, $i=1,2$, 
such that $\alpha_i\cap \partial D_j$
is a point when $i=j$ and empty otherwise,
and $\IN{\alpha_i}{\partial D_i}=1$.

The set $\{[\alpha_1],[\partial D_1],[\alpha_2],[\partial D_2]\}$
is a basis of $H_1(\partial\HK_A)$ such that
the images of $[\alpha_i]$, $i=1,2$,
under $j:H_1(\partial\HK_A)\rightarrow H_1(\ComplHKA)$
generate $H_1(\ComplHKA)$, and $[\partial D_i]$, $i=1,2$,
form a basis of $\op{Ker}(j)$.
In terms of the basis, $[l_+]=(\pm 1, c_+, \pm 1,d_+)$
and $[l_-]=(\mp 1, c_-, \mp 1,d_-)$,
for some $c_\pm,d_\pm\in\mathbb{Z}$, 
since $I_{D_1}=I_{D_2}=(\pm 1, \mp 1)$.
This implies that $[l_+]=-[l_-]$\ in $H_1(\ComplHKA)$,
contradicting Corollary \ref{cor:l_pm}.
\end{proof} 

Denote by $\cl{l_+},\cl{l_-}$ the conjugacy classes
determined by $l_+,l_-$ in $\pi_1(\ComplHKA)$, respectively.
We say $l_+$ (resp.\ $l_-$) \emph{represents} the $n$-th power
of a primitive element $x\in \pi_1(\ComplHKA)$ if there exists $a\in \cl{l_+}$
(resp.\ $b\in \cl{l_-}$) such that 
$a=x^n$ (resp.\ $b=x^n$), and say $\{l_+,l_-\}$ 
\emph{represents a basis} if there exist $a\in \cl{l_+}, b\in \cl{l_-}$
such that $\{a,b\}$ forms a basis of the free group $\pi_1(\ComplHKA)$.

\begin{lemma}\label{lm:separating_disjoint}
Suppose there exists an essential separating disk $D\subset \ComplHKA$ disjoint from $l_+\cup l_-$. Then 
\begin{enumerate} 
\item $l_1,l_2$
are either both $(m,n)$-torus knots with $mn=p$ or both trivial knots in $\sphere$;
\item 
there exists a basis $\{a,b\}$ of $\pi_1(\ComplHKA)$
such that 
$a^{\vert m\vert}\in\cl{l_+}$
and $b^n\in\cl{l_-}$
when $l_1,l_2$ are $(m,n)$-torus knots in $\sphere$, or 
$a$ is in one of 
$\cl{l_+}, \cl{l_-}$ and $b^{\absp}$ in the other
when $l_1,l_2$ are trivial knots in $\sphere$.
\end{enumerate} 
\end{lemma} 
\begin{proof}
If such a disk $D$ exists, then it
separates $\ComplHKA$ into two solid tori $W_+,W_-$. 
$l_+,l_-$ are not in the same solid torus because $[l_+]\neq \pm [l_-]\in H_1(\ComplHKA)$ by Corollary \ref{cor:l_pm}.
It may be assumed that $l_+\subset \partial W_+$ and
$l_-\subset \partial W_-$. 
Suppose the slopes of $l_+,l_-$ 
in $W_+,W_-$ are $\frac{m}{n}, \frac{m'}{n'}$, $n,n'>0$
respectively. Then $mn=m'n'=p$ since $\lk{l_1}{l_2}=p$.

Now, the longitudes of $W_+,W_-$ 
induces a basis 
of $H_1(\ComplHKA)$,
in terms of which $[l_+]=(n,0)$ and $[l_-]=(0,n')$.
On the other hand, in terms of a normalized basis,
$[l_+]=(p_1,p_2)$ and $[l_-]=(p_1-1,p_2+1)$, $p_1+p_2=p$.
Therefore, we have
\begin{equation}\label{eq:change_of_basis_determinant}
n n'=
\begin{vmatrix}
n &0\\
0&n'
\end{vmatrix}
=\pm 
\begin{vmatrix}
p_1&p_2\\
p_1-1&p_2+1
\end{vmatrix}
=\pm p.
\end{equation} 
Especially, the slope of $l_-\subset W_-$ is $\frac{n}{m}$.

Let $L_+,L_-$ be cores of $W_+,W_-$, respectively.
Then $(\sphere, L_+\cup L_-)$ is a non-simple link
since there exists an essential annulus $A^\perp\subset \rnbhd{A}
\subset \Compl{W_+\cup W_-}$ bounded by $l_+\cup l_-$.
In addition, because $W_+\cup W_-\cup \rnbhd{D}=\ComplHKA$
is a handlebody, $(\sphere, L_+\cup L_-)$ is a tunnel number one link.

Consider first the case $n=1$, and hence $m=p$,
by the classification of tunnel number one non-simple links
in \cite{EudUch:96}, $(\sphere, L_-)$ is trivial, and
$L_+$ is a $(1,p)$-curve on a regular neighborhood of 
$L_-$. This implies $l_+,l_-$ and hence $l_1,l_2$ are all trivial
knots in $\sphere$. 
Furthermore, there exists
a basis $\{a,b\}$ of $\pi_1(\ComplHKA)$ 
induced by the longitudes of $W_+,W_-$,
respectively, such that $a\in\cl{l_+}$ and $b^{\vert p\vert}\in\cl{l_-}$.
The same argument applies to the case where $n=\vert p\vert$
and $m=\pm 1$.

Suppose $n>1,\vert m\vert>1$. 
Again by the classification of 
tunnel number one non-simple links \cite{EudUch:96}, 
$(\sphere, L_+\cup L_-)$ 
is a Hopf link. Thus $l_+, l_-$ 
and therefore $l_1,l_2$ are $(m,n)$-torus knots in $\sphere$ with $mn=p$.
As in the previous case, there exists
a basis $\{a,b\}$ of $\pi_1(\ComplHKA)$ 
induced by the longitudes of $W_+,W_-$ 
such that $a^n\in \cl{l_+}$, $b^{\vert m\vert}\in \cl{l_-}$. 
\end{proof}

\begin{remark}\label{rmk:lpm_primtive_trivial_HK}
Any separating essential disk in $\ComplHKA$
disjoint from $l_+\cup l_-$ 
induces a 
non-separating disk in $\ComplHK$,
and implies that $\pair$ is reducible,
Furthermore, since $\ComplHKA\cup\rnbhd{A}=\ComplHK$ and $\ComplHKA\cap\rnbhd{A}=A_+\cup A_-$, the fundamental group $\pi_1(\ComplHK)$
is an HNN-extension of $\pi_1(\ComplHKA)$.
Thus by Lemma \ref{lm:separating_disjoint}, if $l_1,l_2$
are trivial knots, then $\pi_1(\ComplHK)$ is free, and $\pair$
is trivial.
\end{remark}

Lemma \ref{lm:separating_disjoint} implies the following 
algebro-geometric obstruction to the existence of 
an essential, separating disk in $\ComplHKA$
disjoint from $l_+\cup l_-$. 
 
\begin{corollary}\label{cor:separating_disjoint}
Suppose one of the following holds:
\begin{itemize}
\item $\vert p\vert=1$, and   
$\{l_+, l_-\}$ does not represent a basis of  $\pi_1(\ComplHKA)$.
\item $l_1,l_2$ 
are not $(m,n)$-torus knot in $\sphere$ with 
$mn=p$, and neither of $l_+,l_-$ represents
the $\vert p\vert$-th power of some primitive element of $\pi_1(\ComplHKA)$.
\item $l_1,l_2$ 
are $(m,n)$-torus knots in $\sphere$ with 
$mn=p$, and 
one of $l_+,l_-$
does not represent the $n$-th or $\vert m\vert$-th power
of any primitive element in $\pi_1(\ComplHKA)$.
\end{itemize} 
Then there is no essential separating disk $D\subset\ComplHKA$ disjoint from $l_+\cup l_-$.
\end{corollary}

\begin{corollary}\label{cor:non-separating_one_intersection}
Under the same conditions as in Corollary \ref{cor:separating_disjoint}, 
there is no disk $D\subset\ComplHKA$ that intersects $l_+\cup l_-$
at one point. 
\end{corollary}
\begin{proof}
Suppose such a disk $D$ exists;
it may be assumed that $D\cap l_+$ is a point.
Then the boundary of a regular neighborhood
of $D\cap l_+$ disjoint from $l_-$
is an essential, separating disk disjoint from $l_+\cup l_-$,
contradicting Corollary \ref{cor:separating_disjoint}.
\end{proof}

\begin{lemma}\label{lm:non-separating_two_intersections}
Suppose $p$ is odd, $\vert p\vert> 1$, and neither of
$[l_+],[l_-]$ is the $p$-th multiple of some element in $H_1(\Compl{\HK_A})$. Then there exists no essential disk $D\subset\ComplHKA$
with $I_D=(\pm k,0)$, $(0,\pm k)$, $(1,1)$ or $(-1,-1)$, 
where $k=1,2$.
\end{lemma}
\begin{proof}
In terms of an excellent basis of $H_1(\partial\HK_A)$,
we have $[l_+]=(p_1,1,p_2,1)$, $[l_-]=(p_1-1,1,p_2+1,1)$ 
and $[\partial D]=(0,d_1,0,d_2)$ with
either $0<p_1\leq p$ or $p<p_1\leq 0$ and $d_1,d_2\in\mathbb{Z}$. 
The condition that neither of $[l_+],[l_-]$ is the $p$-th multiple of some element in $H_1(\ComplHKA)$ implies that 
$1<p_1<p$ and $0<p_2<p-1$ (resp.\ 
$p+1<p_1<0$ and $p<p_2<-1$) when $p>0$ (resp. $p<0$).

Suppose $I_D=(m,n)$. Then we have the system of equations
\begin{equation}\label{eq:system_for_I_D}
\begin{cases}
p_1d_1+p_2d_2=m\\
(p_1-1)d_1+(p_2+1)d_2=n\nonumber.
\end{cases}
\end{equation}  
Consider first the case $(m.n)=(\pm k,0)$. Then by \eqref{eq:system_for_I_D},
\begin{equation}\label{eq:k_0_disk}
(d_1, d_2)=(\pm\frac{k(p_2+1)}{p}, \pm\frac{k(p_1-1)}{p}).
\end{equation} 
The constraints 
$1<p_1<p$ and $0<p_2<p-1$ (resp.\ $p+1<p_1<0$ and $p<p_2<-1$), 
imply $\vert p_1-1\vert, \vert p_2+1\vert$ 
are not zero or $\vert p \vert$.
Moreover, since $p$ is odd and $p_1-1,p_2+1$
are of the same sign, one of $\vert p_1-1\vert , \vert p_2+1\vert$
is smaller than $\vert\frac{p}{2}\vert$.
Therefore \eqref{eq:k_0_disk} 
is not an integral solution when $1\leq k\leq 2$. 

Similarly by \eqref{eq:system_for_I_D}, 
if $I_D=(0,k)$ or $(\pm 1,\pm1)$,
then $(d_1,d_2)=(\frac{kp_2}{p},\frac{kp_1}{p})$
or $(\pm\frac{1}{p},\pm\frac{1}{p})$, respectively, but 
none is an integral solution,
given the constraints on $p_1,p_2$ and $p$.
\end{proof}
\begin{lemma}\label{lm:trivial_intersection_not_primitive}
Suppose $\pair$ is irreducible, and 
there exists a separating essential disk
$D$ with $I_D=(0,0)$ and $D\cap (l_+\cup l_-)$ two points.
Then $\vert p\vert>1$,
and either $l_+,l_-\subset\ComplHKA$ are primitive loops
or
$l_1,l_2$ are $(m,n)$-torus knots in $\sphere$ with $mn=p$
and one of $l_+, l_-$ representing the
$n$-th power of some primitive element of $\pi_1(\ComplHKA)$.
\end{lemma}
\begin{proof}
Note first that 
$D$ is separating by Lemma \ref{lm:trivial_intersection}; 
it may be assumed that $D\cap l_+=\emptyset$ and   
$D\cap l_-$ are two points.
Denote by $W_1,W_2$ the solid tori in 
$\Compl{\HK_A}-\openrnbhd{D}$ with $l_+\subset W_1$, 
and by $\alpha_i\subset W_i$, $i=1,2$, 
the two subarcs of $l_-$ cut off by $\openrnbhd{D}$. 
Note that $\alpha_i$ must be essential in 
$\partial W_i\cap \partial\ComplHKA$,
$i=1,2$, for otherwise $D$ could be isotoped such that
it is disjoint from $l_+\cup l_-$, contradicting
the irreducibility of $\pair$ (Remark \ref{rmk:lpm_primtive_trivial_HK}).
Consider an arc $\beta_i$ in the boundary  
of the disk $\rnbhd{D}\cap W_i$ with $\partial \beta_i =\partial \alpha_i$, and denote by $\hat{\alpha}_i$ 
the loop $\alpha_i\cup \beta_i$, $i=1,2$.

Observe that by Corollary \ref{cor:l_pm},
$l_+$ and hence $\hat{\alpha}_1$ 
have a finite slope of $\frac{m}{n}$ in $W_1$, $n>0$,
since $[l_+]$ is not trivial in $H_1(\ComplHKA)$,
while $\hat{\alpha}_2\subset W_2$ also has 
a finite slope of $\frac{m'}{n'}$, $n'>0$ 
since $[l_+]\neq \pm [l_-]$ in $H_1(\ComplHKA)$.
There exists a basis
of $H_1(\ComplHKA)$ given by the longitudes of $W_1,W_2$,
in terms of which
$[l_-]=(n,n'),[l_+]=(n,0)$.
As with \eqref{eq:change_of_basis_determinant}, 
we have $nn'=\pm p$; 
on the other hand, $mn=p$ due to $\lk{l_1}{l_2}=p$, so $m=\pm n'$.

On the homotopy level, 
the longitudes of $W_1,W_2$ induce a basis 
$\{u_1,u_2\}$ 
such that $l_+,l_-$ represent $u_1^n,u_1^n u_2^{n'}$; 
especially $l_+$ is the $n$-th power of some primitive element in 
$\pi_1(\ComplHKA)$. Since $\pi_1(\ComplHK)$
is an HNN-extension of $\pi_1(\ComplHKA)$ relative
to the isomorphism between $\pi_1(A_+), \pi_1(A_-)$ 
induced by $\rnbhd{A}$, 
we have the following presentation
\[\pi_1(\ComplHK)=\{u_1,u_2,t \mid  t u_1^{\pm n} t^{-1}=u_1^n u_2^{n'}\}.\]

Suppose $n=\absp$, and hence $m=\pm 1$. Then substitute $u_2$ with $w=u_1^{\absp} u_2$ 
gives us
\[\pi_1(\ComplHK)=
\{u_1,w,t \mid  t u_1^{\pm p} t^{-1}=w\},
\]
which implies that $\pair$ is trivial.
In the same way, one can show that 
$\pair$ is trivial if $p=\pm 1$ by replacing $u_2$ with $w=u_1 u_2$. 
Since $\pair$ is irreducible, we conclude that 
$\vert p\vert>1$ and $n<\vert p\vert$ and hence $\vert m\vert>1$.

Suppose $n=1$. Then there exists a meridian disk of $W_1$
that meets $l_+,l_-$ at one point each, and hence
both are primitive loops of $\ComplHKA$. 
On the other hand, if   
one of  $l_+,l_-\subset \ComplHKA$ is not primitive, 
then $l_+$ is an $(m,n)$-curve on the boundary of $W_1$ with 
$\vert m\vert, n>1$. 
Since $\pairA$ is trivial, the dual arc of $D_A$ in $\rnbhd{D_A}$
implies that $(\sphere, l_+)$ is a
tunnel number one knot.
Let $K$ be a core of $W_1$. 
If $(\sphere, K)$ is non-trivial,
then $(\sphere, l_+)$ is a tunnel number one $(m,n)$-cable knot,
contradicting 
the classification of tunnel number one satellite knots in \cite{MorSak:91}, for $W_1$ cannot be reembeded in $\sphere$ 
to make $l_+$ into an unknot in $\sphere$. 
Thus $(\sphere,K)$ is trivial, and $l_+, l_-$ and 
therefore $l_1,l_2$
are $(m,n)$-torus knots in $\sphere$. 
Since $A_+\subset \partial W_1$, $mn=\lk{l_1}{l_2}=p$.
\end{proof}


\section{Irreducibility and atoroidality}\label{sec:irre_atoro}

Throughout the section, $\pair$ is a handlebody-knot, and 
$A\subset\ComplHK$ is an annulus 
whose boundary components $l_1,l_2$ are
essential in $\partial \HK$, for instance, a type $3$-$3$ annulus\footnote{The condition on $l_1,l_2$ holds for all types of
annuli defined in \cite[Section $3$]{KodOzaGor:15}.}.  
As before, 
$\HK_A$ denotes the union of $\HK$ and a regular neighborhood $\rnbhd{A}$ of $A\subset\ComplHK$.

\subsection{Criteria for irreducibility and atoroidality}\label{subsec:criteria_irre_atoro}
 
\begin{lemma}\label{lm:irreducible_atoroidal_HK}
Suppose $\pair$ is irreducible and atoroidal. Then
$A$ is incompressible and $\ComplHKA$ is atoroidal.
If furthermore, $l_1,l_2$ are not parallel in $\partial\HK$,
then $A$ is essential.
\end{lemma}  
\begin{proof}
If there exists a compressing disk $D$ of $A$, then $D$ 
induces a disk $D'$ in $\Compl\HK$ bounded by $l_1$ or $l_2$;
$D'$ is essential since $l_1,l_2$ are essential in $\partial\HK$, contradicting the irreducibility of $\pair$; hence $A$ is incompressible.
To see $\ComplHKA$ is atoroidal, 
we let $T$ be a torus in $\ComplHKA$. 
By the assumption, it is compressible in $\ComplHK$,
and there exists a compressing disk $D$ of $T$ in $\ComplHK$.
Since $\partial D\subset T$ and $A$ is incompressible, $D\cap A$ contains only 
circles inessential in $A$.
Thus, one can isotope $D$ away from $A$, so $T$ is compressible in $\ComplHKA$.

Now suppose additionally that $l_1,l_2$ are non-parallel loops in $\partial \HK$. If there exists a $\partial$-compressing disk $D$
of $A$, then the disk component of 
the boundary of a regular neighborhood of 
$A\cup D$ in $\ComplHK$ 
is an essential disk in $\Compl\HK$,
contradicting the irreducibility of $\pair$.
Therefore, $A$ is essential.
\end{proof}
\begin{corollary}\label{cor:irreducible_atoroidal_HK}
Suppose $A$ is a type $3$-$3$ annulus with $\{\frac{p}{q},pq\}$ its slope pair, and $\pair$ is irreducible and atoroidal. 
Then $A$ is essential, and
$\pairA$ is either irreducible or trivial.
\end{corollary} 
\begin{proof}  
The first assertion follows directly from Lemma \ref{lm:irreducible_atoroidal_HK}. To see the second assertion,
note first that $\pairA$ is atoroidal by Lemma \ref{lm:irreducible_atoroidal_HK}. 
Now, if $\pairA$ is reducible, then there exists
an essential separating disk $D\subset \ComplHKA$.
The boundary of a regular neighborhood of $D\cup \HK_A$ in $\sphere$ 
consists of two tori $T_1,T_2$.
If $\pairA$ is non-trivial, one of $T_1,T_2$ 
bounds a non-trivial knot exterior in $\ComplHKA$, and 
is therefore incompressible in $\ComplHKA$,
contradicting the atoroidality of $\pairA$.  
\end{proof}

Conversely, the irreducibility and atoroidality of $\pair$
can be inferred from topological properties of
$\pairA$.
\begin{lemma}\label{lm:irreducibility:boundary_irreducible_ComplHKA} 
Suppose $\ComplHKA$ is $\partial$-irreducible.
Then $A$ is essential and $\pair$ is irreducible.
\end{lemma}
\begin{proof}
Observe first that every disk bounded by $l_+$ or $l_-$ in
$\partial\HK_A$ induces a disk bounded by $l_1,l_2$
in $\partial\HK$, and therefore, the assumption of $l_1,l_2$
being essential in $\partial\HK$ implies that
$l_+,l_-$ are essential in $\partial\HK_A$.

Suppose $A$ is compressible, and $D$
is a compressing disk; it may be assumed 
that $D\cap \rnbhd{A}$ is a regular neighborhood
of $\partial D$ in $D$, and the disk  
$D':=D-\openrnbhd{A}$ in $\ComplHKA$ is bounded by either $l_+$ or $l_-$.
By the irreducibility of $\pairA$, $\partial D'$
bounds a disk $D''$ in $\partial \HK_A$, contradicting that $l_+,l_-$ are essential in $\partial\HK_A$. 
Suppose $A$ is $\partial$-compressible,
and $D$ is a $\partial$-compressing disk; 
it may be assumed that $D\cap \rnbhd{A}$
is a regular neighborhood of $D\cap A$.
Then the disk $D':=D-\openrnbhd{A}$ 
intersects $l_+\cup l_-$ at one point, and hence
is essential in $\ComplHKA$,
contradicting the $\partial$-irreducibility of $\ComplHKA$.
Thus $A$ is essential.

Suppose $\pair$ is reducible, and $D$
is an essential disk in $\Compl\HK$.
Isotope $A$ such that $\# D\cap A$
is minimized. By the $\partial$-irreducibility of $\ComplHKA$,
$D\cap A\neq \emptyset$, and
since $A$ is essential
and $\Compl\HK$ is irreducible,
$D\cap A$ contains only arcs that
are inessential in $A$.

Let $\alpha\subset D\cap A$ be an outermost arc in $D$,
and $D_\alpha\subset D$ be an outermost disk cut off by $\alpha$.
It may be assumed that $D_\alpha\cap \rnbhd{A}\subsetneq D_\alpha$
and is a regular neighborhood of $\alpha$ in $D_\alpha$.
Thus, $D_\alpha':=D_\alpha-\openrnbhd{A}$ is a disk in $\Compl{\HK_A}$.
Since 
$\Compl{\HK_A}$ is $\partial$-irreducible, 
$\partial D_\alpha'$ bounds a disk $D_\alpha''$ in $\partial \HK_A$.
Isotoping $A$ through the $3$-ball bounded by $D_\alpha'\cup D_\alpha''$ removes $\alpha$ from $A\cap D$,
contradicting the minimality.

\end{proof}

\begin{corollary}\label{cor:irreducibility:irreducible_HKA} 
Suppose $A$ is a type $3$-$3$ annulus with $\{\frac{p}{q},pq\}$ its slope pair, and $\pairA$ is irreducible.
Then $A$ is essential and $\pair$ is irreducible.
\end{corollary}
 
\begin{remark}\label{rmk:irreducibility_other_types}
Lemma \ref{lm:irreducibility:boundary_irreducible_ComplHKA} 
still holds with $A$ replaced by a M\"obius band $M\subset\Compl\HK$
without conditions on $\partial M$ since $M$ is always incompressible in $\ComplHK$. 
Similarly, Corollary \ref{cor:irreducibility:irreducible_HKA} remains
valid if $A$ is replaced with a type $1$-$2$ M\"obius band $M$ \cite[Section $4$]{KodOzaGor:15}; note that $\pairM$ is also a handlebody-knot, where $\HK_M=\HK\cup\rnbhd{M}$. 
\end{remark}
  
\begin{lemma}\label{lm:irreducibility:trivial_HKA}  
Suppose $A$ is of type $3$-$3$ with a non-trivial boundary slope of $p$,
and $\pairA$ is trivial. 
If one of the following holds:
\begin{itemize}
\item $\vert p\vert=1$, and   
$\{l_+, l_-\}$ does not represent a basis of  $\pi_1(\ComplHKA)$.
\item 
$A$ satisfies the condition \eqref{cond:torus},
and neither of $l_+,l_-$ represents
the $\vert p\vert$-th power of some primitive element of $\pi_1(\ComplHKA)$.
\item $l_1,l_2$ 
are $(m,n)$-torus knots in $\sphere$ with $mn=p$, namely, failing the condition \eqref{cond:torus}, 
and one of $l_+,l_-$
does not represent the $n$-th or $\vert m\vert$-the power
of any primitive element in $\pi_1(\ComplHKA)$.
\end{itemize} 
then $A$ is essential and $\pair$ is irreducible. 
\end{lemma}
\begin{proof}
If there exists a compressing disk of $A$, then it
induces a disk in $\ComplHKA$
bounded by either $l_+$ or $l_-$,
contradicting that $[l_+], [l_-]\in H_1(\ComplHKA)$ 
are not trivial by Corollary \ref{cor:l_pm}.
If there exists a $\partial$-compressing disk of $A$,
then it induces a disk in $\ComplHKA$
that intersects $l_+\cup l_-$ at one point,
contradicting Corollary \ref{cor:non-separating_one_intersection}.
Thus, $A$ is essential.

Suppose $\pair$ is reducible, and $D\subset\Compl\HK$ is an essential disk.
Isotope $A$ such that $\#A\cap D$ is minimized.
Since $A$ is essential, $D\cap A$ contains only arcs
that are inessential in $A$.
Let $\alpha\subset A\cap D$ be an outermost arc
in $D$ and $D_\alpha\subset D$ an outermost disk cut off by $\alpha$.
It may be assumed that
$D_\alpha\cap \rnbhd{A}\subsetneq D_\alpha$ and 
is a regular neighborhood
of $\alpha$ in $D_\alpha$. Thus $D'_\alpha:=D_\alpha-\openrnbhd{A}$ is a disk
in $\ComplHKA$. Since $\alpha$ is inessential in $A$,
it may be assumed that $D_\alpha'$ is disjoint from $l_+,l_-$.
By the minimality, $D'_\alpha$ is essential
in $\ComplHKA$, but this is not possible by Corollaries \ref{cor:non_separating_disjoint} or \ref{cor:separating_disjoint}. 
\end{proof}

\cout{ 
In view of Remark \ref{rmk:lpm_primtive_trivial_HK}, 
the first condition can be replaced with $\pair$ being non-trivial.
}
In many cases, the homology version of Lemma \ref{lm:irreducibility:trivial_HKA} 
is sufficient to detect irreducibility.
\begin{corollary}\label{cor:irreducibility:trivial_HKA}
Let $A$ and $\pairA$ be as in Lemma \ref{lm:irreducibility:trivial_HKA}.  
Suppose $\vert p\vert>1$, 
$A$ satisfies the condition \eqref{cond:torus},
and neither of
$[l_+],[l_-]$ is 
the $\vert p\vert$-th multiple of some generator of $H_1(\ComplHKA)$.
Then $A$ is essential, and $\pair$ is irreducible. 
\end{corollary}

The following atoroidality criterion is 
a corollary of \cite[Theorem $3.3$]{Wan:21}.
\begin{lemma}\label{lm:atoroidality:trivial_HKA}
Suppose $A$ is a type $3$-$3$ annulus with $\{\frac{p}{q},pq\}$ its slope pair, and $\pair$ is irreducible. If $\pairA$ is trivial,
then $\pair$ is atoroidal.
\end{lemma}

The proof of Lemma \ref{lm:atoroidality:trivial_HKA} implies 
the following result for general annuli. For the sake of completeness, 
we recall the argument in \cite[Theorem $3.3$]{Wan:21} below.

\begin{lemma}\label{lm:atoroidality_general}
Suppose $A\subset \Compl\HK$ is incompressible and $\ComplHKA$ atoroidal. 
If $\Compl\HK$ is toroidal, 
then $\partial A$ is an $(2m,2n)$-torus link in $\sphere$, $\vert m\vert, n>1$. 
If furthermore $A$ is of type $3$-$3$,
then $A$ has a non-trivial boundary slope of $p=mn$.
\end{lemma}
\begin{proof}
Let $T$ be an incompressible torus in
$\Compl\HK$ that minimizes
\[\{\# T\cap A\mid T\subset\Compl\HK \text{ an incompressible torus}\}.\]
Denote by $U$ the solid torus bounded by $T$; note that 
$\HK\subset U$.
By the incompressibility of $A,T$, 
every circle in $A\cap T$ is essential in both $A$ and $T$.

\textbf{Case $1$: $T\cap A\subset U$ is meridional.} 
There exists an annulus $B\subset A$ with $B\cap T=\partial B$
and $B\not\subset U$. Let $B'\subset T$ be an annulus 
cut off by $\partial B$ and 
$T_B:=B\cup B'$.   
If $V$ is the component of $\sphere-T_B$ not containing $\HK$,
then $V$ is a solid torus
since $T_B$ has less intersection with $A$ than $T$ does. 
On the other hand, because 
$T\cap A$ are meridional in $\partial U$, 
any essential loop of $B'$ bounds a disk in $U$,
and is therefore a longitude of $V$. In particular, 
one can isotope $T$ through $V$ to decrease
$\# A\cap T$, a contradiction.

\textbf{Case $2$: $T\cap A\subset U$ is non-meridional.}
We first prove the following claim:

\centerline{\textbf{$\# T\cap A$ is at most $2$}.}

\begin{figure}[t]
\centering
\def\svgwidth{.42 \columnwidth}
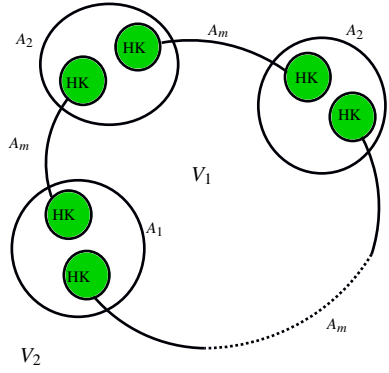
\caption{Schematic diagram of $A\cap T$.} 
\label{fig:intersection_T_and_A}
\end{figure}

Suppose $\# T\cap A>2$. Then there is an annulus $B\subset A$
with $B\cap T=\partial B$ and $B\subset U$.
Since $\partial B\subset U$ is not meridional,
$B\subset U$ is incompressible, 
and divides $U$ into two solid tori $U_1, U_2$. 
Let $U_1$ be the solid torus containing $\HK$.
If $\Compl{U_1}$ is not a solid torus,
then $\partial U_1$ induces an incompressible torus
having less intersection with $A$, contradicting the minimality.
On the other hand, if $\Compl{U_1}$ is a solid torus, 
then 
$\Compl{U}$ is also a solid torus since 
$\Compl{U}$ can be obtained by 
cutting $\Compl{U_1}$ along 
the incompressible 
annulus $B':=\Compl{U}\cap U_2=\partial U_2-\mathring{B}
\subset \Compl{U_1}$, contradicting the incompressibility of $T$. 
 
Therefore, $\# T\cap A=2$, and $T\cap A$
are two parallel 
(resp.\ parallel non-meridional) essential loops in $A$ (resp.\ in $T$).
Especially, $T\cap A$ cuts $T$ into two annuli $A_1,A_2$.
Denote by $A_m\subset A$ the annulus with $A_m\cap T=\partial A_m$.
Note that $A_m$ is necessarily in $\Compl U$.
By the atoroidality of $\HK_A$, the components $V_i$ bounded by 
$A_m\cup A_i$, $i=1,2$, with $\HK\not\subset V_i$ are solid tori (see Fig.\ \ref{fig:intersection_T_and_A}).

Suppose one of $\pi_1(A_m)\rightarrow \pi_1(V_i)$, $i=1,2$, is an isomorphism. Then $V_1\cup V_2=\Compl U$,
is a solid torus, contradicting the incompressibility of $T=\partial U$.
Now, because neither of $A_m\rightarrow V_i$, $i=1,2$,
induces an isomorphism on $\pi_1$, the core of $U$
is an $(m,n)$-torus knot in $\sphere$
by the classification of Seifert fiber structure  
of $\sphere$ \cite{Sei:33}.
Since the link $\partial A=l_1\cup l_2\subset \sphere$ is isotopic to $\partial A_m\subset\sphere$,
$l_1\cup l_2$ is an $(2m,2n)$-torus link in $\sphere$ 
with $\vert m\vert, n>1$.

Suppose $A$ is of type $3$-$3$ with $\{\frac{p}{q},pq\}$ its slope pair, 
$q>0$, and $\disk_A\subset \HK$
is a meridian disk disjoint from $A$.
By Lemma \ref{lm:separating_disks}, $\HK-\openrnbhd{\disk_A}$
are two solid tori $W_1,W_2$ with $l_i\subset W_i$, $i=1,2$.
It may be assumed that the slope of $l_1\subset W_1$ is $\frac{p}{q}$,
and therefore $[l_1]$ is the $q$-th multiple of some generator of $H_1(W_1)$.
On the other hand, its image under the homomorphism 
$H_1(W_1)\rightarrow H_1(U)$ induced by
the inclusion is a generator of $H_1(U)$, so $q=1$, 
and the slope pair of $A$ is $\{p,p\}$ with
$p=\lk{l_1}{l_2}=mn\neq 0$. 
\end{proof}

\begin{corollary}\label{cor:atoroidality_mobius_band}
Given a M\"obius band $M$ in $\ComplHK$.
Suppose $\Compl{\HK_M}$ is atoroidal, and $\Compl\HK$
is toroidal.
Then $\partial M$ is an $(m,n)$-torus knot. 
\end{corollary}
\begin{proof}
The annulus $A:=\rnbhd{M}\cap\partial\Compl{\HK_M}$ 
is incompressible in $\ComplHK$.   
\end{proof} 
 

\subsection{Examples}\label{subsec:examples_irreatoro}
Here we present a 
construction of irreducible, atoroidal handlebody-knots 
whose exteriors contain a type $3$-$3$ annulus,
and prove the irreducibility and 
atoroidality of handlebody-knots in Fig.\ \ref{fig:intro:nonuniqueness},
employing criteria developed in Section \ref{subsec:criteria_irre_atoro}. 


Recall that given a knot $K\subset\sphere$, a \emph{tunnel} of $K$  
is an arc $\tau$ in $\sphere$ with $\tau\cap K=\partial \tau$; 
$\tau$ is called \emph{unknotting} if $\pairKtau$ is a
trivial handlebody-knot, where
$\rnbhd{K\cup\tau}$ is a regular neighborhood of $K\cup\tau$ in $\sphere$.

A $p$-annulus $\mathcal{A}$ associated to $(K,\tau)$, $p\in\mathbb{Z}$,
is an annulus in $\sphere$ that satisfies the following conditions:
\begin{itemize}
\item $\lk{l_1}{l_2}=p$, 
where $l_1,l_2$ are the components of $\partial \mathcal{A}$;
\item $K\subset \mathring{\mathcal{A}}$ is an essential loop of $\mathcal{A}$; 
\item $\mathcal{A}\cap \tau$ is a regular neighborhood of $\partial \tau$ in $\tau$;
\item $\tau$ meets both $l_1$ and $l_2$.
\end{itemize}
Given a $p$-annulus $\mathcal{A}$ associated to $(K,\tau)$,
the handlebody-knot $\pairAtau$ is given by a regular neighborhood
\[\HKAtau:=\rnbhd{(\tau-\mathring{\mathcal{A}})\cup l_1\cup l_2}\]
of $(\tau-\mathring{\mathcal{A}})\cup l_1\cup l_2$.
It may be assumed that $\HKAtau \cap \mathcal{A}$
is a regular neighborhood of $\partial \mathcal{A}$ in $\mathcal{A}$,
and hence $A:=\overline{\mathcal{A}-\big(
\HKAtau\cap \mathcal{A}
\big)}\subset
\Compl{\HKAtau}$ 
is a type $3$-$3$ annulus with a boundary slope of $p$,
and $\pairAtauA=\pairKtau$. 

\subsubsection{Examples: $\pairAtauA$ is trivial}  
\begin{figure}[h]
\begin{subfigure}{0.45\textwidth}
\centering
\def\svgwidth{.7\columnwidth}
\begingroup%
  \makeatletter%
  \providecommand\color[2][]{%
    \errmessage{(Inkscape) Color is used for the text in Inkscape, but the package 'color.sty' is not loaded}%
    \renewcommand\color[2][]{}%
  }%
  \providecommand\transparent[1]{%
    \errmessage{(Inkscape) Transparency is used (non-zero) for the text in Inkscape, but the package 'transparent.sty' is not loaded}%
    \renewcommand\transparent[1]{}%
  }%
  \providecommand\rotatebox[2]{#2}%
  \newcommand*\fsize{\dimexpr\f@size pt\relax}%
  \newcommand*\lineheight[1]{\fontsize{\fsize}{#1\fsize}\selectfont}%
  \ifx\svgwidth\undefined%
    \setlength{\unitlength}{992.12598425bp}%
    \ifx\svgscale\undefined%
      \relax%
    \else%
      \setlength{\unitlength}{\unitlength * \real{\svgscale}}%
    \fi%
  \else%
    \setlength{\unitlength}{\svgwidth}%
  \fi%
  \global\let\svgwidth\undefined%
  \global\let\svgscale\undefined%
  \makeatother%
  \begin{picture}(1,1)%
    \lineheight{1}%
    \setlength\tabcolsep{0pt}%
    \put(0,0){\includegraphics[width=\unitlength,page=1]{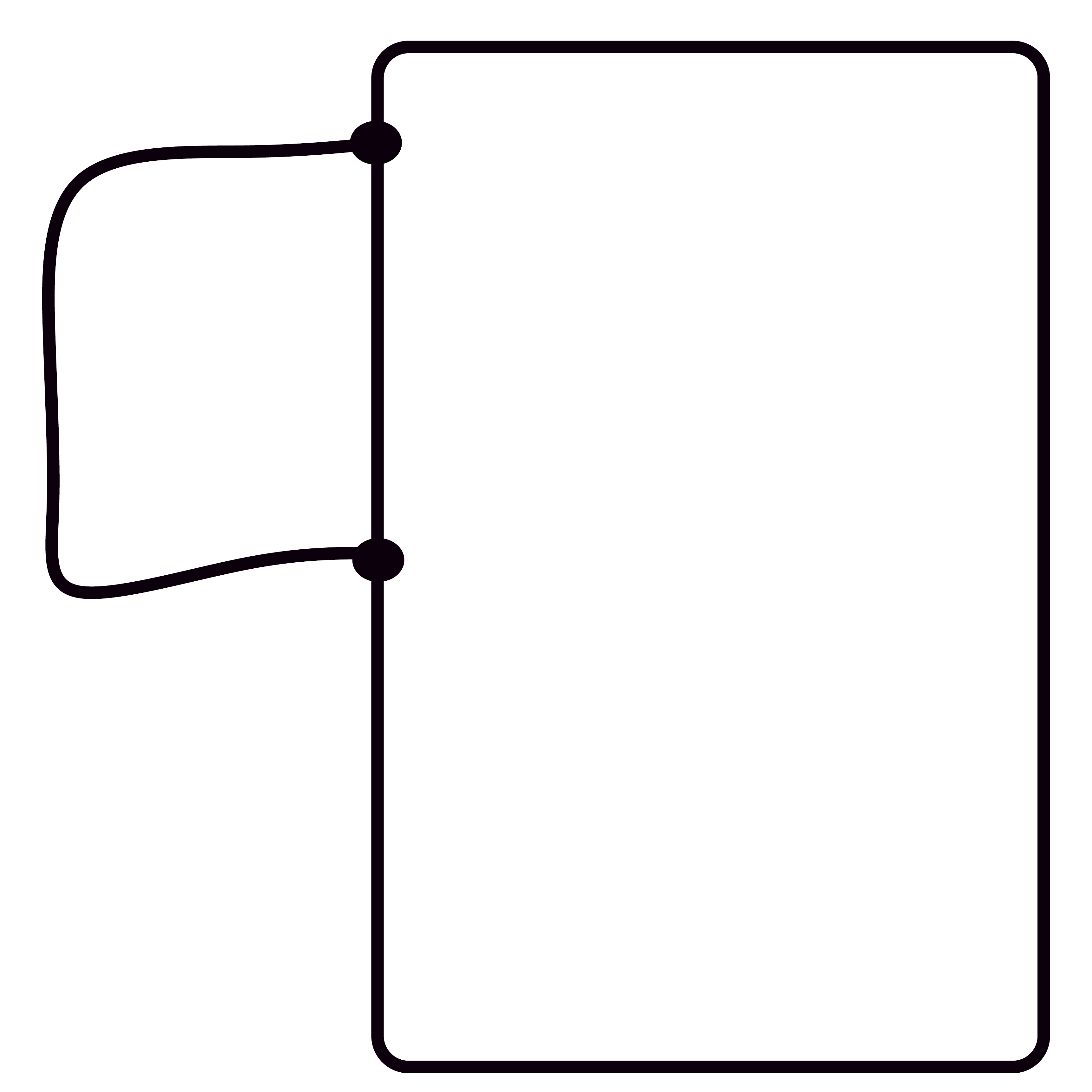}}%
    \put(0.22974427,0.11232944){\color[rgb]{0,0,0}\makebox(0,0)[lt]{\lineheight{1.25}\smash{\begin{tabular}[t]{l}{$K$}\end{tabular}}}}%
    \put(0.05817515,0.69299937){\color[rgb]{0,0,0}\makebox(0,0)[lt]{\lineheight{1.25}\smash{\begin{tabular}[t]{l}{\footnotesize $\tau$}\end{tabular}}}}%
    \put(0,0){\includegraphics[width=\unitlength,page=2]{unknotting_tunnel_trivial_knot.pdf}}%
    \put(0.39908558,0.80686637){\color[rgb]{0,0,0}\makebox(0,0)[lt]{\lineheight{1.25}\smash{\begin{tabular}[t]{l}{\footnotesize $x_1$}\end{tabular}}}}%
    \put(0.39897638,0.34798092){\color[rgb]{0,0,0}\makebox(0,0)[lt]{\lineheight{1.25}\smash{\begin{tabular}[t]{l}{\footnotesize $x_2$}\end{tabular}}}}%
  \end{picture}%
\endgroup%

\caption{Tunnel $\tau$ of the trivial knot $K$.}
\label{fig:tunnel_trivial_knot}
\end{subfigure}
\begin{subfigure}{0.5\textwidth}
\centering
\def\svgwidth{.7 \columnwidth}
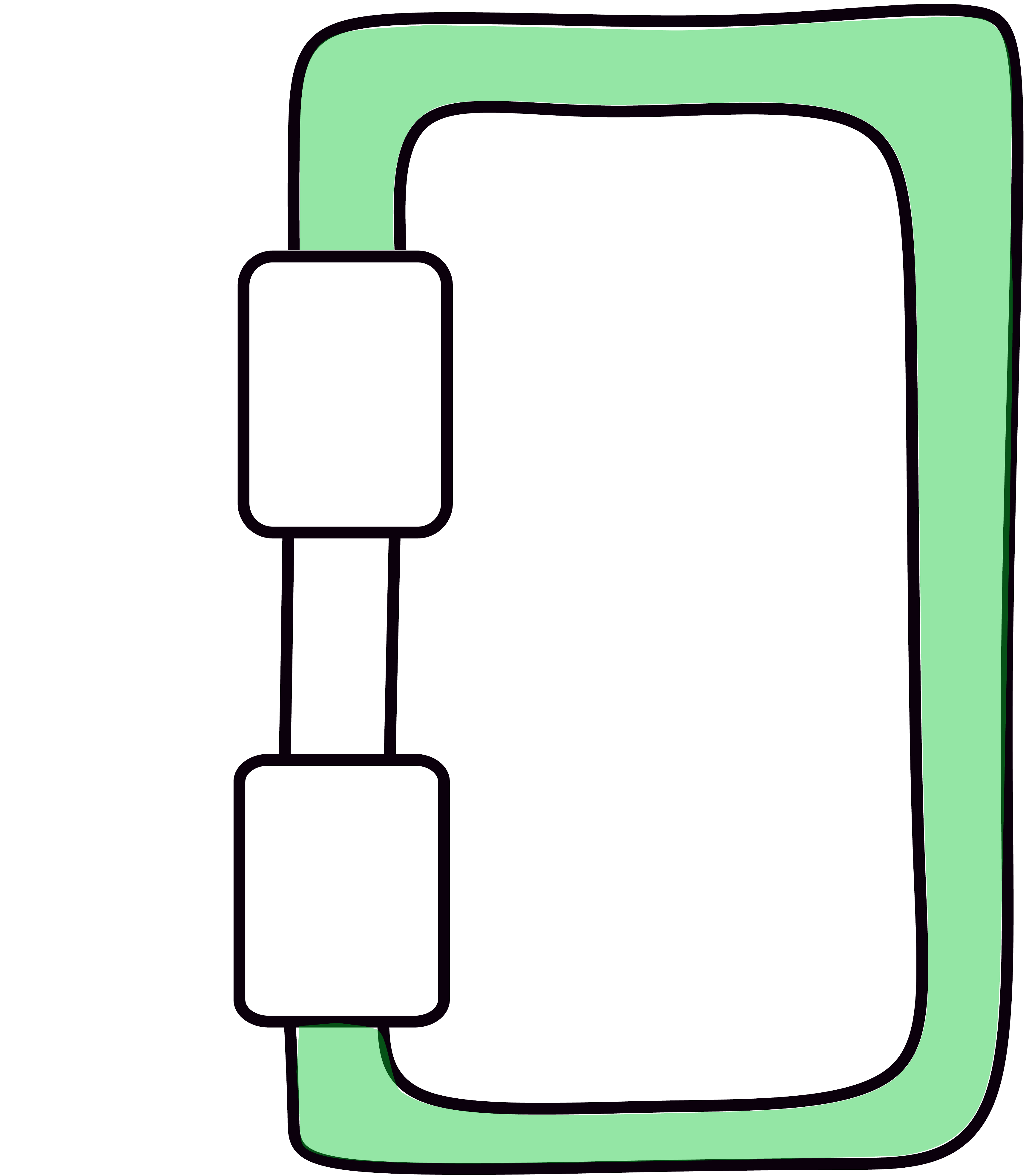
\caption{$p$-annulus $\mathcal{A}$ associated to $(K,\tau)$.}
\label{fig:annulus_tunnel_trivial_knot}
\end{subfigure}
\medskip
\begin{subfigure}{0.47\textwidth}
\centering
\def\svgwidth{.7 \columnwidth}
\begingroup%
  \makeatletter%
  \providecommand\color[2][]{%
    \errmessage{(Inkscape) Color is used for the text in Inkscape, but the package 'color.sty' is not loaded}%
    \renewcommand\color[2][]{}%
  }%
  \providecommand\transparent[1]{%
    \errmessage{(Inkscape) Transparency is used (non-zero) for the text in Inkscape, but the package 'transparent.sty' is not loaded}%
    \renewcommand\transparent[1]{}%
  }%
  \providecommand\rotatebox[2]{#2}%
  \newcommand*\fsize{\dimexpr\f@size pt\relax}%
  \newcommand*\lineheight[1]{\fontsize{\fsize}{#1\fsize}\selectfont}%
  \ifx\svgwidth\undefined%
    \setlength{\unitlength}{992.12598425bp}%
    \ifx\svgscale\undefined%
      \relax%
    \else%
      \setlength{\unitlength}{\unitlength * \real{\svgscale}}%
    \fi%
  \else%
    \setlength{\unitlength}{\svgwidth}%
  \fi%
  \global\let\svgwidth\undefined%
  \global\let\svgscale\undefined%
  \makeatother%
  \begin{picture}(1,1.14285714)%
    \lineheight{1}%
    \setlength\tabcolsep{0pt}%
    \put(0,0){\includegraphics[width=\unitlength,page=1]{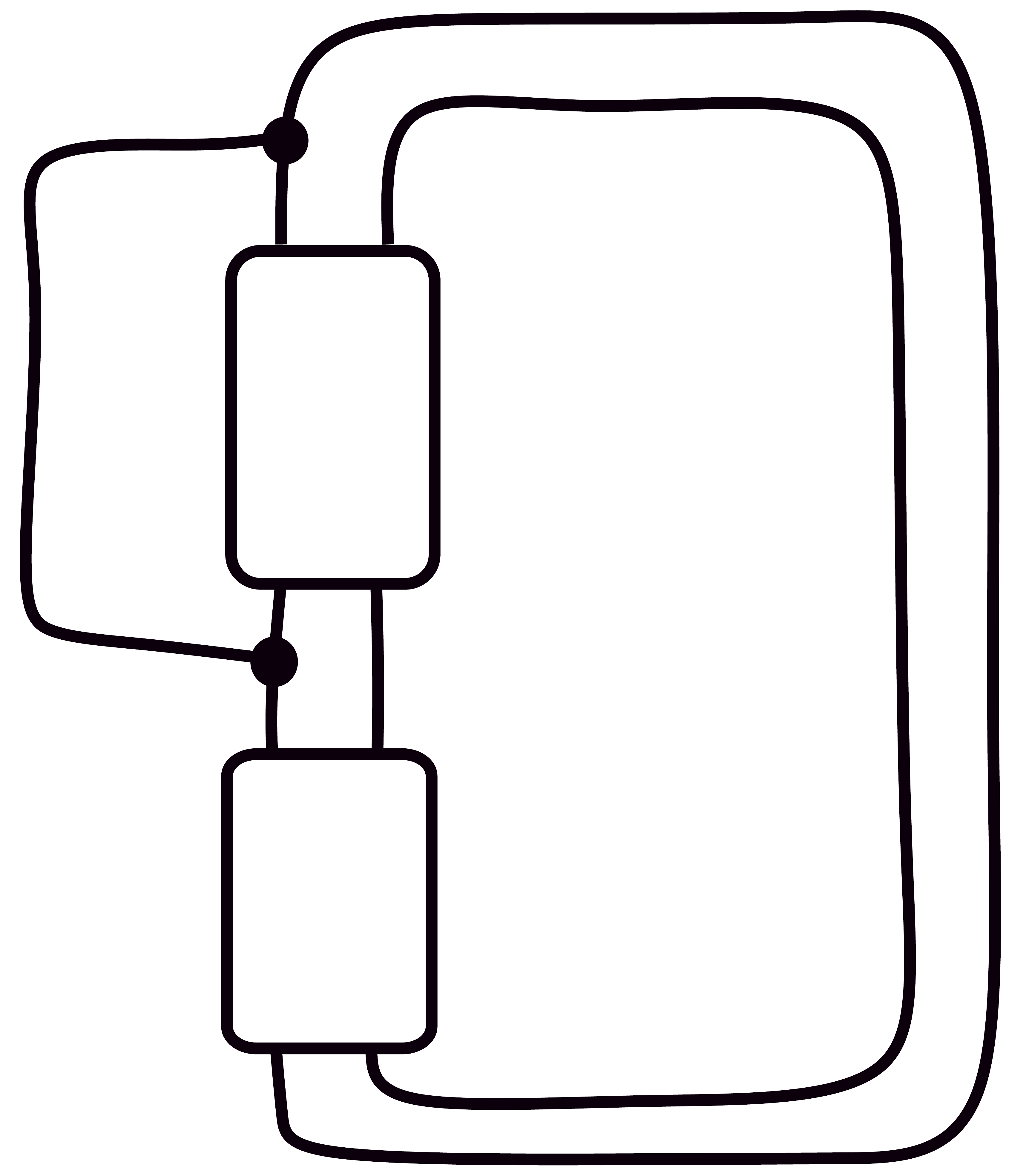}}%
    \put(0.28084636,0.80375734){\color[rgb]{0,0,0}\makebox(0,0)[lt]{\lineheight{1.25}\smash{\begin{tabular}[t]{l}{\footnotesize $\mu$}\end{tabular}}}}%
    \put(0.26005953,0.74011751){\color[rgb]{0,0,0}\makebox(0,0)[lt]{\lineheight{1.25}\smash{\begin{tabular}[t]{l}{\tiny  half }\end{tabular}}}}%
    \put(0.24844219,0.67447122){\color[rgb]{0,0,0}\makebox(0,0)[lt]{\lineheight{1.25}\smash{\begin{tabular}[t]{l}{\tiny  twists }\end{tabular}}}}%
    \put(0.279392,0.31497598){\color[rgb]{0,0,0}\makebox(0,0)[lt]{\lineheight{1.25}\smash{\begin{tabular}[t]{l}{\footnotesize $\nu$}\end{tabular}}}}%
    \put(0.25860517,0.25738382){\color[rgb]{0,0,0}\makebox(0,0)[lt]{\lineheight{1.25}\smash{\begin{tabular}[t]{l}{\tiny  half }\end{tabular}}}}%
    \put(0.24698783,0.19929701){\color[rgb]{0,0,0}\makebox(0,0)[lt]{\lineheight{1.25}\smash{\begin{tabular}[t]{l}{\tiny  twists }\end{tabular}}}}%
  \end{picture}%
\endgroup%

\caption{$\pairTmn$.}
\label{fig:example_Tmn}
\end{subfigure}
\begin{subfigure}{0.47\textwidth}
\centering
\def\svgwidth{.7 \columnwidth}
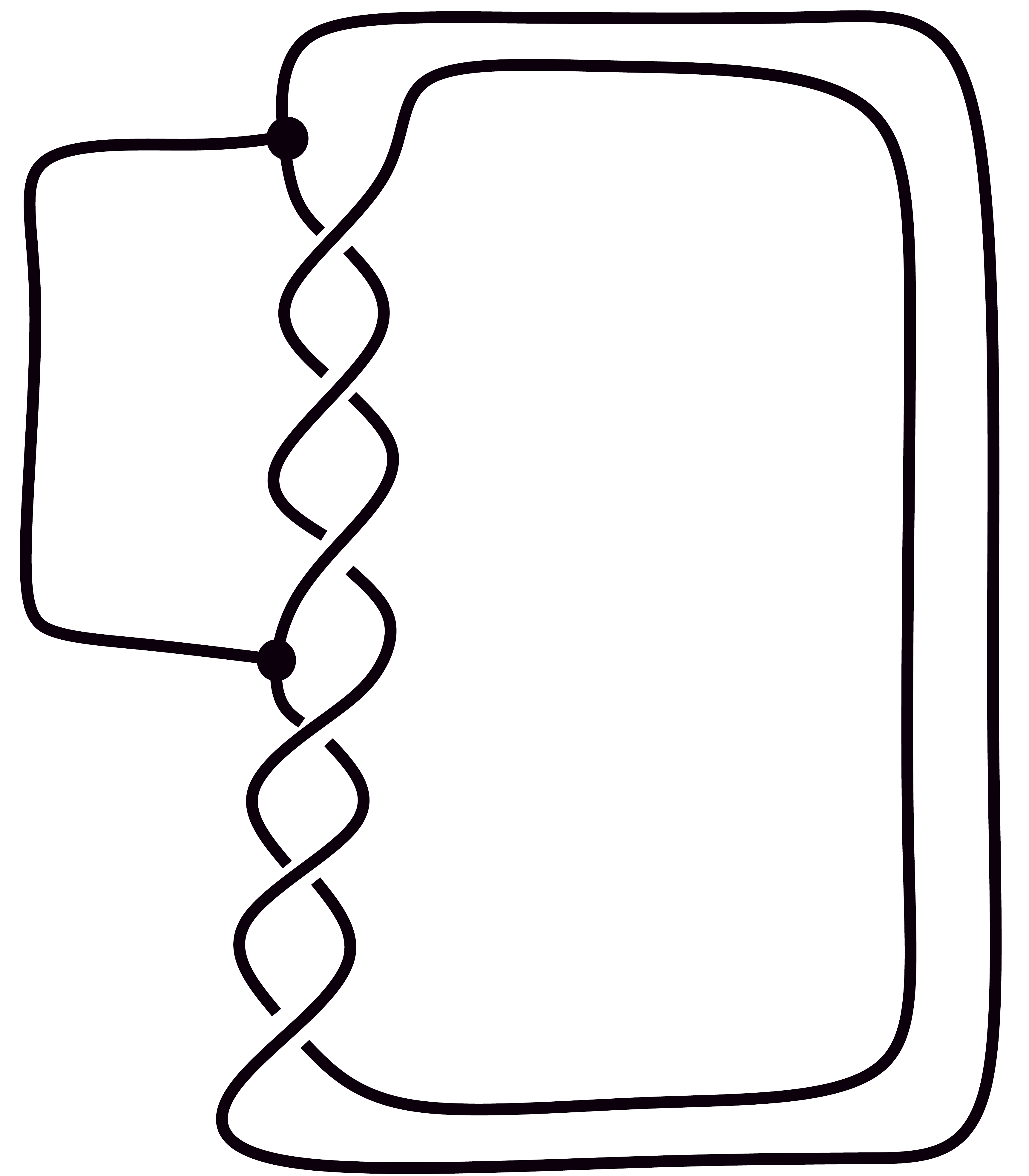
\caption{$(\sphere,\mathcal{T}_{3,3})$.}
\label{fig:example_T_33}
\end{subfigure}
\caption{Construction of the handlebody-knot family $\pairTmn$.}
\end{figure}
 
We construct here
a family of irreducible, atoroidal handlebody-knots
whose exteriors admit an unknotting type $3$-$3$ annulus. 
Start with a trivial knot $K$
and an unknotting tunnel $\tau$ as shown in Fig.\ \ref{fig:tunnel_trivial_knot}.
Take a $p$-annulus associated to $(K,\tau)$
as in Fig.\ \ref{fig:annulus_tunnel_trivial_knot},
where $\mu,\nu$ are odd integers and $\mu+\nu=2p$.
Denote by $\Tmn$ the handlebody $\HKAtau$ 
produced with the data (see Fig.\ \ref{fig:example_Tmn}; 
a spine of $\pairTthreethree$ 
is illustrated in 
Fig.\ \ref{fig:example_T_33}).

Orient $A$ as shown in Fig.\ \ref{fig:annulus_tunnel_trivial_knot}. 
Then 
\[[l_+]=(\frac{\mu+1}{2},\frac{\nu-1}{2}), \quad 
[l_-]=(\frac{\mu-1}{2},\frac{\nu+1}{2})\in H_1(\ComplTmnA),\]
in terms of the meridional basis 
of $H_1(\ComplTmnA)$ 
given by $x_1,x_2$ in Fig.\ \ref{fig:tunnel_trivial_knot}.
Since $l_+,l_-$, being isotopic to $K$, are trivial knots,
Corollary \ref{cor:irreducibility:trivial_HKA}
and Lemma \ref{lm:atoroidality:trivial_HKA} imply the following. 
\begin{corollary}\label{cor:example_Tmn_irre_atoro}
If $\frac{\mu\pm 1}{2},\frac{\nu\pm 1}{2}$
are not divisible by $p$, 
then the handlebody-knot $\pairTmn$ 
is irreducible and atoroidal.
\end{corollary}
We remark that $\pairTmn$ is trivial when $\mu=\pm1$ or 
$\nu =\pm 1$.  

\subsubsection{Examples: $\pairAtauA$ is irreducible} 
\begin{figure}[t]
\begin{subfigure}{0.48\textwidth}
\centering
\def\svgwidth{.7\columnwidth}
\begingroup%
  \makeatletter%
  \providecommand\color[2][]{%
    \errmessage{(Inkscape) Color is used for the text in Inkscape, but the package 'color.sty' is not loaded}%
    \renewcommand\color[2][]{}%
  }%
  \providecommand\transparent[1]{%
    \errmessage{(Inkscape) Transparency is used (non-zero) for the text in Inkscape, but the package 'transparent.sty' is not loaded}%
    \renewcommand\transparent[1]{}%
  }%
  \providecommand\rotatebox[2]{#2}%
  \newcommand*\fsize{\dimexpr\f@size pt\relax}%
  \newcommand*\lineheight[1]{\fontsize{\fsize}{#1\fsize}\selectfont}%
  \ifx\svgwidth\undefined%
    \setlength{\unitlength}{1133.85826772bp}%
    \ifx\svgscale\undefined%
      \relax%
    \else%
      \setlength{\unitlength}{\unitlength * \real{\svgscale}}%
    \fi%
  \else%
    \setlength{\unitlength}{\svgwidth}%
  \fi%
  \global\let\svgwidth\undefined%
  \global\let\svgscale\undefined%
  \makeatother%
  \begin{picture}(1,0.75)%
    \lineheight{1}%
    \setlength\tabcolsep{0pt}%
    \put(0,0){\includegraphics[width=\unitlength,page=1]{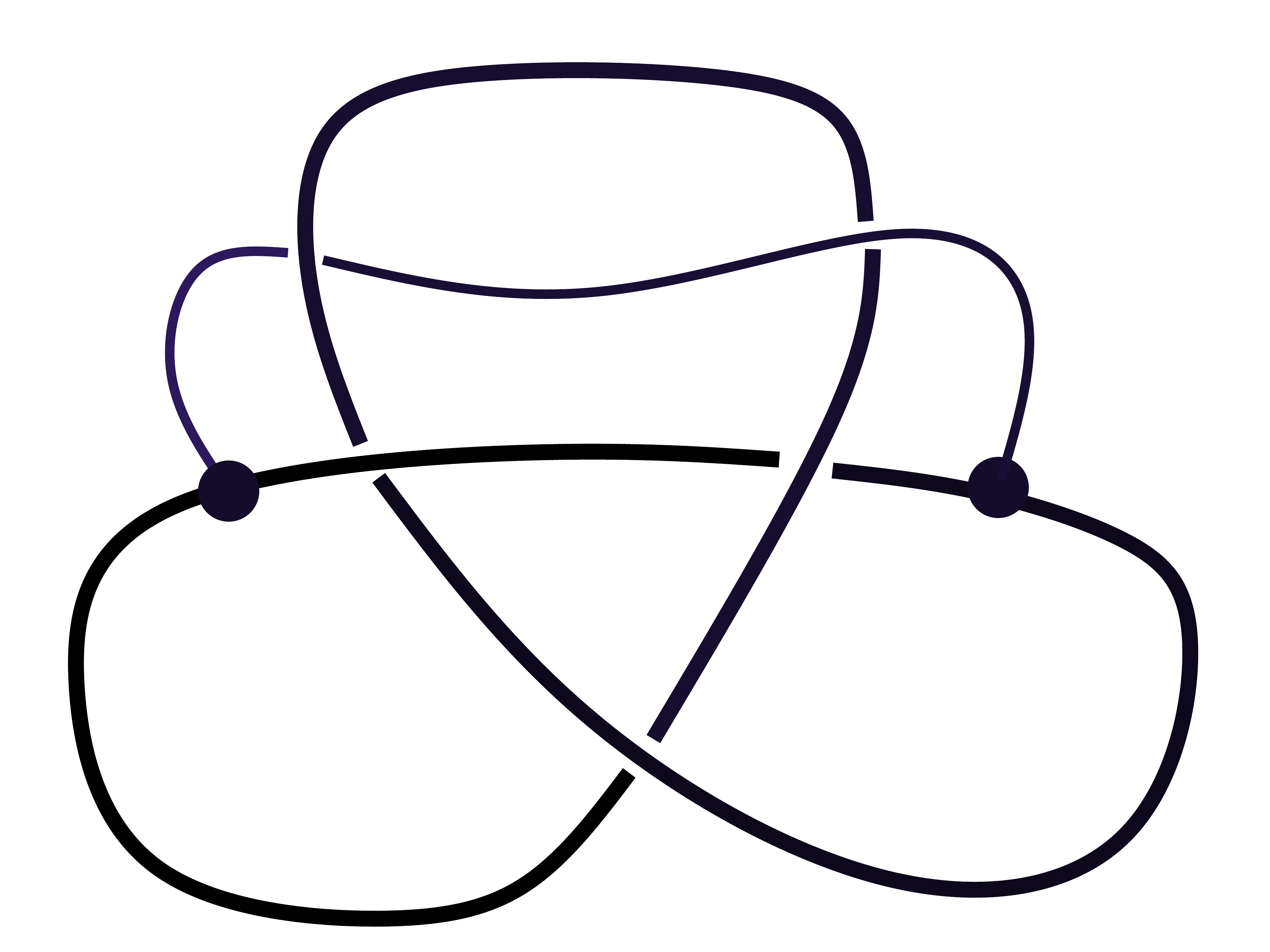}}%
    \put(0.4341568,0.02411308){\color[rgb]{0,0,0}\makebox(0,0)[lt]{\lineheight{1.25}\smash{\begin{tabular}[t]{l}{\footnotesize $K$}\end{tabular}}}}%
    \put(0.55253937,0.49202921){\color[rgb]{0,0,0}\makebox(0,0)[lt]{\lineheight{1.25}\smash{\begin{tabular}[t]{l}{\footnotesize $\tau$}\end{tabular}}}}%
    \put(0,0){\includegraphics[width=\unitlength,page=2]{spine_5_2.pdf}}%
    \put(0.33430557,0.71121596){\color[rgb]{0,0,0}\makebox(0,0)[lt]{\lineheight{1.25}\smash{\begin{tabular}[t]{l}{\footnotesize $x_1$}\end{tabular}}}}%
    \put(0.33193695,0.41102535){\color[rgb]{0,0,0}\makebox(0,0)[lt]{\lineheight{1.25}\smash{\begin{tabular}[t]{l}{\footnotesize $x_2$}\end{tabular}}}}%
  \end{picture}%
\endgroup%

\caption{$K,\tau\subset\Theta$.}
\label{fig:spine_5_2}
\end{subfigure}
\begin{subfigure}{0.48\textwidth}
\centering
\def\svgwidth{.7 \columnwidth}
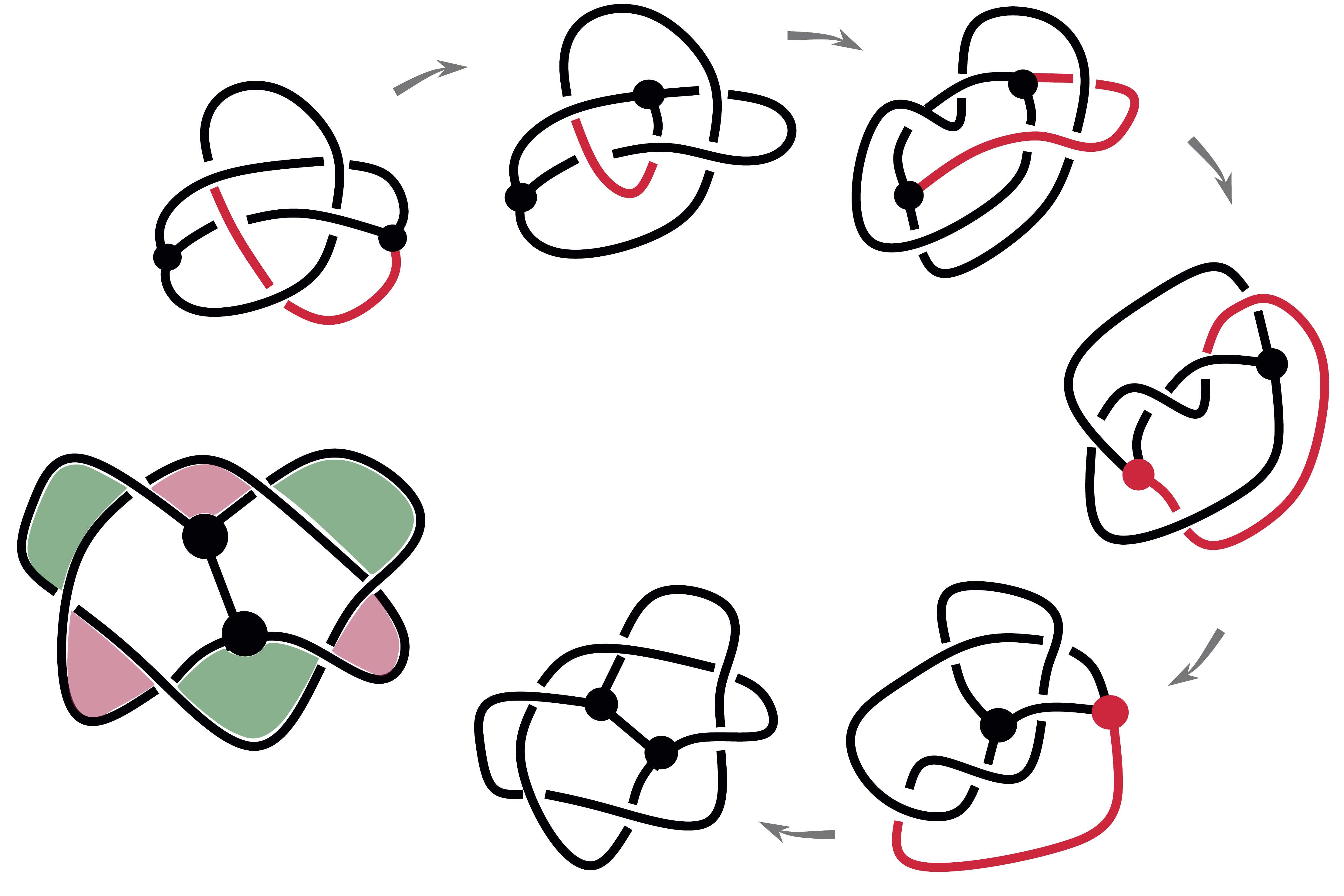
\caption{Equivalence.}
\label{fig:equivalence_to_m5_2}
\end{subfigure}
\medskip

\begin{subfigure}{0.48\textwidth}
\centering
\def\svgwidth{.7\columnwidth}
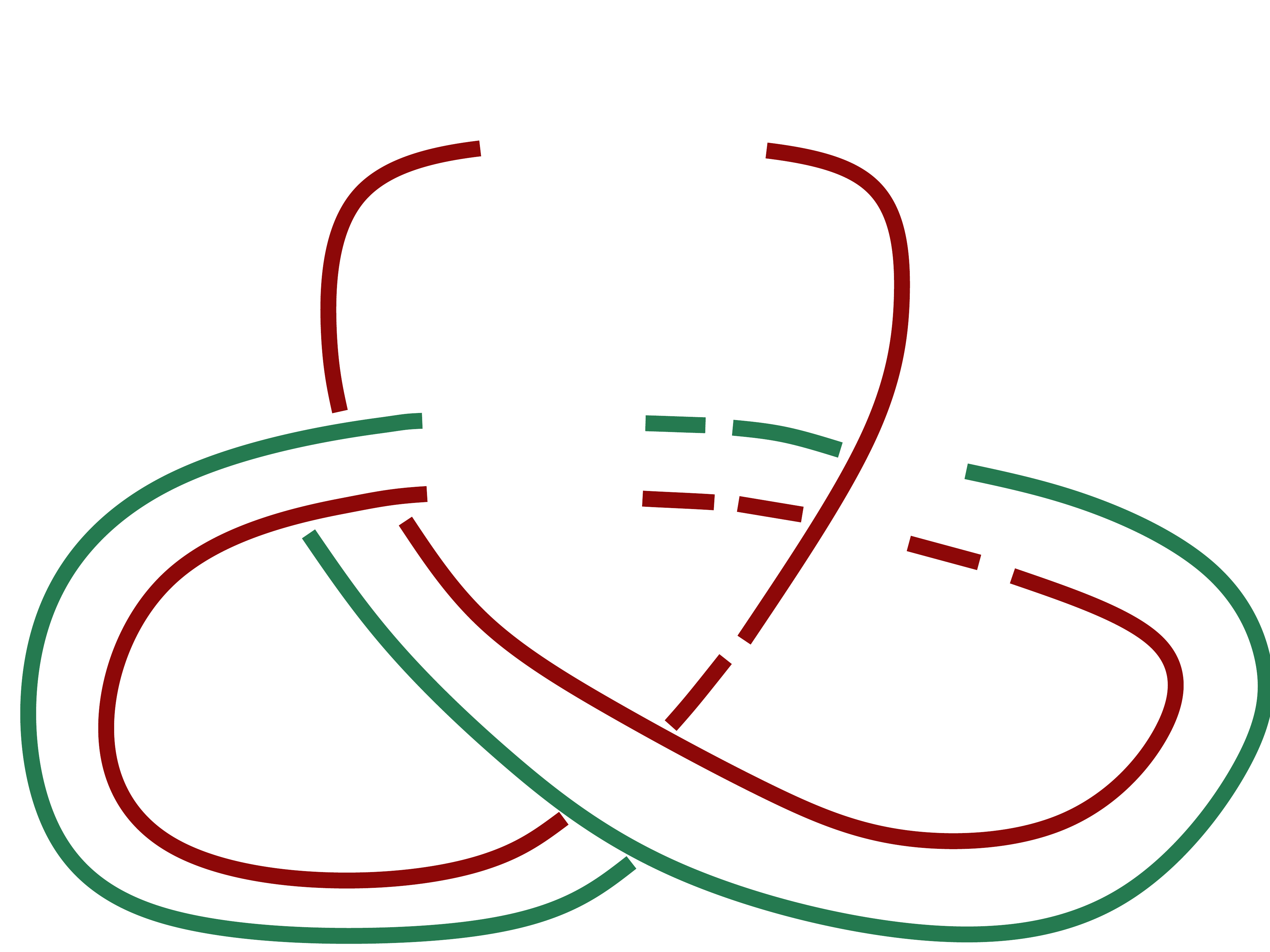
\caption{$p$-annulus $\mathcal{A}$ associated to $(K,\tau)$.}
\label{fig:mn_twisted_annulus}
\end{subfigure}
\begin{subfigure}{0.48\textwidth}
\centering
\def\svgwidth{.7 \columnwidth}
\begingroup%
  \makeatletter%
  \providecommand\color[2][]{%
    \errmessage{(Inkscape) Color is used for the text in Inkscape, but the package 'color.sty' is not loaded}%
    \renewcommand\color[2][]{}%
  }%
  \providecommand\transparent[1]{%
    \errmessage{(Inkscape) Transparency is used (non-zero) for the text in Inkscape, but the package 'transparent.sty' is not loaded}%
    \renewcommand\transparent[1]{}%
  }%
  \providecommand\rotatebox[2]{#2}%
  \newcommand*\fsize{\dimexpr\f@size pt\relax}%
  \newcommand*\lineheight[1]{\fontsize{\fsize}{#1\fsize}\selectfont}%
  \ifx\svgwidth\undefined%
    \setlength{\unitlength}{1133.85826772bp}%
    \ifx\svgscale\undefined%
      \relax%
    \else%
      \setlength{\unitlength}{\unitlength * \real{\svgscale}}%
    \fi%
  \else%
    \setlength{\unitlength}{\svgwidth}%
  \fi%
  \global\let\svgwidth\undefined%
  \global\let\svgscale\undefined%
  \makeatother%
  \begin{picture}(1,0.75)%
    \lineheight{1}%
    \setlength\tabcolsep{0pt}%
    \put(0,0){\includegraphics[width=\unitlength,page=1]{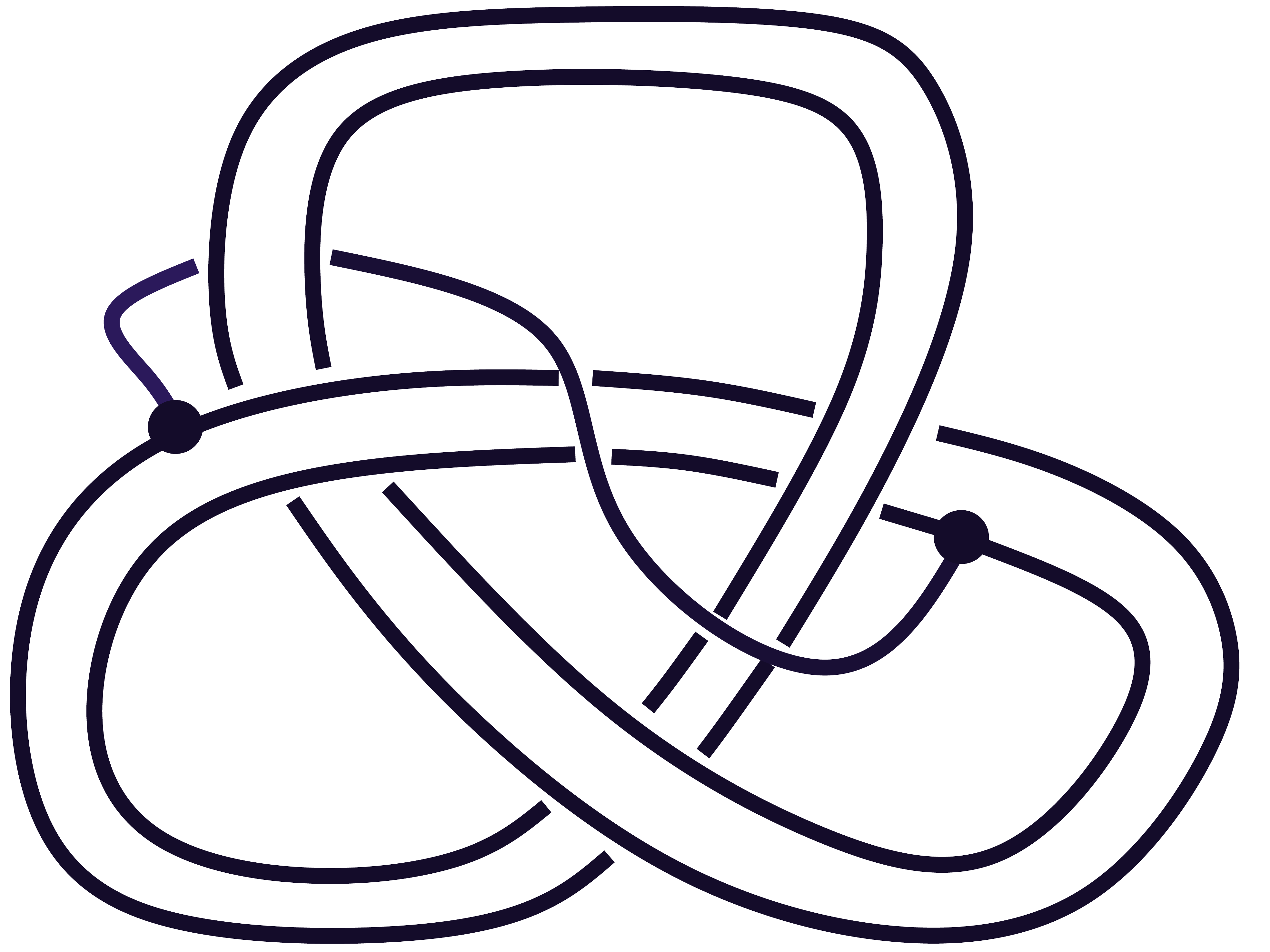}}%
    \put(0.55749683,0.62134687){\color[rgb]{0,0,0}\makebox(0,0)[lt]{\lineheight{1.25}\smash{\begin{tabular}[t]{l}{\footnotesize $l_1$}\end{tabular}}}}%
    \put(0.779095,0.58637479){\color[rgb]{0,0,0}\makebox(0,0)[lt]{\lineheight{1.25}\smash{\begin{tabular}[t]{l}{\footnotesize $l_2$}\end{tabular}}}}%
  \end{picture}%
\endgroup%

\caption{$(\sphere,\mathcal{I}_{0,0})$.}
\label{fig:example_I_00}
\end{subfigure}
\caption{Construction of handlebody family $\pairImn$.}
\end{figure}

Consider the handlebody-knot
$(\sphere,5_2)$ in the handlebody-knot table 
\cite{IshKisMorSuz:12}, whose spine $\Theta$ 
is depicted in Fig.\ \ref{fig:spine_5_2}. 
Take $K$ to be the constituent trefoil knot
and $\tau$ the other arc in $\Theta$.
Choose a $p$-annulus $A$ associated to $(K,\tau)$
as in Fig.\ \ref{fig:mn_twisted_annulus},
where $p=\mu+\nu+3$. Denote by $\Imn$ the resulting handlebody
$\HKAtau$. 

Note that the irreducibility and atoroidality 
$(\sphere,5_2)$ follow from the fact that
its mirror image is equivalent to $(\sphere,\mathcal{T}_{3,3})$
as illustrated in Fig.\ \ref{fig:equivalence_to_m5_2} 
and Corollary \ref{cor:example_Tmn_irre_atoro} (see
also \cite[Table $3$]{IshKisMorSuz:12}, \cite[Figures $4c$ and $13$]{LeeLee:12}, and \cite[Table $2$]{BelPaoWan:20}).  
Thus we have the following corollary of Corollary \ref{cor:irreducibility:irreducible_HKA} and Lemma \ref{lm:atoroidality_general}.
\begin{corollary}
$\pairImn$ is irreducible for every $\mu,\nu$,
and is atoroidal if $\mu+\nu+3\neq 6$.
\end{corollary}  
A spine of $\pairImn$ with $\mu=\nu=0$ is shown in
Fig.\ \ref{fig:example_I_00}. 

\subsubsection{Examples in Fig.\ \ref{fig:intro:nonuniqueness}}\label{subsubsec:irre_atoro_of_non_uniqueness_examples}
\begin{figure}[t]
\begin{subfigure}{0.48\textwidth}
\centering
\def\svgwidth{.7\columnwidth}
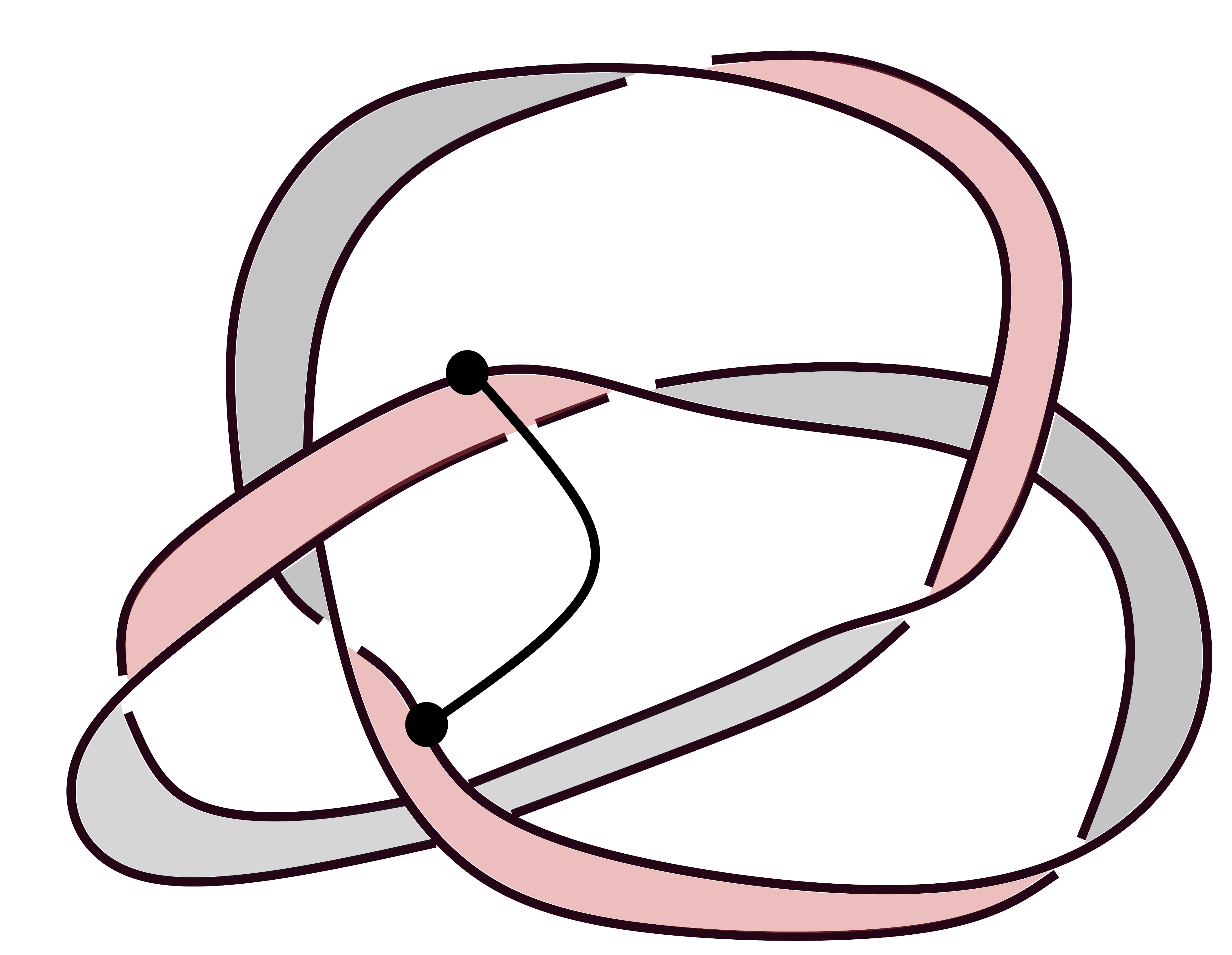
\caption{$A$ in $\ComplHKt$ and generators $x_1,x_2$.}
\label{fig:example_HKt_w_A_basis}
\end{subfigure}
\begin{subfigure}{0.48\textwidth}
\centering
\def\svgwidth{.7\columnwidth}
\begingroup%
  \makeatletter%
  \providecommand\color[2][]{%
    \errmessage{(Inkscape) Color is used for the text in Inkscape, but the package 'color.sty' is not loaded}%
    \renewcommand\color[2][]{}%
  }%
  \providecommand\transparent[1]{%
    \errmessage{(Inkscape) Transparency is used (non-zero) for the text in Inkscape, but the package 'transparent.sty' is not loaded}%
    \renewcommand\transparent[1]{}%
  }%
  \providecommand\rotatebox[2]{#2}%
  \newcommand*\fsize{\dimexpr\f@size pt\relax}%
  \newcommand*\lineheight[1]{\fontsize{\fsize}{#1\fsize}\selectfont}%
  \ifx\svgwidth\undefined%
    \setlength{\unitlength}{1133.85826772bp}%
    \ifx\svgscale\undefined%
      \relax%
    \else%
      \setlength{\unitlength}{\unitlength * \real{\svgscale}}%
    \fi%
  \else%
    \setlength{\unitlength}{\svgwidth}%
  \fi%
  \global\let\svgwidth\undefined%
  \global\let\svgscale\undefined%
  \makeatother%
  \begin{picture}(1,0.775)%
    \lineheight{1}%
    \setlength\tabcolsep{0pt}%
    \put(0,0){\includegraphics[width=\unitlength,page=1]{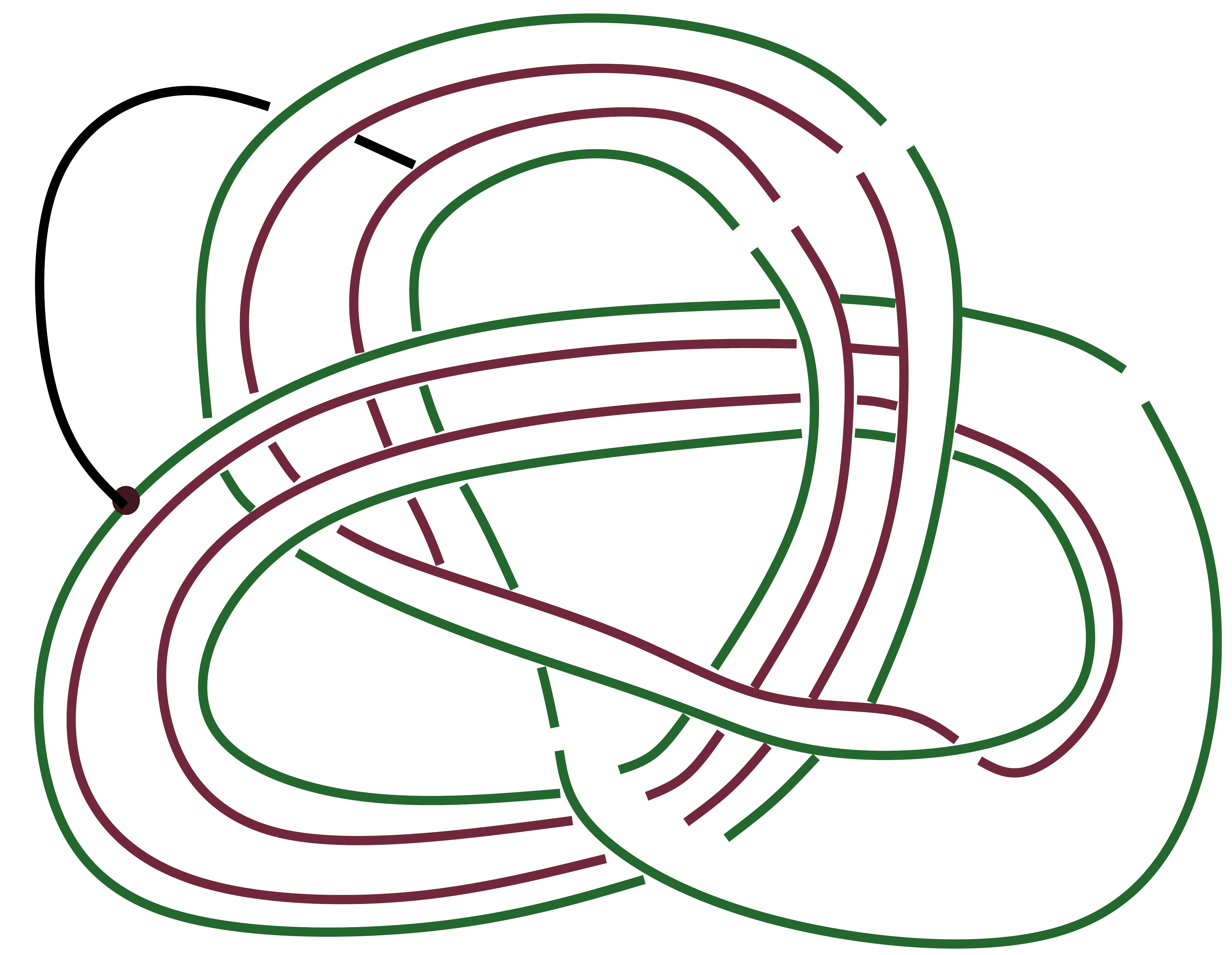}}%
    \put(0.68718842,0.7322929){\color[rgb]{0,0,0}\makebox(0,0)[lt]{\lineheight{1.25}\smash{\begin{tabular}[t]{l}{\footnotesize $\pairHKc$}\end{tabular}}}}%
    \put(0,0){\includegraphics[width=\unitlength,page=2]{example_HKc_w_A.pdf}}%
    \put(0.74721428,0.01924606){\color[rgb]{0,0,0}\makebox(0,0)[lt]{\lineheight{1.25}\smash{\begin{tabular}[t]{l}{\footnotesize $A$}\end{tabular}}}}%
    \put(0.29110493,0.51213984){\color[rgb]{0,0,0}\makebox(0,0)[lt]{\lineheight{1.25}\smash{\begin{tabular}[t]{l}{\footnotesize $A$}\end{tabular}}}}%
  \end{picture}%
\endgroup%

\caption{$A$ in $\ComplHKc$.}
\label{fig:example_HKc_w_A}
\end{subfigure}
\medskip

%
\begin{subfigure}{0.48\textwidth}
\centering
\def\svgwidth{.7\columnwidth}
\begingroup%
  \makeatletter%
  \providecommand\color[2][]{%
    \errmessage{(Inkscape) Color is used for the text in Inkscape, but the package 'color.sty' is not loaded}%
    \renewcommand\color[2][]{}%
  }%
  \providecommand\transparent[1]{%
    \errmessage{(Inkscape) Transparency is used (non-zero) for the text in Inkscape, but the package 'transparent.sty' is not loaded}%
    \renewcommand\transparent[1]{}%
  }%
  \providecommand\rotatebox[2]{#2}%
  \newcommand*\fsize{\dimexpr\f@size pt\relax}%
  \newcommand*\lineheight[1]{\fontsize{\fsize}{#1\fsize}\selectfont}%
  \ifx\svgwidth\undefined%
    \setlength{\unitlength}{1133.85826772bp}%
    \ifx\svgscale\undefined%
      \relax%
    \else%
      \setlength{\unitlength}{\unitlength * \real{\svgscale}}%
    \fi%
  \else%
    \setlength{\unitlength}{\svgwidth}%
  \fi%
  \global\let\svgwidth\undefined%
  \global\let\svgscale\undefined%
  \makeatother%
  \begin{picture}(1,0.775)%
    \lineheight{1}%
    \setlength\tabcolsep{0pt}%
    \put(0,0){\includegraphics[width=\unitlength,page=1]{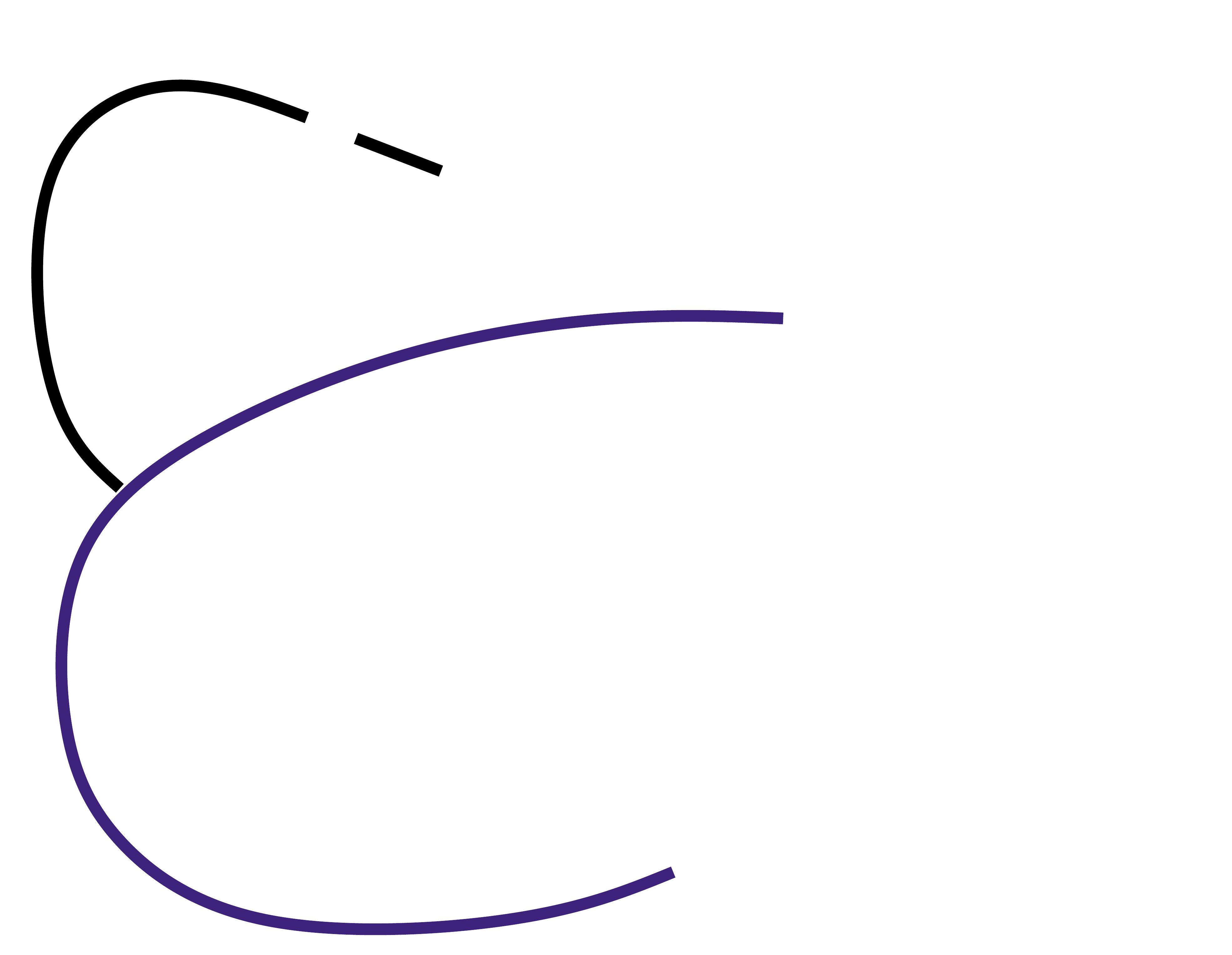}}%
    \put(0.68718842,0.7322929){\color[rgb]{0,0,0}\makebox(0,0)[lt]{\lineheight{1.25}\smash{\begin{tabular}[t]{l}{\footnotesize $\pairHKcA$}\end{tabular}}}}%
    \put(0,0){\includegraphics[width=\unitlength,page=2]{example_HKcA.pdf}}%
  \end{picture}%
\endgroup%

\caption{$\pairHKcA$.}
\label{fig:example_HKcA}
\end{subfigure}
\begin{subfigure}{0.48\textwidth}
\centering
\def\svgwidth{.7 \columnwidth}
\begingroup%
  \makeatletter%
  \providecommand\color[2][]{%
    \errmessage{(Inkscape) Color is used for the text in Inkscape, but the package 'color.sty' is not loaded}%
    \renewcommand\color[2][]{}%
  }%
  \providecommand\transparent[1]{%
    \errmessage{(Inkscape) Transparency is used (non-zero) for the text in Inkscape, but the package 'transparent.sty' is not loaded}%
    \renewcommand\transparent[1]{}%
  }%
  \providecommand\rotatebox[2]{#2}%
  \newcommand*\fsize{\dimexpr\f@size pt\relax}%
  \newcommand*\lineheight[1]{\fontsize{\fsize}{#1\fsize}\selectfont}%
  \ifx\svgwidth\undefined%
    \setlength{\unitlength}{1133.85826772bp}%
    \ifx\svgscale\undefined%
      \relax%
    \else%
      \setlength{\unitlength}{\unitlength * \real{\svgscale}}%
    \fi%
  \else%
    \setlength{\unitlength}{\svgwidth}%
  \fi%
  \global\let\svgwidth\undefined%
  \global\let\svgscale\undefined%
  \makeatother%
  \begin{picture}(1,0.775)%
    \lineheight{1}%
    \setlength\tabcolsep{0pt}%
    \put(0.68718842,0.7322929){\color[rgb]{0,0,0}\makebox(0,0)[lt]{\lineheight{1.25}\smash{\begin{tabular}[t]{l}{\footnotesize $\pairHKcA$}\end{tabular}}}}%
    \put(0,0){\includegraphics[width=\unitlength,page=1]{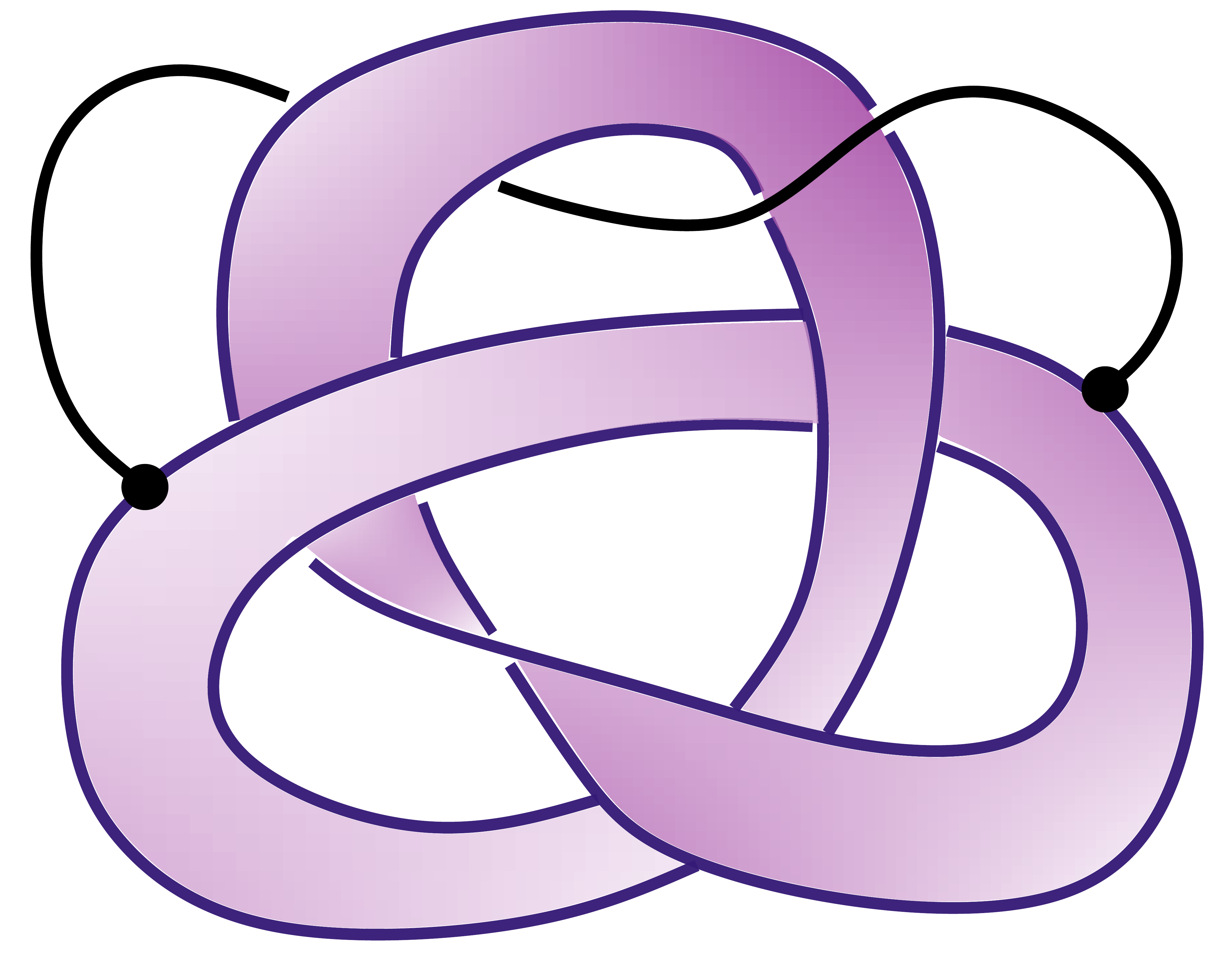}}%
    \put(0.22719584,0.60156425){\color[rgb]{0,0,0}\makebox(0,0)[lt]{\lineheight{1.25}\smash{\begin{tabular}[t]{l}{\footnotesize $M$}\end{tabular}}}}%
    \put(0.65083733,0.09173355){\color[rgb]{0,0,0}\makebox(0,0)[lt]{\lineheight{1.25}\smash{\begin{tabular}[t]{l}{\footnotesize $M$}\end{tabular}}}}%
  \end{picture}%
\endgroup%

\caption{M\"obius band $M\subset \ComplHKcA$.}
\label{fig:example_HKcA_w_M}
\end{subfigure}
\caption{Annuli and M\"obius band in handlebody-knot exteriors.}
\end{figure}

In Section \ref{sec:intro}, two handlebody-knots
are given in Fig.\ \ref{fig:intro:nonuniqueness} 
to show how torus and cable links are used to
construct irreducible, atoroidal 
handlebody-knots whose exteriors contain
non-isotopic type $3$-$3$ annuli.
We now verify their irreducibility and atoroidality.
\begin{corollary}
$\pairHKt$ in Fig.\ \ref{fig:intro:example_two_annuli_torus_link}
is irreducible and atoroidal.
\end{corollary}
\begin{proof}
Let $A$ be the oriented type $3$-$3$ annulus in $\ComplHKt$
indicated in Fig.\ \ref{fig:example_HKt_w_A_basis}.  
Then note that first $\pairHKtA$ is trivial,
and secondly, in terms of the meridional basis 
of $H_1(\ComplHKtA)$ given by $x_1,x_2$ 
in Fig.\ \ref{fig:example_HKt_w_A_basis},
$[l_+]=(5,1)$ and $[l_-]=(4,2)$.
 
Since $\partial A$ is a $(3,2)$-torus link with $3\cdot 2=6$
and $(5,1)\in H_1(\ComplHKtA)$ is a generator, 
the third criterion of
Lemma \ref{lm:irreducibility:trivial_HKA} implies 
that $\pairHKt$ is irreducible. As $A$ is unknotting, 
$\pairHKt$ is atoroidal by Lemma \ref{lm:atoroidality:trivial_HKA}. 
\end{proof}

\begin{corollary}
$\pairHKc$ in Fig.\ \ref{fig:intro:example_two_annuli_torus_link}
is irreducible and atoroidal.
\end{corollary}
\begin{proof}
Let $A\subset \ComplHKc$ be the type $3$-$3$ annulus 
in Fig.\ \ref{fig:example_HKc_w_A}, and observe that 
$\pairHKcA$ is the handlebody-knot
given in Fig.\ \ref{fig:example_HKcA}. 
There exists an essential M\"obius band $M$ in $\ComplHKcA$
as shown in Fig.\ \ref{fig:example_HKcA_w_M}. 
The M\"obius band $M\subset \ComplHKcA$ induces 
the handlebody-knot
$\pairHKcAM$, which is equivalent to $(\sphere,5_2)$ 
in \cite{IshKisMorSuz:12} 
and hence equivalent to 
the mirror image of $(\sphere,\mathcal{T}_{3,3})$ (see Fig.\ \ref{fig:equivalence_to_m5_2}). In particular, $\pairHKcAM$
is irreducible and atoroidal by Corollary \ref{cor:example_Tmn_irre_atoro}.

Now, by Remark \ref{rmk:irreducibility_other_types},
$\pairHKcA$ is irreducible, and since 
$\partial M$ is a $(5,2)$-cable of a trefoil, 
$\pairHKcA$ is atoroidal by Corollary \ref{cor:atoroidality_mobius_band}.
Applying Corollary \ref{cor:irreducibility:irreducible_HKA}
and Lemma \ref{lm:atoroidality_general} 
to $\pairHKc$ and $A$, 
we obtain that $\pairHKc$ is irreducible and atoroidal as 
$\partial A$ is a $(10,4)$-cable of a torus knot.
\end{proof}


\section{Uniqueness}\label{sec:uniqueness}
Throughout the section,  
$\pair$ is assumed to be irreducible, atoroidal,
unless otherwise specified, and its exterior admits 
a type $3$-$3$ annulus $A$ with a non-trivial boundary slope of $p$. 
We fix an orientation of $A$, and let $l_1,l_2, l_+\subset A_+,l_-\subset A_-$ be as in Section \ref{sec:disks}, and
denote by $l_1^+, l_2^+$ (resp.\ $l_1^-, l_2^-$) 
the components of $\partial A_+$ (resp.\ $\partial A_-$).

Note that the intersection $\rnbhd{A}\cap \partial\HK$
is a regular neighborhood $\rnbhd{l_1\cup l_2}$ 
of $l_1\cup l_2$ in $\partial \HK$;
we denote by $\rnbhd{l_i}$ the component of $\rnbhd{l_1\cup l_2}$
containing $l_i,i=1,2$.
Recall from Section \ref{subsec:typethreethree} 
that the annulus $A$ determines a meridian-disk system
$\dsystem_A:=\{\disk_A,\disk_1,\disk_2\}$ of $\HK$,
where $\disk_A$ is a separating meridian disk
disjoint from $l_1,l_2$, and $\disk_i$ is a meridian disk of the solid torus in
$\HK-\rnbhd{D_A}$ containing $l_i$, $i=1,2$. 
$\dsystem_A$ induces a handcuff spine $\Gamma_A$ 
of $\HK$. 

\subsection{Uniqueness criteria}
\begin{lemma}\label{lm:parallel}
Let $A'$ be another type $3$-$3$ annulus in $\ComplHKA$ 
disjoint from $A$. 
Suppose one of the following holds:
\begin{itemize}
\item $\vert p\vert=1$; 
\item $\vert p\vert>1$, $A$ satisfies the condition \eqref{cond:torus_cable}, and 
neither of $l_+,l_-$ represents the $\absp$-th power of 
some element in $\pi_1(\ComplHKA)$; 
\item $\vert p\vert>1$, $\pairA$ is trivial,
$A$ satisfies the condition \eqref{cond:torus}, and 
neither of $l_+,l_-$ represents the $\absp$-th power of 
some primitive element in $\pi_1(\ComplHKA)$   
\end{itemize}
Then $A,A'$ are parallel in $\Compl\HK$.
\end{lemma}
\begin{proof}
Without loss of generality, it may be assumed that
$A'\cap \rnbhd{A}=\emptyset$.
Denote by $l_i',i=1,2$, the components of $\partial A'$.
Then $l'_i, i=1,2$, are in the surface 
\[S_{0,4}:=\partial \HK-\openrnbhd{l_1\cup l_2}.\]
$S_{0,4}$ is homeomorphic to a $2$-sphere with $4$ open disks removed
and $\partial S_{0,4}= l_1^+\cup l_1^-\cup l_2^+\cup l_2^-$. 
Since $l_1'$ is not parallel to $l_2'$ in $\partial \HK$, 
one of $l_i',i=1,2,$ is parallel to one of 
$l_1^\pm,l_2^\pm$ in $S_{0,4}$. Without loss of generality, it may be
assumed that $l_1'$ is parallel to $l_1^+$, denoted by 
$l_1'\varparallel l_1^+$, 
and hence $l_2'\nvarparallel l_1^-$ in $S_{0,4}$.

Now, $\{\partial \disk_1,\partial \disk_2\}$ induces 
a basis $\{u_1,u_2\}$ of $H_1(\Compl\HK)$ in terms of which  
$[l_i]=(p,p)$, $i=1,2$. Orient $A',l_1',l_2'$
so that $\partial A'=l_1'\cup -l_2'$.
Suppose $l_2'\nvarparallel l_2^+$ or $l_2\nvarparallel l_2^-$ in $S_{0,4}$.
Then $[l_2']=(\pm 2p,\pm 2p)$ or $(0,0)$ in $H_1(\ComplHK)$ in terms of 
the basis $\{u_1,u_2\}$, contradicting 
the fact $[l_1']=[l_2']\in H_1(\Compl{\HK})$.
Suppose $l_2'\varparallel l_2^-$; isotope $A'$ so that $l_1'=l_1^+$ and $l_2'=l_2^-$. 
Then the union $A_+\cup A'\cup\rnbhd{l_2}$ is a torus $T$
separating the two loops of $\Gamma_A$ but disjoint from $\Gamma_A$,
contradicting the connectedness of $\Gamma_A$. 
 
As a result, the only possible case is 
$l_2'\varparallel l_2^+$. Isotope $A'$ so that
$\partial A_+=\partial A'$. Then the union $T:=A_+\cup A'$ 
is a torus in $\Compl\HK$, which, by the atoroidality
of $\pair$, bounds a solid torus $W\subset\Compl\HK$ disjoint from $A$.
Since $\lk{l_1}{l_2}=p$,
and $A_+\subset \partial W\subset W$, 
if the slops of an essential loop $l'$ of $A'$ 
in $W$ is $\frac{m}{n}$, then $mn=p$.

Suppose $\vert p \vert=1$.
Then $\pi_1(A')\rightarrow  \pi_1(W)$
is an isomorphism, and thus
$A'$ is parallel through $W$ to $A_+$ and hence parallel to $A$.

Suppose $\vert p\vert >1$.
Consider first the case where $(\sphere,W)$
is trivial, namely $\Compl{W}$ is a solid torus.
By the condition \eqref{cond:torus_cable} or \eqref{cond:torus},   
$l_+\subset \sphere$ is not an $(m,n)$-torus knot, and therefore
the slope of $l'$ is either $p$ or $\frac{1}{p}$.
The former implies $A,A'$ are parallel through $W$.
For the latter, we 
consider the complement $U:=\ComplHKA-\mathring{W}$; note that 
$U\cap W=A'$. Since 
$A'$ is incompressible in $\ComplHK$,
it is incompressible in $\ComplHKA$.
In particular, the homomorphisms 
\[
\pi_1(A')\rightarrow \pi_1(U)
\quad\text{and}\quad
\pi_1(A')\xrightarrow{\phi}  \pi_1(W)
\]
induced by inclusions are injective with $\phi$ sending the generator 
of $\pi_1(A')$ to the $p$-th power of a generator of 
$\pi_1(W)$. By the van Kampen theorem, the homomorphism 
\[\pi_1(W)\xrightarrow{\psi} \pi_1(\ComplHKA)\]
induced by $W\subset\ComplHKA$
is also an injection. 
This implies that $l'$
and hence $l_+$ represent the $\vert p\vert$-th power of some element
in $\pi_1(\ComplHKA)$, contradicting the 
assumption. 

If $\pairA$ is trivial, then 
$A'$ is a incompressible, separating, 
non-boundary-parallel annulus in 
the handlebody $\ComplHKA$. 
By Corollary \ref{cor:incompressible_non_b_parallel_primitive_element},
the image of the generator of 
$\pi_1(W)$ under $\psi$ is a primitive element of 
the free group $\pi_1(\ComplHKA)$. 
Thus, $l'$
and hence $l_+$ represent the $\vert p\vert$-th power of some primitive element
in $\pi_1(\ComplHKA)$, contradicting 
the third criterion.

Consider now the case $(\sphere,W)$ is non-trivial. 
The condition \eqref{cond:torus_cable} implies that
the slope of $l'$ can only be $p$; thus $A,A'$ are parallel.
On the other hand, if $\pairA$ is trivial, 
the cable knot condition in \eqref{cond:torus_cable} can be dropped.
To see this, we note that, since $\pairA$ is trivial,
$l_+\subset \sphere$ 
is a tunnel number one satellite knot. By the classification
theorem of tunnel number one satellite knots in \cite{MorSak:91}, 
the slope of $l_+$ and hence $l'$ on $W$ can only be $p$ or $\frac{1}{p}$.
The latter, as in the previous case, 
implies that $l_+$
represents the $p$-th power of some primitive element of $\pi_1(\ComplHKA)$,
contradicting the third criterion.
Therefore the slope of $l'$ can only be $p$, and hence $A,A'$ are parallel.
\end{proof}

\begin{theorem}\label{teo:uniqueness_HKA_irreducible}
Suppose $\pairA$ is irreducible, and 
$A$ satisfies the condition \eqref{cond:torus_cable}.
If one of the following conditions holds:
\begin{itemize}
\item $\vert p\vert=1$;  
\item $\vert p\vert>1$, and 
none of $l_+,l_-$ represents the $\vert p\vert$-th power of 
some element in $\pi_1(\ComplHKA)$.   
\end{itemize}
then, up to isotopy, $A$ is the unique type $3$-$3$
annulus in $\ComplHK$.
\end{theorem}
\begin{proof}
Suppose $A'$ is another type $3$-$3$ annulus in $\Compl\HK$. 
Isotope $A$ such that   
$\# A'\cap A$ is minimized.
If $\# A'\cap A=\emptyset$, the assertion follows from Lemma \ref{lm:parallel}, so we assume $\# A'\cap A\neq\emptyset$. 
By the essentiality of $A,A'$ and ($\partial$-)irreducibility of $\Compl\HK$, any arc or circle in $A'\cap A$ is essential in both $A$ and $A'$,
and therefore $A\cap A'$ are either some circles or some arcs.

Suppose $A\cap A'$ are some circles.
It may be assumed that $\partial A'\cap \rnbhd{A}=\emptyset$,
and as argued in the proof of Lemma \ref{lm:parallel}, since 
the components $l_1',l_2'$ of $\partial A'$ are not parallel in 
\[S_{0,4}:=\partial \HK-\openrnbhd{l_1\cup l_2},\]
and represent the same element in $H_1(\Compl\HK)$, up to sign,
$l_1',l_2'$ are parallel to some components of 
$\partial S_{0,4}$. It may be assumed that
$l_1'$ is parallel to $l_1^+$, and hence 
$l_2'$ is parallel to either $l_2^+$ or $l_2^-$.

Let $\rho\subset A'$ be an outermost circle in $A'$ and
$P'$ the outermost annulus cut off by $\rho$ from $A'$.
Without loss of generality, it may be assumed that
$\partial P'=\rho\cup l_1'$. Let $P\subset A$ be the annular
component cut off by $\rho$ with $\partial P=\rho\cup l_2$.
Then $P\cup P'$ induces a type $3$-$3$ annulus $A''$,
which is disjoint from $A$, and therefore isotopic to $A$
by Lemma \ref{lm:parallel}. 
$A''$ however has less intersection
with $A'$ than $A$ does, contradicting the minimality.

Suppose $A\cap A'$ are some arcs.
It may be assumed that 
$A'\cap\rnbhd{A}$ is a regular neighborhood of $A\cap A'$.
Thus, $A'-\openrnbhd{A}$ consists of some disks   
in $\ComplHKA$, 
each of which meets $l_+\cup l_-$ at two points.
Let $D$ be a disk in $A'-\openrnbhd{A}$.
Then, by the irreducibility of $\pairA$,
$D\subset \ComplHKA$ is inessential 
and $\partial D$ bounds a disk $E\subset\partial\HK_A$,
which cobuonds a $3$-ball $B$ with $D$.
Note that this implies that $D$ meets either $l_+$ or $l_-$.
Isotoping $A$ through $B$ decreases $\# A'\cap A$, and 
contradicts the minimality.  
Therefore $A\cap A'=\emptyset$. 
 
\end{proof}

\begin{theorem}\label{teo:uniqueness_HKA_trivial_not_primitive}
Suppose $\pairA$ is trivial, and $A$ satisfies the condition \eqref{cond:torus}.
If one of the following holds:
\begin{enumerate}
\item $\vert p\vert=1$;  
\item $\vert p\vert>1$, 
at least one of $l_+,l_-\subset\ComplHKA$ is not primitive, and neither of $l_+,l_-$ represents the $\vert p\vert$-th power of some primitive element in $\pi_1(\ComplHKA)$,
\end{enumerate}
then, up to isotopy, $A$ is the unique type $3$-$3$
annulus in $\ComplHK$.
\end{theorem}
\begin{proof}
Let $A'$ be another type $3$-$3$ annulus in $\Compl\HK$, and 
isotope $A$ such that   
$\# A'\cap A$ is minimized. 
It suffices to show that $A'\cap A=\emptyset$
in view of Lemma \ref{lm:parallel}.

Suppose $A'\cap A\neq \emptyset$.
The same argument in the proof of Theorem \ref{teo:uniqueness_HKA_irreducible} implies that $A'\cap A$
only contains arcs essential in both $A'$ and $A$;
it may be assumed that
$A'\cap\rnbhd{A}$ is a regular neighborhood of $A'\cap A$ in $A'$; 
thus by the minimality, $A'-\openrnbhd{A}$ 
consists of essential disks   
in $\ComplHKA$, each meeting $l_+\cup l_-$ at two points.
Let $D$ be one of the disks.
Properly orient $D$. Then $I_D$ is one of the following:      
\begin{equation}\label{eq:five_cases}
(1, 1),(1,-1),(2,0),(0,2), (0,0).
\end{equation}
 
Any of the 
first four cases in \eqref{eq:five_cases} implies that
$D$ cuts $l_+\cup l_-$ into two arcs $\alpha,\beta$. 
The boundary $D'$ of a regular 
neighborhood $\rnbhd{D\cup \alpha}$
is an essential disk which meets either $l_+$ or $l_-$
at two points with $I_{D'}=(0,0)$, but
this contradicts Lemma \ref{lm:trivial_intersection_not_primitive}.
Thus $I_D=(0,0)$ with $D\subset\ComplHKA$ inessential 
in the only possibility, but in this case, 
one can isotope $A$ via the $3$-ball bounded by $D$
and the disk in $\partial\HK_A$ bounded by $\partial D$  
to decrease $\# A'\cap A$, contradicting the minimality. 
\end{proof}
 
The following simplified version of Theorems \ref{teo:uniqueness_HKA_irreducible} and \ref{teo:uniqueness_HKA_trivial_not_primitive}
is sufficient for many applications.
\begin{corollary}\label{cor:uniqueness}
Suppose $A$ satisfies the condition \eqref{cond:torus_cable}, and 
none of $l_+,l_-$
represents the $\absp$-th multiple of some element in $H_1(\ComplHKA)$,
and if $\pairA$ is trivial, at least one of $l_+,l_-\subset\ComplHKA$ is not primitive. 
Then, up to isotopy, $A$ is the unique type $3$-$3$
annulus in $\ComplHK$.
\end{corollary}

Note that Corollary \ref{cor:uniqueness} 
fails to include the case $\vert p\vert\leq 2$. 
On the other hand, when $p=\pm 1$,
the condition \eqref{cond:torus_cable} 
is automatically satisfied, and the existence
of such a type $3$-$3$ annulus $A$ turns out to impose strong constraints
on $\pair$ as well as on $A$ itself.
Theorems \ref{teo:uniqueness_HKA_irreducible} and \ref{teo:uniqueness_HKA_trivial_not_primitive}, 
along with 
Corollary \ref{cor:irreducibility:irreducible_HKA} and 
Lemmas \ref{lm:irreducibility:trivial_HKA} and \ref{lm:atoroidality_general}, imply the following.

\begin{corollary}\label{cor:rigidity_slope_pm1}
Given a handlbody-knot $\pair$ 
and a type $3$-$3$ annulus $A\subset\ComplHK$ with
a boundary slope of $\pm 1$. 
Suppose $\pairA$ is atoroidal, and if $\pairA$ is trivial, $\{l_+,l_-\}$
does not represent a basis of $\pi_1(\ComplHKA)$.
Then $\pair$ is irreducible, atoroidal, and $A$ is 
the unique type $3$-$3$
annulus in $\ComplHK$, up to isotopy.
\end{corollary}
 
In view of Remark \ref{rmk:lpm_primtive_trivial_HK}, 
Corollary \ref{cor:rigidity_slope_pm1} still holds
if the condition
of $\{l_+,l_-\}$
not representing a basis of $\pi_1(\ComplHKA)$
is replaced with $\pair$ being non-trivial.
In the event that $p\neq \pm 1$ and both $l_+,l_-\subset\ComplHKA$
are primitive, the next criterion comes in handy.

\begin{theorem}\label{teo:uniqueness_HKA_trivial_odd_p}
Suppose $\pairA$ is trivial and $A$ satisfies
the condition \eqref{cond:torus}.
If $\vert p\vert>1$, $p$ is odd, 
and neither of $l_+,l_-$ represents the $\absp$-th multiple of some generator 
of $H_1(\ComplHKA)$, 
then, up to isotopy, $A$ is the unique type $3$-$3$
annulus in $\ComplHK$.   
\end{theorem}
\begin{proof}
Suppose $A'$ is another type $3$-$3$ annulus, 
and $\# A'\cap A$ is minimized in the isotopy classes of $A,A'$. 
By Lemma \ref{lm:parallel}, 
it suffices to consider the case where $A'\cap A\neq\emptyset$. 
As before, we may assume $A'\cap\rnbhd{A}$ is 
a regular neighborhood of $A'\cap A$ in $A'$. 
Then every disk component $D$ of 
$A'-\openrnbhd{A}$ is essential in $\ComplHKA$ 
by the minimality, and
$I_D$ is one of five cases in \eqref{eq:five_cases}. 

Apply Lemma \ref{lm:non-separating_two_intersections}
to rule out the first, third and fourth cases in \eqref{eq:five_cases}. 
Then observe that, 
since $\# A_+\cap A'= \# A_-\cap A'$,
there exists a disk $D\subset A'-\openrnbhd{A}$
with $I_D=(0,0)$ and $D\cap l_+=\emptyset$
if and only if 
there exists a disk $D' \subset A'-\openrnbhd{A}$
with $I_{D'}=(0,0)$, $D' \cap l_-= \emptyset$.
By Lemma \ref{lm:trivial_intersection}, $D',D$ are separating essential disks; since $D,D'$ are disjoint,
they are parallel. 
Thus, one can isotope $D,D'$ 
away from $l_+\cup l_-$,
contradicting Corollary \ref{cor:separating_disjoint}.
%


Consider now the remaining case: $I_D=(1,-1)$.
Suppose $A'-\openrnbhd{A}$ consists of $n$ disks $D_1,\dots, D_n$.
Label the arcs in $A\cap A'\subset A$ consecutively from
$\alpha_1$ to $\alpha_n$, that is, $\alpha_i,\alpha_j$
cutting off a disk $E\subset A$ with $E\cap A'=\alpha_i\cup \alpha_j$ 
wherever $j\equiv i+1$ (mod $n$).
Since $I_{D_i}=(1,-1)$, $i=1,\dots, n$, 
the disks $D_i\subset A',i=1,\dots,n$, 
induce a permutation $\sigma$ on $\{1,\dots,n\}$
defined as follows: $\sigma(r)=s$ if there exists 
a disk $D\in\{D_1,\dots, D_n\}$   
with $\alpha_r^+=D \cap A_+$ and $\alpha_s^-=D\cap A_-$,
where $\alpha_i^\pm\subset A_\pm\cap A'$ are the arcs 
corresponding to $\alpha_i$, namely, $\alpha_i^\pm$ parallel to $\alpha_i$
in $A'\cap \rnbhd{A}$ (Fig.\ \ref{fig:permutation}).
Because $A'$ is connected, $\sigma$ is of order $n$. 
On the other hand, 
by Lemma \ref{lm:plus_minus_one_parallel}, 
$D_i,i=1,\dots,n$,
are parallel in $\ComplHKA$ (Fig.\ \ref{fig:parallel_disks}), and
therefore if $\sigma(1)=k$, then
\begin{align*}
\sigma(i)&=k-i+1, \quad 1\leq i\leq k\\
\sigma(i)&=n-i+k+1, \quad k+1\leq i\leq n. 
\end{align*}
Particularly, we have $\sigma^2=\id$, and hence $n=2$.

\begin{figure}[t] 
\begin{subfigure}{0.49\textwidth}
\centering
\def\svgwidth{.9\columnwidth}
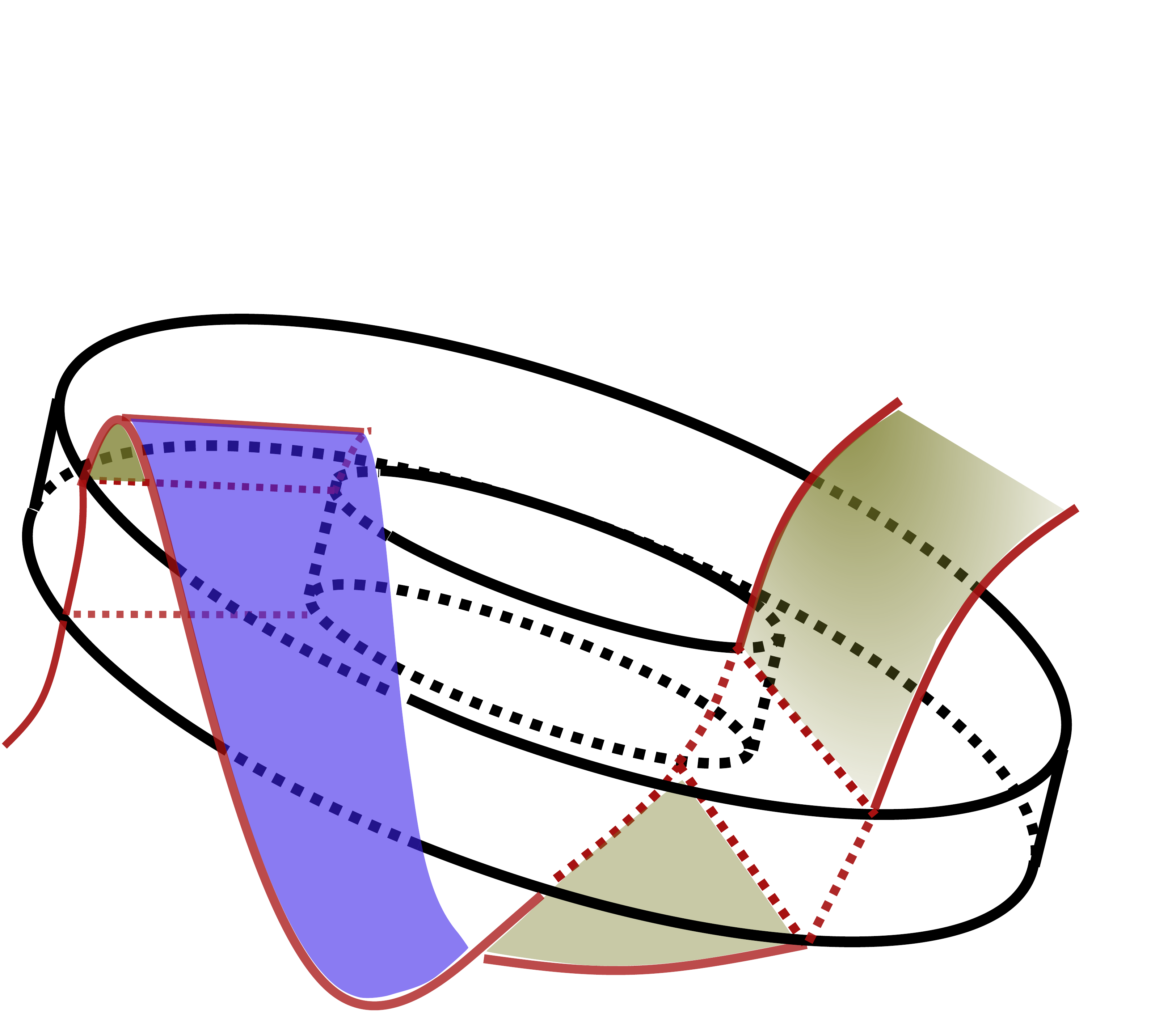
\caption{The induced permutation $\sigma$.}
\label{fig:permutation} 
\end{subfigure}
\begin{subfigure}{0.49\textwidth}
\centering
\def\svgwidth{.95\columnwidth}
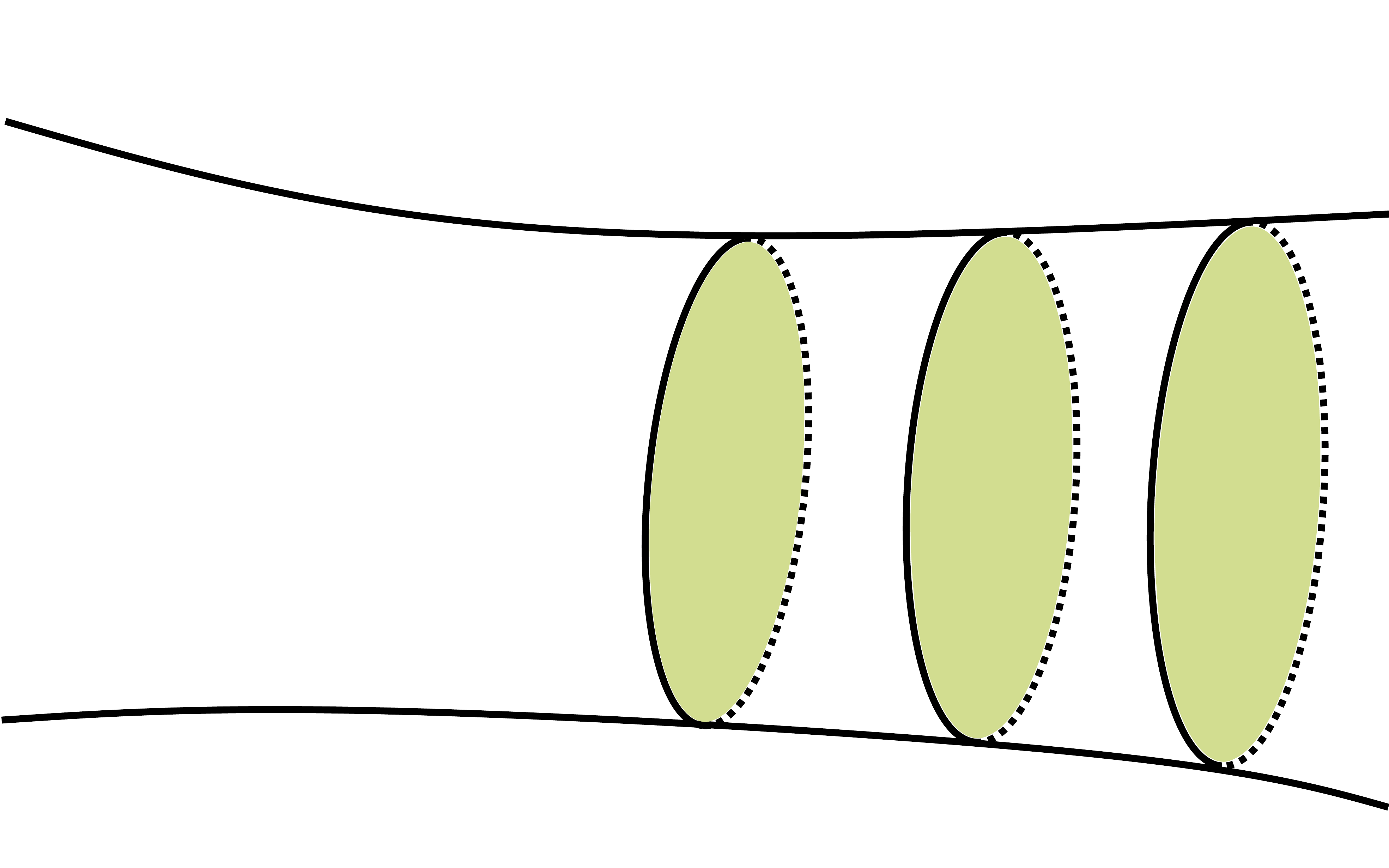
\caption{Parallel $D_1,\dots, D_n\subset\ComplHKA$.}
\label{fig:parallel_disks} 
\end{subfigure}
\begin{subfigure}{0.4\textwidth}
\centering
\def\svgwidth{.95\columnwidth}
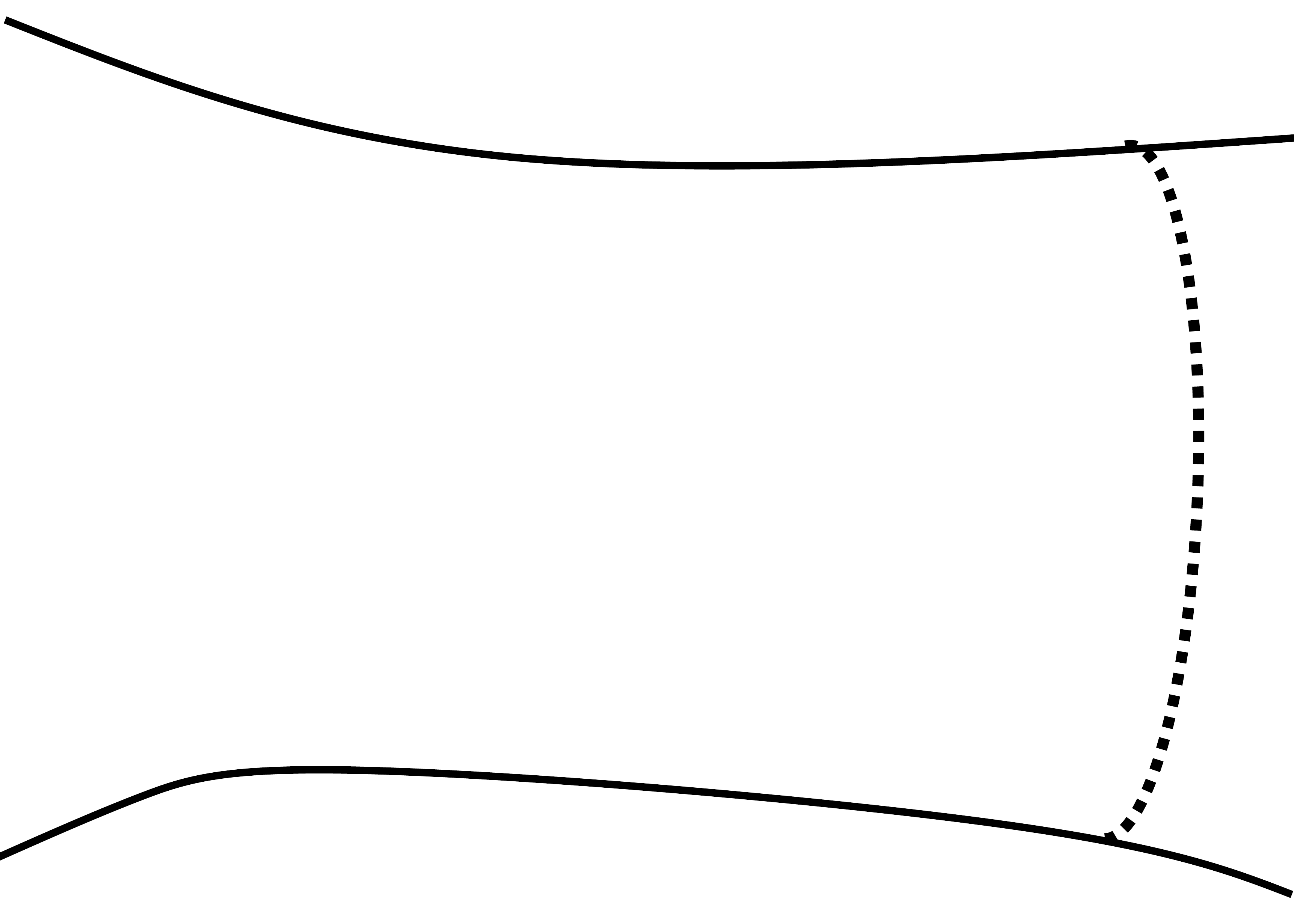
\caption{$P_1,P_2\subset \overline{\partial\ComplHKA-A_+\cup A_-}$.}
\label{fig:cylinder_P1_P2} 
\end{subfigure}
\begin{subfigure}{0.49\textwidth}
\centering
\def\svgwidth{.9\columnwidth}
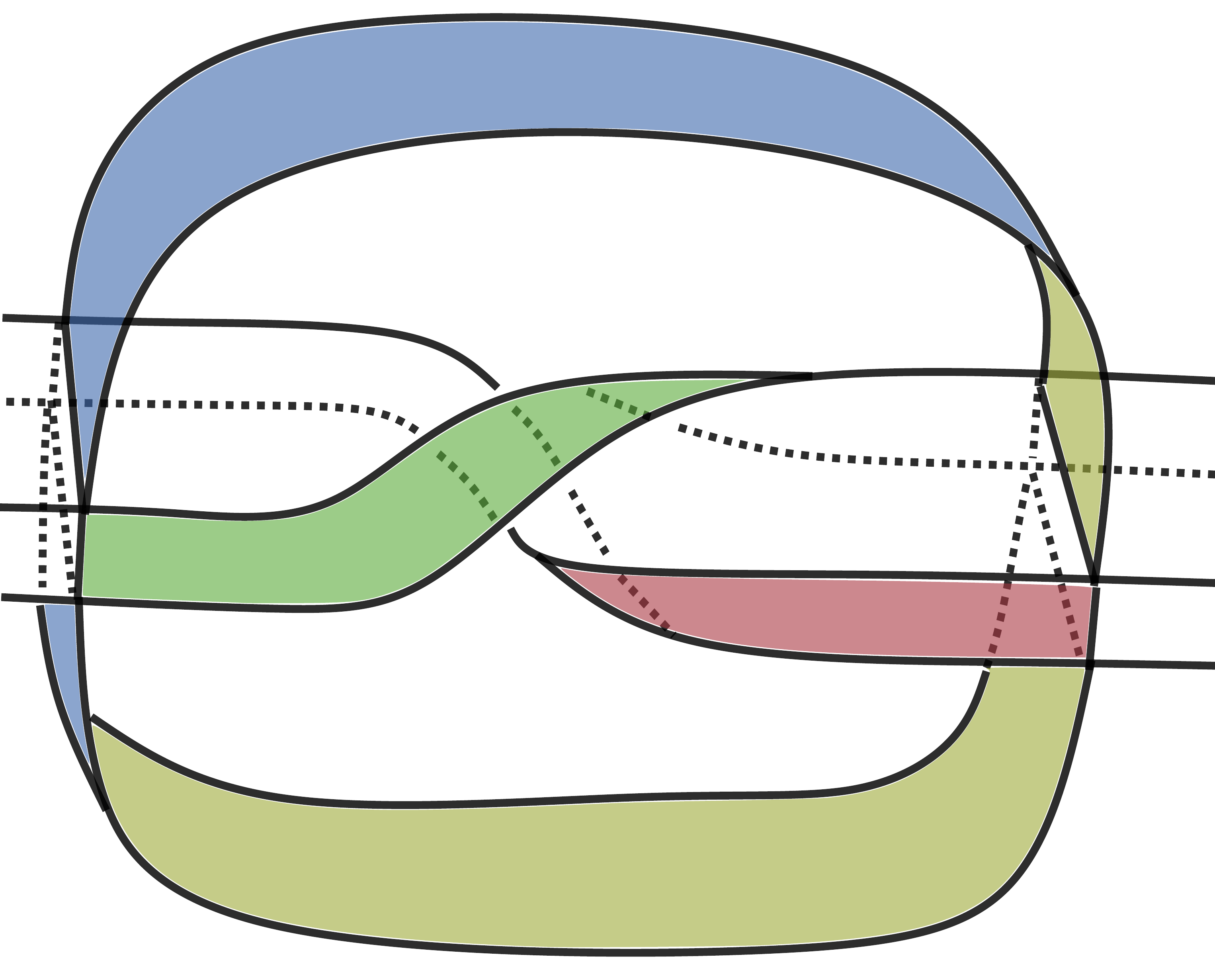
\caption{$Q_1,Q_2\subset \overline{\partial \rnbhd{A}-A_+\cup A_-}$.}
\label{fig:cylinder_Q1_Q2} 
\end{subfigure}
\caption{The case $I_D=(1,-1)$.}
\end{figure}

Since $D_1,D_2$ are parallel, $\partial (D_1\cup D_2)$
cuts off a cylinder $C_0$ from $\partial\ComplHKA$.
The closure of $C_0-(A_+\cup A_-)$
consists of two disks $P_1,P_2$ (Fig.\ \ref{fig:cylinder_P1_P2}).
At the same time, 
$\partial A'$ cuts the annulus 
$\rnbhd{l_1}$ (resp.\ $\rnbhd{l_2}$)
into two disks, 
one of which, denoted by $Q_1$ (resp.\ $Q_2$),
meets each of $P_1, P_2$
at an arc---one in $l_1^+$ (resp.\ $l_2^-$) and the other in $l_1^-$ (resp.\ $l_2^+$); see Fig.\ \ref{fig:cylinder_Q1_Q2} for 
an illustration.  
The union $P_1\cup Q_1\cup P_2\cup Q_2$
is a cylinder $C$ in $\partial\HK$ 
with $\partial C=\partial A'$, contradicting that
$A'$ is of type $3$-$3$.
Therefore $A'\cap A=\emptyset$. 
\end{proof}

\subsection{Examples}\label{subsec:examples_uniqueness}
Consider first the handlebody-knot family 
$\pairTmn$ in Section \ref{subsec:examples_irreatoro}, 
where $\mu,\nu$ are odd integers; we assume $\frac{\mu+\nu}{2}\neq 0$.
Let $A\subset \ComplTmn$ be
the canonical type $3$-$3$ annulus with a boundary slope of $\frac{\mu+\nu}{2}$ 
given by the construction. Orient $A$  
as in Fig.\ \ref{fig:annulus_tunnel_trivial_knot}. Then, 
in terms of the meridional basis  
of $H_1(\ComplTmnA)$ induced by 
$x_1,x_2$ in Fig.\ \ref{fig:tunnel_trivial_knot},  
$[l_+]=(\frac{\mu+1}{2}, \frac{\nu-1}{2} )$,  
$[l_-]=(\frac{\mu-1}{2}, \frac{\nu+1}{2} )$, and hence, 
Theorem \ref{teo:uniqueness_HKA_trivial_odd_p} implies the following. 
\begin{corollary}\label{cor:example_T_mn_unique_A_p_odd}
If $\frac{\mu+\nu}{2}$ is odd, and not equal to $\pm 1$,
and $\frac{\mu\pm 1}{2},\frac{\nu\pm 1}{2}$
are not divisible by $p$,
then $A\subset\ComplTmn$ is the unique type $3$-$3$ annulus, up to isotopy. 
\end{corollary}
To produce examples with even boundary slope,
we apply Corollary \ref{cor:uniqueness}.
\begin{corollary}\label{cor:example_T_mn_unique_A_non_primitive}
If at most one of $\frac{\mu+ 1}{2},\frac{\mu- 1}{2}, 
\frac{\nu+1}{2}, \frac{\nu- 1}{2}$ equals to $1$ or $-1$, and 
none of them is divisible by $p$, 
then $A\subset\ComplTmn$ is the unique type $3$-$3$ annulus, up to isotopy. 
\end{corollary}
\begin{proof} 
The first condition ensures that one of $l_+,l_-\subset \ComplTmnA$ is not a primitive loop by \cite{Zie:70},
\cite{CohMetZim:81}, \cite{OsbZie:81}.
\end{proof}
Criterion in Corollary \ref{cor:example_T_mn_unique_A_non_primitive},
though covers many cases, does fail to include some
important small crossing handlebody-knots, such as $(\sphere, 5_3),
(\sphere,6_4)$ in the handlebody-knot table \cite{IshKisMorSuz:12}.
The former, being equivalent to the mirror image of 
$(\sphere,\mathcal{T}_{3,3})$ (Fig.\ \ref{fig:equivalence_to_m5_2}), is covered in Corollary \ref{cor:example_T_mn_unique_A_p_odd}, 
whereas the latter,
equivalent to the mirror image of $(\sphere,\mathcal{T}_{-3,5})$ (Fig.\ \ref{fig:equivalence_to_m6_4}),
is considered in Section \ref{subsec:homotopy_criteria}.
 
Corollaries \ref{cor:example_T_mn_unique_A_p_odd}
and \ref{cor:example_T_mn_unique_A_non_primitive},
together with the slope type of $A$,
give us the following. 
 
\begin{corollary}
The family of handlebody-knots
\begin{equation}\label{eq:family_Tmn}
\{\pairTmn\mid \mu,\nu \text{ odd integers }\}
\end{equation}
contains infinitely many irreducible,
atoroidal handlebody-knots whose exterior
admit a unique type $3$-$3$ annulus, which is unknotting.
\end{corollary}
\begin{proof}
Consider the subfamily of \eqref{eq:family_Tmn}
\[\mathcal{P}_T:=\{\pairTmn\mid  \mu\geq \nu>  1 \text{ or } -1>\mu\geq \nu\}.\]
The condition $\mu\geq \nu>  1$ or $-1>\mu\geq \nu$
implies that the meridional basis of $H_1(\ComplHKA)$
induced by $x_1,x_2$ is normalized, and 
the orientation of $A$ given in Fig.\ \ref{fig:annulus_tunnel_trivial_knot}
is a preferred one. Thus, $(\frac{\mu+1}{2},\frac{\nu-1}{2})$ is the slope type of $A$.
On the other hand, given $\pairTmn\in\mathcal{P}_T$, 
the uniqueness of $A\subset \ComplTmn$ is guaranteed
by Corollary \ref{cor:example_T_mn_unique_A_p_odd} 
when $\mu=\nu=\pm 3$ and by Corollary \ref{cor:example_T_mn_unique_A_non_primitive}  
for the other cases. Hence, the slope type $(\frac{\mu+1}{2},\frac{\nu-1}{2})$
depends only on the isotopy class of $\pairTmn$,
so members in $\mathcal{P}_T$ are all inequivalent. 
\end{proof}

Similarly, for the handlebody-knot family $\pairImn$ in Section \ref{subsec:examples_irreatoro}, we have the following 
by Corollary \ref{cor:uniqueness}.
\begin{corollary}\label{cor:example_Imn_unique_A}
Suppose $\mu+\nu+3\neq 6$, and 
$\mu+1,\mu+2,\nu+1,\nu+2$ are not divisible by $p$
Then $A\subset\ComplImn$ is the unique type $3$-$3$ annulus, up to isotopy. 
\end{corollary}
\begin{proof}
The first criterion is equivalent to saying $\vert p\vert=1$,
whereas the second implies that $l_+,l_-$
do not represent the $\absp$-th multiples of some elements 
in $H_1(\ComplImnA)$, for 
\[[l_+]=(\mu+2,\nu+1),\quad [l_-]=(\mu+1,\nu+2)\in 
H_1(\ComplImnA)
\]
in terms of the meridional basis of $H_1(\ComplImn)$ 
given by $x_1,x_2$ in Fig.\ \ref{fig:spine_5_2} with  
$A$ oriented as in Fig.\ \ref{fig:mn_twisted_annulus}. 
\end{proof}

\begin{corollary}\label{cor:infinite_family_Imn}
The handlebody-knot family
\begin{equation}\label{eq:family_Imn}
\{\pairImn\mid \mu,\nu\in\mathbb{Z}\}
\end{equation}
contains infinitely many irreducible, atoroidal handlebody-knots
whose exteriors admit a unique type $3$-$3$ annulus.
\end{corollary}
\begin{proof}
Consider the subfamily of \eqref{eq:family_Imn}
\[\mathcal{P}_I:=\{\pairImn\mid \mu\geq \nu> -1 \text{ or } -2>\mu\geq \nu, \text{ and } \mu+\nu+3\neq 6\}.\]
The condition
$\mu\geq \nu> -1$ or $-2>\mu\geq \nu$
implies that the meridional basis of $H_1(\ComplImnA)$
given in Fig.\ \ref{fig:spine_5_2} is normalized and the orientation
of $A$ in Fig.\ \ref{fig:mn_twisted_annulus}
is a preferred one. Thus, $(\mu+2,\nu+1)$
is the slope type of $A$, and depends only on the isotopy class
of $\pairImn$ by the uniqueness of $A$
following from Corollary \ref{cor:example_Imn_unique_A}. 
Members in $\mathcal{P}_I$ are therefore mutually inequivalent. 
\end{proof}

\subsection{When homology criteria fail}\label{subsec:homotopy_criteria} 
For handlebody-knots 
in Corollaries \ref{cor:example_Tmn_irre_atoro}, \ref{cor:example_T_mn_unique_A_p_odd}, 
and \ref{cor:example_T_mn_unique_A_non_primitive},
homology criteria in Corollaries \ref{cor:irreducibility:trivial_HKA}
and \ref{cor:uniqueness}
provide a simple way to detect their irreducibility and the uniqueness
of $A$.
In some cases though the homology criteria are
not strong enough to ``see'' the irreducibility and uniqueness of $A$,
and homotopy criteria in Lemma \ref{lm:irreducibility:trivial_HKA}
and Theorem \ref{teo:uniqueness_HKA_irreducible}
are called for.

As an example, consider the handlebody-knots $\pairTmn$ 
with $(\mu,\nu)=(-2p+1,4p-1)$,
$\vert\frac{\mu+\nu}{2}\vert=\absp>1$.  
In terms of the meridional basis given 
in Fig.\ \ref{fig:tunnel_trivial_knot},
\[[l_+]=(-p+1,2p-1),\hspace*{.5em} [l_-]=(-p,2p)\text{ in } H_1(\ComplTmnA).\]
As $[l_-]$
is the $\absp$-th multiple of $(-1,2)$, Corollary \ref{cor:irreducibility:trivial_HKA} is not applicable here.

On the other hand, $l_+,l_-$
determine the conjugacy classes
of the cyclically reduced words 
$x_1^{-p+1}x_2^{2p-1},x_1^{-p}x_2^{2p}$,
respectively; neither is 
the $\absp$-th power of an element in $\pi_1(\ComplTmnA)$,
so applying Lemmas \ref{lm:irreducibility:trivial_HKA}
and \ref{lm:atoroidality:trivial_HKA}
and Theorem \ref{teo:uniqueness_HKA_trivial_not_primitive},
we obtain the following.
\begin{corollary}\label{cor:example_w_homotopy_criteria_1}
$\{(\sphere,\mathcal{T}_{-2p+1,4p-1}\mid \absp>1\}$
is an infinite family of irreducible, atoroidal
handlebody-knots whose exteriors admit a unique type $3$-$3$ annulus.
\end{corollary}
\begin{proof}
Members in the family
are mutually inequivalent by the uniqueness of $A$.
\end{proof}

It is worth noting that the homology criterion 
in Corollary \ref{cor:irreducibility:trivial_HKA} 
exclude the case $\vert p\vert<3$; 
many small crossing handlebody-knots, however, fall 
into this category. 
To see their irreducibility and atoroidality, and the uniqueness of $A$, one can employ Lemma \ref{lm:irreducibility:trivial_HKA} and Corollary \ref{cor:rigidity_slope_pm1}.

For instance, consider the handlebody-knot family 
\begin{equation}\label{eq:family_Tmu}
\{\pairTmu\mid \mu< -1\}.
\end{equation} 
Note that because $\pairTmu=(\sphere,\mathcal{T}_{2-\mu,\mu})$, we have 
\[\{\pairTmu\mid \mu< -1\}=\{\pairTmu\mid \mu>3\},\]
and when $\mu= \pm 1$ or $3$, $\pairTmu$ is trivial; also,
the canonical type $3$-$3$ annulus $A\subset\ComplTmu$ 
has a boundary slope of $1$, for every $\pairTmu$. 
In terms of the basis $x_1,x_2$ of $\pi_1(\ComplTmuA)$ in Fig.\ \ref{fig:tunnel_trivial_knot},
$l_+,l_-$
determine the conjugacy classes
of the cyclically reduced words 
\[x_1^{\frac{\mu+1}{2}}x_2^{\frac{1-\mu}{2}},\quad x_1^{\frac{\mu-1}{2}}x_2^{\frac{3-\mu}{2}},\]
respectively. Since none of the exponents is $0$,
and exponents of $x_1$ (resp.\ $x_2$) 
are not $\pm 1$ simultaneously,
$\{l_+,l_-\}$ does not represent a basis of $\pi_1(\ComplTmuA)$
by \cite{CohMetZim:81}, and therefore we have the following by
Corollary \ref{cor:rigidity_slope_pm1}.
\begin{corollary}
$\{\pairTmu\mid \mu<-1\}$
is an infinite family of irreducible, atoroidal
handlebody-knots whose exteriors admit a unique type $3$-$3$ annulus. 
\end{corollary}
\begin{proof}
The irreducibility, atoroidality, and 
the uniqueness of $A$ follows from Corollary \ref{cor:rigidity_slope_pm1}, so
it suffices to show that members in the family are mutually inequivalent.
Suppose $f$ is an equivalence between
$\pairTmu$
and $(\sphere,\mathcal{T}_{\mu',2-\mu'})$ with 
$\mu\neq \mu'$ and $\mu,\mu'  <-1$. 
Let $A\subset\ComplTmu, A'\subset\Compl{\mathcal{T}_{\mu',2-\mu'}}$ 
be the type $3$-$3$ annuli 
given by the construction in Section \ref{subsec:examples_irreatoro}. 
By the uniqueness of $A,A'$, we may assume 
$f(\rnbhd{A})=\rnbhd{A'}$. 

Let $l_\pm$ (resp.\ $l_\pm'$)
be essential loops of the annular components of 
$\rnbhd{A}\cap \partial\mathcal{T}_{\mu,2-\mu,A}$ 
(resp.\ $\rnbhd{A'}\cap\partial\mathcal{T}_{\mu',2-\mu',A'}$), respectively.
Then it may be assumed that $f$ sends $l_\pm$ either to $l_\pm'$
or to $l_\mp'$. Therefore, 
the induced homomorphism
\[f_\ast:\pi_1(\ComplTmuA)\rightarrow 
\pi_1(\Compl{\mathcal{T}_{\mu',2-\mu',A'}})\]
sends the conjugacy class
of $x_1^{\frac{\mu+1}{2}}x_2^{\frac{1-\mu}{2}}$ (resp.\ $x_1^{\frac{\mu-1}{2}}x_2^{\frac{3-\mu}{2}}$) 
to the conjugacy class of either 
$x_1^{\frac{\mu'+1}{2}}x_2^{\frac{1-\mu'}{2}}$ or $x_1^{\frac{\mu'-1}{2}}x_2^{\frac{3-\mu'}{2}}$.

Denote by $X_\mu^+$ (resp.\ $X_\mu^-$) 
the $3$-manifold 
obtained by attaching a $2$-handle 
to $\ComplTmuA$ 
along $l_+$ (resp.\ $l_-$). Observe that 
when $\mu<-3$,
its fundamental group 
is isomorphic to the torus knot group
$\pi_1(\ComplKpq)$ with $(p,q)=(\frac{-1-\mu}{2},\frac{1-\mu}{2})$
(resp.\ $(p,q)=(\frac{1-\mu}{2},\frac{3-\mu}{2})$), whereas $\pi_1(X_{-3}^+)\simeq \mathbb{Z}$
and $\pi_1(X_{-3}^+)\simeq \pi_1(\Compl{\mathcal{K}_{2,3}})$.
Thus we may assume $\mu,\mu'<-3$. 

Since torus knots are distinguished by their knot groups, up to mirror 
image.
If $f(l_\pm)=l_\pm'$, then $f_\ast$ implies
that 
$(\frac{-1-\mu}{2},\frac{1-\mu}{2})$-
and 
$(\frac{-1-\mu'}{2},\frac{1-\mu'}{2})$-torus knots
(resp.\ $(\frac{1-\mu}{2},\frac{3-\mu}{2})$-
and $(\frac{1-\mu'}{2},\frac{3-\mu'}{2})$-torus knots) 
are equivalent, up to mirror image; this happens only when $\mu=\mu'$,
contradicting the assumption. 
Similarly, if $f(l_\pm)=l_\mp'$, then we have 
$(\frac{-1-\mu}{2},\frac{1-\mu}{2})$-
and 
$(\frac{1-\mu'}{2},\frac{3-\mu'}{2})$-torus knots
(resp.\ $(\frac{1-\mu}{2},\frac{3-\mu}{2})$-
and $ (\frac{-1-\mu'}{2},\frac{1-\mu'}{2})$-torus knots)
are equivalent, up to mirror image, but this too is an impossibility 
under the assumption $\mu,\mu'<-1$.
\end{proof}

\begin{figure}[h] 
\begin{subfigure}{0.49\textwidth}
\centering
\def\svgwidth{.9\columnwidth}
\begingroup%
  \makeatletter%
  \providecommand\color[2][]{%
    \errmessage{(Inkscape) Color is used for the text in Inkscape, but the package 'color.sty' is not loaded}%
    \renewcommand\color[2][]{}%
  }%
  \providecommand\transparent[1]{%
    \errmessage{(Inkscape) Transparency is used (non-zero) for the text in Inkscape, but the package 'transparent.sty' is not loaded}%
    \renewcommand\transparent[1]{}%
  }%
  \providecommand\rotatebox[2]{#2}%
  \newcommand*\fsize{\dimexpr\f@size pt\relax}%
  \newcommand*\lineheight[1]{\fontsize{\fsize}{#1\fsize}\selectfont}%
  \ifx\svgwidth\undefined%
    \setlength{\unitlength}{1559.05511811bp}%
    \ifx\svgscale\undefined%
      \relax%
    \else%
      \setlength{\unitlength}{\unitlength * \real{\svgscale}}%
    \fi%
  \else%
    \setlength{\unitlength}{\svgwidth}%
  \fi%
  \global\let\svgwidth\undefined%
  \global\let\svgscale\undefined%
  \makeatother%
  \begin{picture}(1,0.67272727)%
    \lineheight{1}%
    \setlength\tabcolsep{0pt}%
    \put(0,0){\includegraphics[width=\unitlength,page=1]{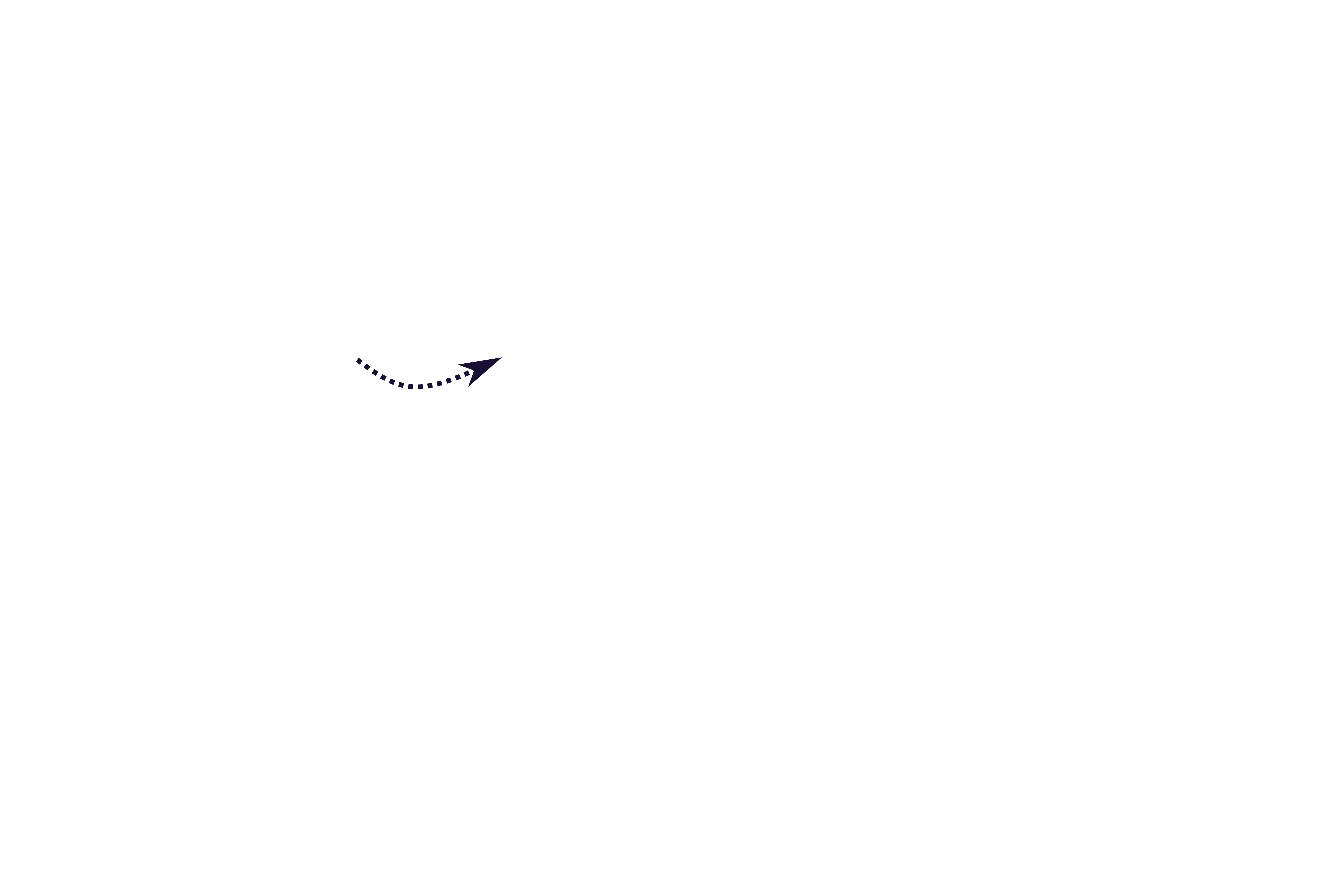}}%
    \put(0.02100209,0.02878374){\color[rgb]{0,0,0}\makebox(0,0)[lt]{\lineheight{1.25}\smash{\begin{tabular}[t]{l}{\footnotesize $(\sphere, \mirror 6_4)$}\end{tabular}}}}%
    \put(0.00009125,0.64264471){\color[rgb]{0,0,0}\makebox(0,0)[lt]{\lineheight{1.25}\smash{\begin{tabular}[t]{l}{\footnotesize $(\sphere,\mathcal{T}_{-3,5})$}\end{tabular}}}}%
    \put(-0.15690161,-0.57813665){\color[rgb]{0,0,0}\makebox(0,0)[lt]{\begin{minipage}{1.39515201\unitlength}\raggedright \end{minipage}}}%
    \put(0,0){\includegraphics[width=\unitlength,page=2]{equivalence_to_m6_4.pdf}}%
  \end{picture}%
\endgroup%

\caption{Equivalence: $(\sphere,\mathcal{T}_{-3,5})$ and $(\sphere,\mirror 6_4)$.}
\label{fig:equivalence_to_m6_4} 
\end{subfigure}
\begin{subfigure}{0.49\textwidth}
\centering
\def\svgwidth{.9\columnwidth}
\begingroup%
  \makeatletter%
  \providecommand\color[2][]{%
    \errmessage{(Inkscape) Color is used for the text in Inkscape, but the package 'color.sty' is not loaded}%
    \renewcommand\color[2][]{}%
  }%
  \providecommand\transparent[1]{%
    \errmessage{(Inkscape) Transparency is used (non-zero) for the text in Inkscape, but the package 'transparent.sty' is not loaded}%
    \renewcommand\transparent[1]{}%
  }%
  \providecommand\rotatebox[2]{#2}%
  \newcommand*\fsize{\dimexpr\f@size pt\relax}%
  \newcommand*\lineheight[1]{\fontsize{\fsize}{#1\fsize}\selectfont}%
  \ifx\svgwidth\undefined%
    \setlength{\unitlength}{1303.93700787bp}%
    \ifx\svgscale\undefined%
      \relax%
    \else%
      \setlength{\unitlength}{\unitlength * \real{\svgscale}}%
    \fi%
  \else%
    \setlength{\unitlength}{\svgwidth}%
  \fi%
  \global\let\svgwidth\undefined%
  \global\let\svgscale\undefined%
  \makeatother%
  \begin{picture}(1,0.65217391)%
    \lineheight{1}%
    \setlength\tabcolsep{0pt}%
    \put(0,0){\includegraphics[width=\unitlength,page=1]{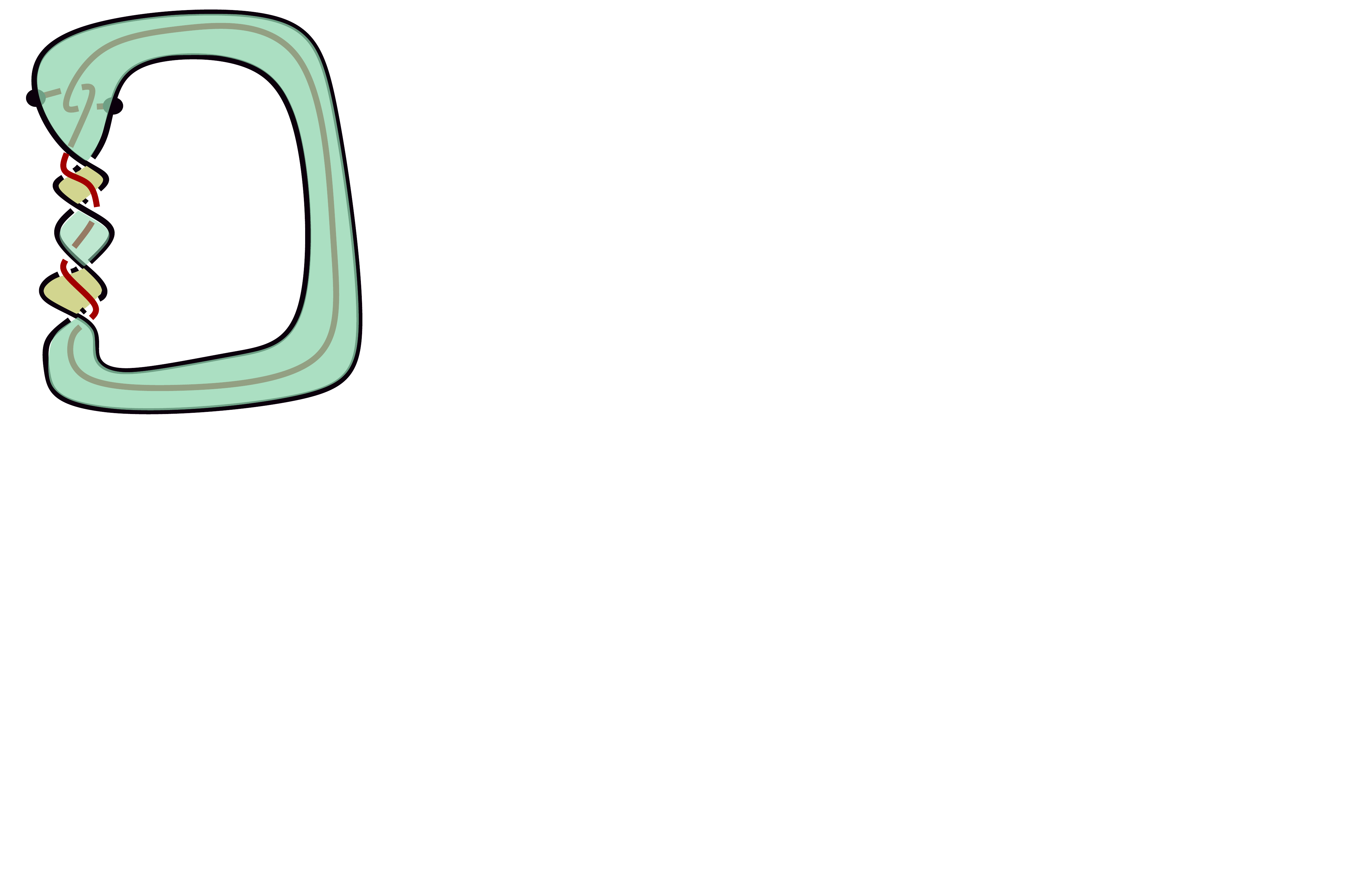}}%
    \put(0.00585832,0.29441333){\color[rgb]{0,0,0}\makebox(0,0)[lt]{\lineheight{1.25}\smash{\begin{tabular}[t]{l}{\footnotesize $(\sphere, \HK_{A_1})$}\end{tabular}}}}%
    \put(0,0){\includegraphics[width=\unitlength,page=2]{equivalence_to_mT3_3.pdf}}%
    \put(0.00768118,0.00988792){\color[rgb]{0,0,0}\makebox(0,0)[lt]{\lineheight{1.25}\smash{\begin{tabular}[t]{l}{\footnotesize $(\sphere,\mathcal{T}_{-3,-3})$}\end{tabular}}}}%
    \put(-0.15033753,0.77093344){\color[rgb]{0,0,0}\makebox(0,0)[lt]{\begin{minipage}{1.60046836\unitlength}\raggedright \end{minipage}}}%
    \put(0,0){\includegraphics[width=\unitlength,page=3]{equivalence_to_mT3_3.pdf}}%
  \end{picture}%
\endgroup%

\caption{Equivalence: $(\sphere,\HK_{A_1})$ and $(\sphere,\mathcal{T}_{-3,-3})$.}
\label{fig:equivalence_to_mT3_3} 
\end{subfigure}
\caption{Equivalences between handlebody-knots.}
\end{figure}

\begin{remark}
$(\sphere,\mathcal{T}_{3,3})$ and $(\sphere,\mathcal{T}_{-3,5})$
are equivalent 
to the mirror images of handlebody-knots $(\sphere, 5_2)$ and $(\sphere, 6_4)$
in \cite[Table $1$]{IshKisMorSuz:12} as demonstrated in Figs.\ \ref{fig:equivalence_to_m5_2} and \ref{fig:equivalence_to_m6_4},
respectively. Hence we obtain an alternative proof
of their irreducibility \cite[Section $4$]{IshKisMorSuz:12}
and $A$ being the unique type $3$-$3$ annulus in their exteriors (compare with \cite[Lemmas
$3.4$ and $4.5$]{LeeLee:12}. 
Also, the type $3$-$3$ annuli in $\Compl{\mathcal{T}_{3,3}},\Compl{\mathcal{T}_{-3,5}}$ having different slope types
is another way to see the inequivalence between $(\sphere,5_2)$ and $(\sphere,6_4)$ \cite[Table $2$]{IshKisMorSuz:12}.
\end{remark}

\subsection{Non-uniqueness}\label{subsec:example_non_uniqueness}
\begin{figure}[h]
\begin{subfigure}{0.32\textwidth}
\centering
\def\svgwidth{.8\columnwidth}
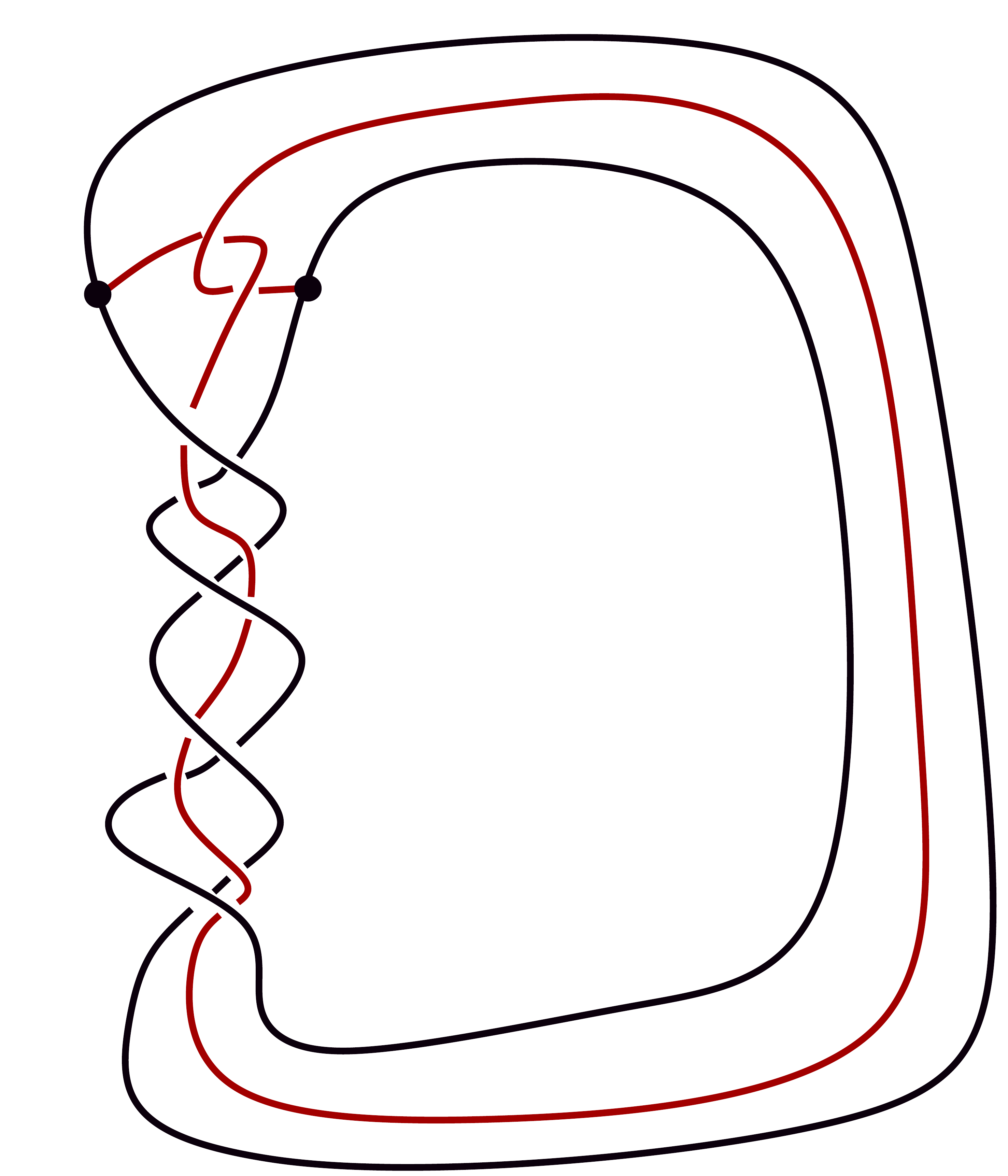
\caption{Spine of $\pair$.}
\label{fig:non_example_spine}
\end{subfigure}
\begin{subfigure}{0.33\textwidth}
\centering
\def\svgwidth{.8 \columnwidth}
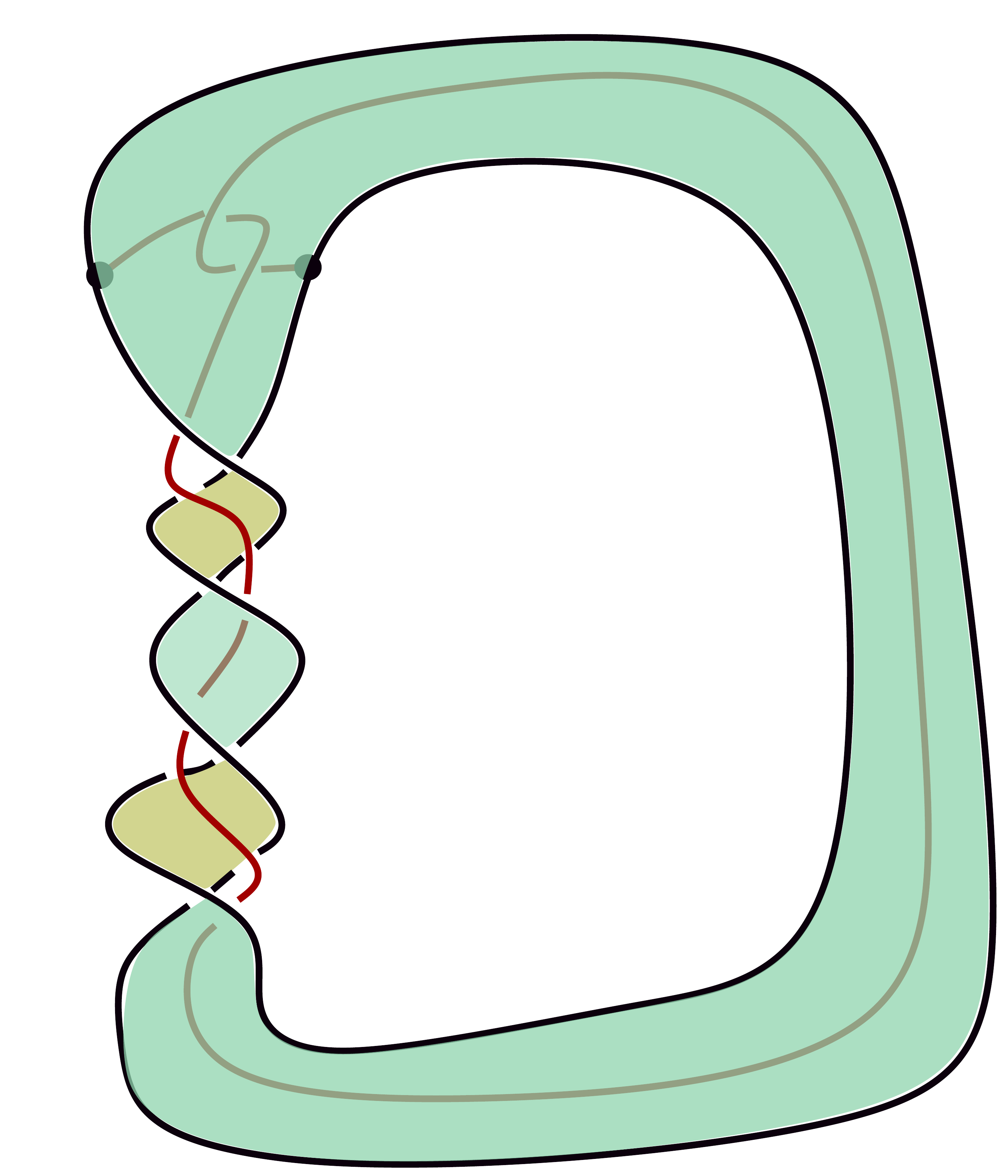
\caption{Annulus $A_1$.}
\label{fig:non_example_annulus1}
\end{subfigure}
\begin{subfigure}{0.33\textwidth}
\centering
\def\svgwidth{.8\columnwidth}
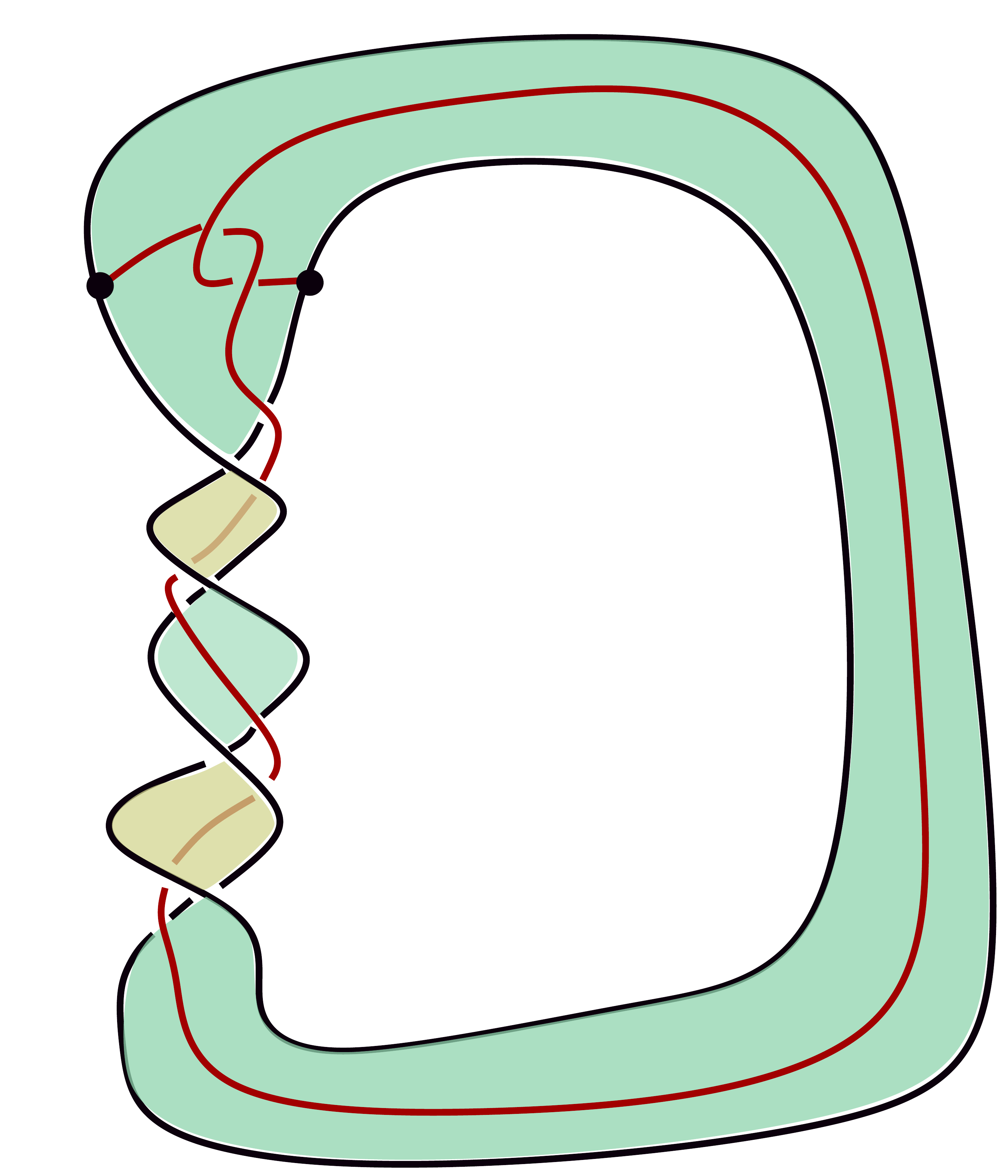
\caption{Annulus $A_2$.}
\label{fig:non_example_annulus2}
\end{subfigure}
\caption{Irreducible, atoroidal $\pair$ with non-isotopic type $3$-$3$ annuli.}
\end{figure}

Let $\pair$ be the handlebody-knot given by
a regular neighborhood of the 
handcuff spatial graph in Fig.\ \ref{fig:non_example_spine}.
$\ComplHK$ contains two non-isotopic
type $3$-$3$ annuli $A_1,A_2$ with a boundary slope of $2$
as shown in Figs.\ \ref{fig:non_example_annulus1} and \ref{fig:non_example_annulus2}, respectively.

The irreducibility and atoroidality of $\pair$ follow from the fact 
that $(\sphere, \HK_{A_1})$ is equivalent to
$(\sphere, \mathcal{T}_{-3,-3})$ 
as shown in Fig.\ \ref{fig:equivalence_to_mT3_3} 
and Corollary \ref{cor:irreducibility:irreducible_HKA}
and Lemma \ref{lm:atoroidality_general}.
%

Let $l_\pm$ be essential loops of the two annular components of $\rnbhd{A_1}\cap \partial\HK_{A_1}$, respectively. Then one of 
them represents the square of some element in $\pi_1(\Compl{\HK_{A_1}})$---failing the criterion in Theorem \ref{teo:uniqueness_HKA_irreducible}.
The same is true for $A_2$, which is unknotting,
so this example also shows that criteria in 
Lemma \ref{lm:irreducibility:trivial_HKA}
are not a necessary condition for a handlebody-knot to be irreducible.


\section{Symmetry}\label{sec:symmetry}
%
\subsection{Symmetry group}
Throughout the subsection $\pair$ 
is an irreducible, atoroidal handlebody-knot
whose exterior contains 
a unique type $3$-$3$ annulus $A$
with a non-trivial boundary slope of $p$, up to isotopy.
Let $(p_1,p_2)$ be the slope type of $A$.
Orient $A$ and its boundary components $l_1,l_2$
so that $l_1\cup -l_2=\partial A$.
Let $\mathbf{D}_A=\{\disk_1,\disk_2,\disk_A\}$ 
be the meridian system induced by $A$ 
and $\Gamma_A$ the associated spine of $\HK$ (see
Section \ref{subsec:typethreethree}). 
We fix an orientation of $\disk_i$, 
$i=1,2$.

\begin{theorem}\label{teo:injection}
The composition
\begin{equation}\label{eq:composition_symmetry_groups}
\pSym{\HK}\simeq \pSym{\HK,A}\xrightarrow{\pi} \MCG{\mathring{A}}\simeq \mathbb{Z}_2\times\mathbb{Z}_2
\end{equation}
is an injection.
\end{theorem}
\begin{proof}
Given $f\in\pAut{\sphere,\HK}$,
suppose $f\vert_{\mathring{A}}$
is isotopic to the identity. 
Then $f(l_i)=l_i$ and $f\vert_{l_i}$
is orientation-preserving, $i=1,2$.
By Corollary \ref{cor:isom_mapping_class_groups},
we have the isomorphism 
\[\pSym{\HK,\cup\mathbf{D}_A}\rightarrow \pSym\HK,\]
and hence $f$ can be isotoped in $\pAut{\sphere,\HK}$
to a homeomorphism $f'$ with $f'(\disk_A)=\disk_A, f'(\disk_i)=\disk_i$, $i=1,2$,
and $f'\vert_{\disk_i}$ orientation-preserving.
%
The homomorphism  
\begin{equation}\label{eq:pSym_w_disk_to_pSym_graph}
\pSym{\HK,\cup\mathbf{D}_A} \rightarrow \pSym{\Gamma_A}
\end{equation} 
in Lemma \ref{lm:symmetry_groups_hk_gamma} given by the Alexander trick
allows us to further isotope $f'$ in
\[ \pAut{\sphere,\HK,\cup\mathbf{D}_A}\] to a homeomorphism $f''\in \pAut{\sphere,\HK, \cup\mathbf{D}_A,\Gamma_A}$
that restricts to the identity on $\Gamma_A$. 
The injectivity of \eqref{eq:composition_symmetry_groups} then follows from Lemma \ref{lm:spine_atoro_irre_hk}
and \eqref{eq:pSym_w_disk_to_pSym_graph} being an 
isomorphism by Lemma \ref{lm:symmetry_groups_hk_gamma}.
\end{proof}

The next two corollaries 
follows readily from Theorem \ref{teo:injection}.

\begin{corollary}\label{cor:symmetry_group_general}
$\Sym\HK=\pSym\HK\leq \mathbb{Z}_2\times \mathbb{Z}_2$,
\end{corollary}
\begin{proof}
It suffices to show that $\pair$ is chiral.
Since $A$ is unique, any $f\in\Aut{\sphere,\HK}$
can be isotoped so that $f(l_1\cup l_2)=l_1\cup l_2$ or $-l_1\cup -l_2$.
If $f$ is orientation-reversing, then
$\lk{f(l_1)}{f(l_2)}=-p$, but at the same time, 
$\lk{f(l_1)}{f(l_2)}=\lk{l_1}{l_2}=p$, contradicting $p\neq 0$.
\end{proof}

Combining Corollary \ref{cor:symmetry_group_general} with Lemma \ref{cor:reversing_orientation_A},
we have the following.
\begin{corollary}\label{cor:symmetry_group_asymmetric_slope_type}
If $(p_1,p_2)\neq(\frac{p+1}{2},\frac{p-1}{2})$,
then $\Sym\HK=\pSym\HK\leq \mathbb{Z}_2$.
\end{corollary}
\begin{proof}
By Lemma \ref{cor:reversing_orientation_A},
there exists no homeomorphism 
$f\in\Aut{\sphere,\HK,A}$
whose restriction $f\vert_A$ on $A$ is orientation-reversing unless
$(p_1,p_2)= (\frac{p+1}{2},\frac{p-1}{2})$, so if
$(p_1,p_2)\neq(\frac{p+1}{2},\frac{p-1}{2})$, the homomorphism \eqref{eq:composition_symmetry_groups}
is not surjective.
\end{proof}

\subsection{Examples: optimal upper bounds}\label{subsec:example_full}
Here we show that
the group $\mathbb{Z}_2\times\mathbb{Z}_2$ 
(resp.\ $\mathbb{Z}_2$)
in Corollary \ref{cor:symmetry_group_general} (resp.\ \ref{cor:symmetry_group_asymmetric_slope_type})
is optimal, in the sense that
there are handlebody-knots
satisfying the given condition and having a symmetry group 
isomorphic to the group.

\subsubsection{$\mathbb{Z}_2\times\mathbb{Z}_2$}

Consider the subfamily 
\[\mathcal{V}:=\{\pairTp\mid p \text{ odd }, \vert p\vert >1\}\]
of the handlebody-knot family $\{\pairTmn\mid \mu,\nu \text{ odd }\}$,
and let $A$ be the type $3$-$3$ annulus given by the construction in Section \ref{subsec:examples_irreatoro}.
By Corollaries \ref{cor:example_Tmn_irre_atoro} and  
\ref{cor:example_T_mn_unique_A_p_odd}. 
every $\pair\in\mathcal{V}$ is irreducible and atoroidal,
and $A\subset \Compl\HK$ is unique type $3$-$3$ annulus, up to isotopy.
 
\begin{figure}[b]
\centering
\def\svgwidth{.7\columnwidth}
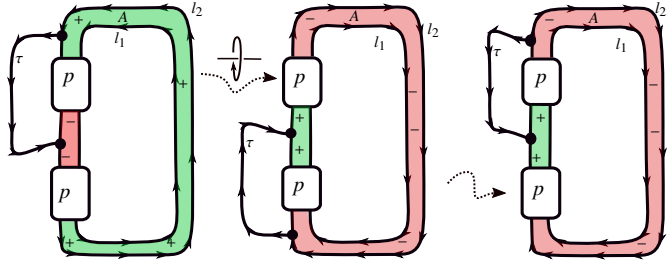
\caption{$g_1$: flip $A$ $+$ reverse $l_1,l_2$.}
\label{fig:V_generator1}
\end{figure}

\begin{figure}[b]
\centering
\def\svgwidth{.9\columnwidth}
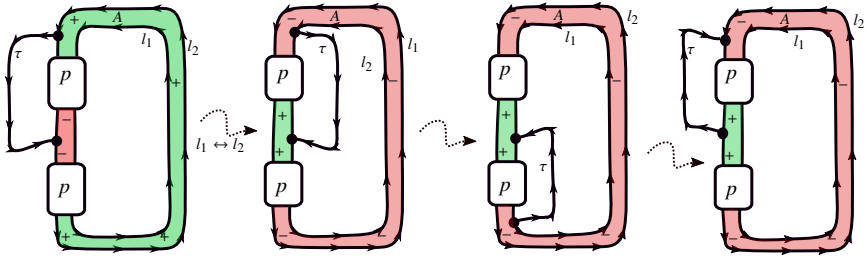
\caption{$g_2$: flip $A$ $+$ swap $l_1,l_2$.}
\label{fig:V_generator2}
\end{figure}

\begin{corollary}\label{cor:symmetry_group_Vp}
$\Sym\HK=\pSym\HK=\mathbb{Z}_2\times \mathbb{Z}_2$, 
for any $\pair\in\mathcal{V}$.
\end{corollary}
\begin{proof}
Note first  
$\Sym\HK=\pSym\HK\leq \mathbb{Z}_2\times \mathbb{Z}_2$
by Corollary \ref{cor:symmetry_group_general}.
For the other direction, 
consider the homeomorphisms $g_1,g_2$ given by 
the isotopies in Figs.\ \ref{fig:V_generator1}, \ref{fig:V_generator2},
respectively.
The isotopy in Fig.\ \ref{fig:V_generator1} 
is given by first rotating the diagram against 
a horizontal
line and then moving the lower end of the arc $\tau$
up counterclockwise along the untwisted part of $l_1$.
In particularly, $g_1\vert_A$ reverses the orientation of $A$,
but does not swap $l_1,l_2$. 
On the other hand, 
the isotopy in Fig.\ \ref{fig:V_generator2} 
is given by first
swapping $l_1,l_2$, and hence flipping $A$, then moving
the upper end of $\tau$ down clockwise, and 
then shifting two ends of
$\tau$ up simultaneously along the twisted part of $\partial A$.
Especially, $g_2\vert_A$ reverses the orientation of $A$,
and swaps $l_1,l_2$.

Since the elements represented by $g_1,g_2$ 
have different non-trivial images under the homomorphism
\begin{equation}\label{eq:restriction_symmetry_groups}
\Sym\HK\simeq \Sym{\HK,A}\rightarrow \MCG{A},
\end{equation} 
we see $\Sym\HK\geq \mathbb{Z}_2\times\mathbb{Z}_2$.
\end{proof}

\subsubsection{$\mathbb{Z}_2$}
\begin{figure}[h]
\centering
\def\svgwidth{.9 \columnwidth}
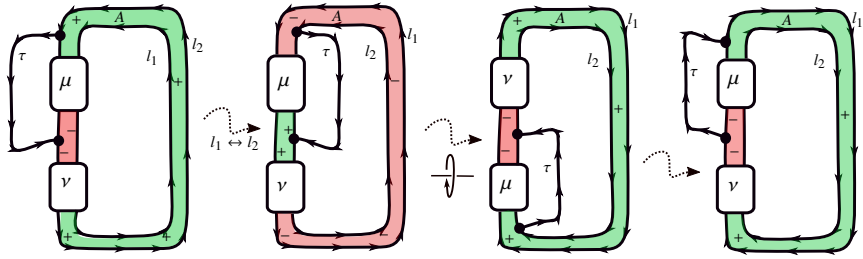
\caption{$g$: swap $l_1,l_2$ + reverse $l_1,l_2$.}
\label{fig:W_generator}
\end{figure}

Denote by $\mathcal{W}$ the subfamily
\[\{\pairTmn\mid \mu,\nu>1 \text{ or }\mu,\nu<-1, \mu\neq \nu, \text{ and }\mu,\nu \text{ odd }\}\] 
of the handlebody-knot family 
$\{\pairTmn\mid \mu,\nu \text{ odd }\}$. 
By Corollaries \ref{cor:example_Tmn_irre_atoro} and \ref{cor:example_T_mn_unique_A_non_primitive},
every $\pair\in\mathcal{W}$ is irreducible and atoroidal,
and the canonical annulus 
$A\subset \Compl\HK$ is the unique type $3$-$3$ annulus, up to isotopy.
\begin{corollary}
$\Sym\HK=\pSym\HK=\mathbb{Z}_2$, for any $\pair\in \mathcal{W}$.
\end{corollary}
\begin{proof}
The condition $\mu,\nu>1$ or $\mu,\nu<-1$ implies 
the slope type of $A$
is either $(\frac{\mu+1}{2},\frac{\nu-1}{2})$
or $(\frac{\nu+1}{2},\frac{\mu-1}{2})$, so 
$\mu = \nu$ if and only if the slope type is 
$(\frac{p+1}{2},\frac{p-1}{2})$; by Corollary \ref{cor:symmetry_group_asymmetric_slope_type}, 
\[\Sym\HK=\pSym\HK\leq \mathbb{Z}_2.\]

To see $\pSym\HK$ is non-trivial, we observe that 
the homeomorphism $g$ given by the isotopy 
in Fig.\ \ref{fig:W_generator}
represents a non-trivial element $g$ in $\pSym\HK$.
The isotopy 
starts with swapping $l_1,l_2$, then rotates the diagram against a horizontal line by $\pi$, and then moves the arc $\tau$ up along the twisted part of $\partial A$ by $\nu$ half crossings. 
In particular, the element represented by $g$ has a non-trivial image under the homomorphism\eqref{eq:restriction_symmetry_groups}.
\end{proof}

\subsection{Examples: not an isomorphism}\label{subsec:examples_not_isomorphism}
Here we show
that the inequalities $\leq$'s 
in Corollaries \ref{cor:symmetry_group_general}
and \ref{cor:symmetry_group_asymmetric_slope_type}
are in general not an isomorphism. 
Denote by $\mathcal{V}'$ the family 
\[\{\pairTmu\mid \mu  <-1\}\] 
of handlebody-knots in Section \ref{subsec:homotopy_criteria}. 
For any $\pair\in\mathcal{V}'$, 
the slope type of the unique type $3$-$3$
annulus $A\subset \Compl\HK$ is $(1,0)$, and
so Corollary \ref{cor:symmetry_group_asymmetric_slope_type} 
does not apply.
\begin{corollary}\label{cor:symmetry_group_Wm}
$\Sym\HK\simeq \pSym\HK\simeq \mathbb{Z}_2$,
for any $\pair\in\mathcal{V}'$.
\end{corollary}
\begin{proof}
Note first that, given $\pairTmu\in\mathcal{V}'$, 
the isotopy in Fig.\ \ref{fig:Tmu_generator} 
represents 
a non-trivial element in $\Sym\Tmu$.
\begin{figure}[h]
\centering
\def\svgwidth{.9 \columnwidth}
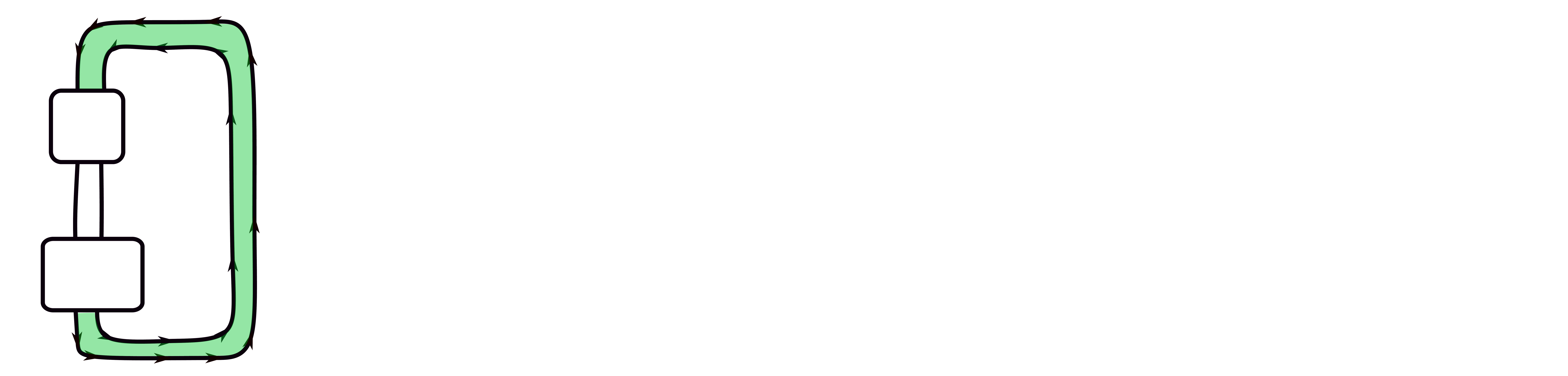
\caption{Swap and reverse $l_1,l_2$.}
\label{fig:Tmu_generator}
\end{figure}
Secondly, recall that $l_+,l_-$ determine 
the conjugacy classes represented by the cyclically reduced words
$x_1^{\frac{\mu+1}{2}}x_2^{\frac{1-\mu}{2}},
x_1^{\frac{\mu-1}{2}}x_2^{\frac{3-\mu}{2}}$, 
respectively, where $x_1,x_2\in \pi_1(\ComplTmuA)$ are generators given in Fig.\ \ref{fig:tunnel_trivial_knot}.  

Suppose the homomorphism
\eqref{eq:composition_symmetry_groups}
is surjective. Then there exists
a homeomorphism 
$f\in\pAut{\sphere,\Tmu, A}$ 
which swaps $l_+,l_-$.
%
Let $X_\mu^\pm$ 
be the spaces obtained by 
attaching a $2$-cell along $l_\pm$, respectively.
Then the induced homomorphism $f_\ast$ on $\pi_1(\ComplTmuA)$ 
gives an isomorphism between 
$\pi_1(X_\mu^+)$ and $\pi_1(X_\mu^-)$;
this contradicts the facts that 
$\pi_1(X_{-3}^+)$ is isomorphic to $\mathbb{Z}$
and $\pi_1(X_{-3}^-)$ isomorphic to the torus knot group $\pi_1(\Compl{\mathcal{K}_{2,3}})$ and that when $\mu<-3$, 
$(\frac{-1-\mu}{2},\frac{1-\mu}{2})$- 
and $(\frac{1-\mu}{2},\frac{3-\mu}{2})$-torus knots
are never equivalent, up to mirror image.
\end{proof}

\begin{remark} 
Since $(\sphere, \mathcal{T}_{3,3})\in\mathcal{V}$
and $(\sphere, \mathcal{T}_{-3,5})\in\mathcal{V}'$ are equivalent to
the mirror images of 
$(\sphere, 5_2)$ 
and $(\sphere, 6_4)$ in 
\cite{IshKisMorSuz:12}, respectively (Figs.\ \ref{fig:equivalence_to_m5_2} and \ref{fig:equivalence_to_m6_4},),   
we obtain 
\[\Sym{5_2}\simeq \pSym{5_2}\simeq \mathbb{Z}_2\times \mathbb{Z}_2, \text{ and }\Sym{6_4}\simeq \pSym{6_4}\simeq \mathbb{Z}_2.\] 
\end{remark}

For our last example, we observe that the proof of Corollary \ref{cor:symmetry_group_asymmetric_slope_type},
that is, Corollary \ref{cor:reversing_orientation_A} and Theorem \ref{teo:injection}, implies the following. 
\begin{corollary}\label{cor:non_invertiable_no_symmetries}
Suppose the slope type of $A$ is not $(\frac{p+1}{2},\frac{p-1}{2})$
and $l_i$, $i=1,2$, are non-invertible knots in $\sphere$.
Then $\Sym\HK\simeq \pSym\HK=1$. 
\end{corollary}   
\begin{proof}
Since the slope type of $A$ is not $(\frac{p+1}{2},\frac{p-1}{2})$,
any non-trivial element in $\pSym\HK$  
is represented by a homeomorphism $f\in\pAut{\sphere,\HK,A}$
whose restriction $f\vert_A$ on $A$ reverses the orientation of $l_1,l_2$.  
\end{proof}

\begin{figure}[b]
\begin{subfigure}{0.48\textwidth}
\centering
\def\svgwidth{.8\columnwidth}
\begingroup%
  \makeatletter%
  \providecommand\color[2][]{%
    \errmessage{(Inkscape) Color is used for the text in Inkscape, but the package 'color.sty' is not loaded}%
    \renewcommand\color[2][]{}%
  }%
  \providecommand\transparent[1]{%
    \errmessage{(Inkscape) Transparency is used (non-zero) for the text in Inkscape, but the package 'transparent.sty' is not loaded}%
    \renewcommand\transparent[1]{}%
  }%
  \providecommand\rotatebox[2]{#2}%
  \newcommand*\fsize{\dimexpr\f@size pt\relax}%
  \newcommand*\lineheight[1]{\fontsize{\fsize}{#1\fsize}\selectfont}%
  \ifx\svgwidth\undefined%
    \setlength{\unitlength}{1133.85826772bp}%
    \ifx\svgscale\undefined%
      \relax%
    \else%
      \setlength{\unitlength}{\unitlength * \real{\svgscale}}%
    \fi%
  \else%
    \setlength{\unitlength}{\svgwidth}%
  \fi%
  \global\let\svgwidth\undefined%
  \global\let\svgscale\undefined%
  \makeatother%
  \begin{picture}(1,0.75)%
    \lineheight{1}%
    \setlength\tabcolsep{0pt}%
    \put(0,0){\includegraphics[width=\unitlength,page=1]{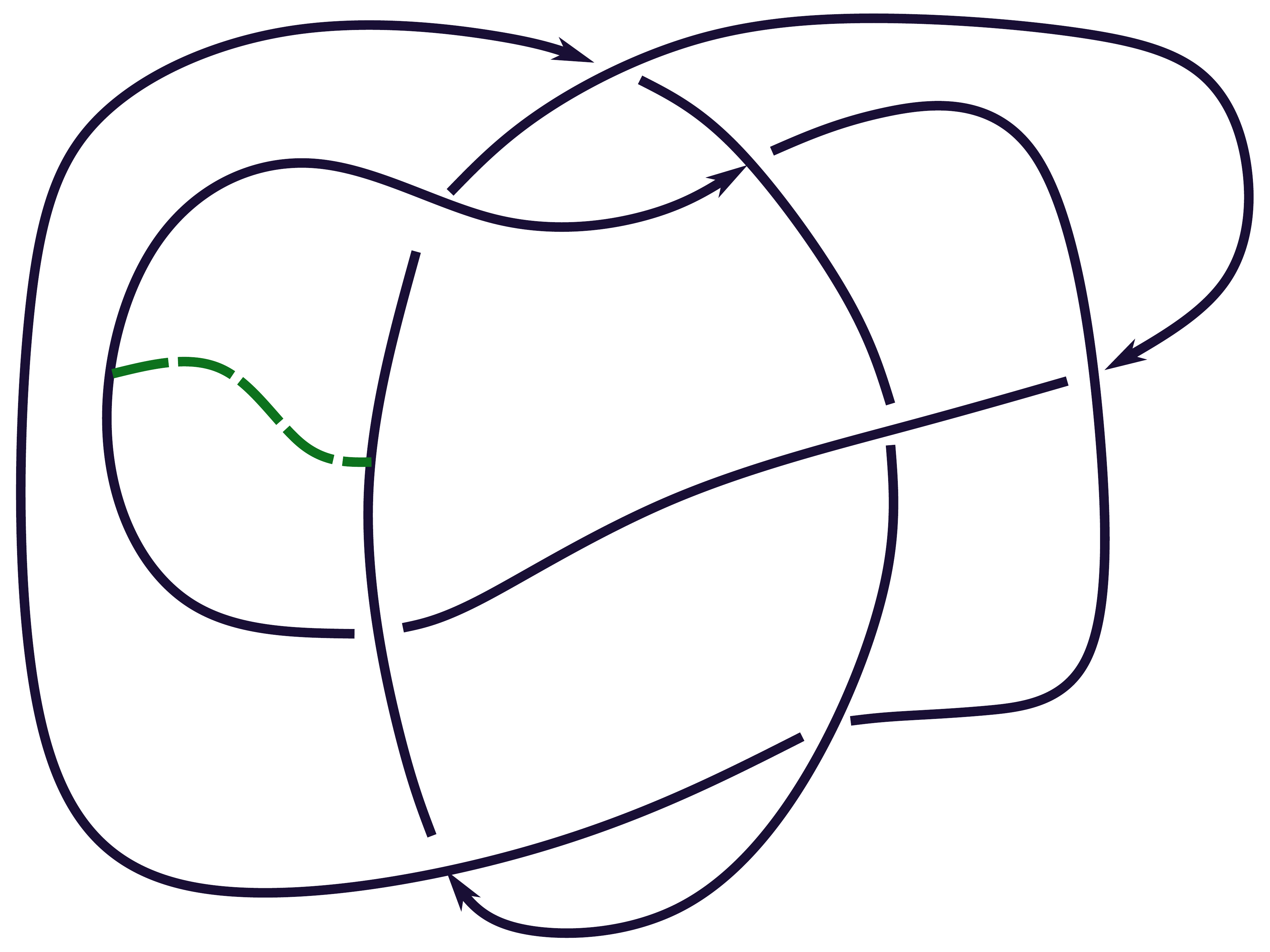}}%
    \put(0.14104192,0.407878){\color[rgb]{0,0,0}\makebox(0,0)[lt]{\lineheight{1.25}\smash{\begin{tabular}[t]{l}{\footnotesize $\tau$}\end{tabular}}}}%
    \put(0,0){\includegraphics[width=\unitlength,page=2]{8_16_w_tunnel.pdf}}%
    \put(0.17724373,0.563792){\color[rgb]{0,0,0}\makebox(0,0)[lt]{\lineheight{1.25}\smash{\begin{tabular}[t]{l}{\tiny $x_1$}\end{tabular}}}}%
    \put(0.03417205,0.30597801){\color[rgb]{0,0,0}\makebox(0,0)[lt]{\lineheight{1.25}\smash{\begin{tabular}[t]{l}{\tiny $x_2$}\end{tabular}}}}%
  \end{picture}%
\endgroup%

\caption{$8_{16}$ and $\tau$.}
\label{fig:8_16_w_tunnel}
\end{subfigure}
\begin{subfigure}{0.48\textwidth}
\centering
\def\svgwidth{.8\columnwidth}
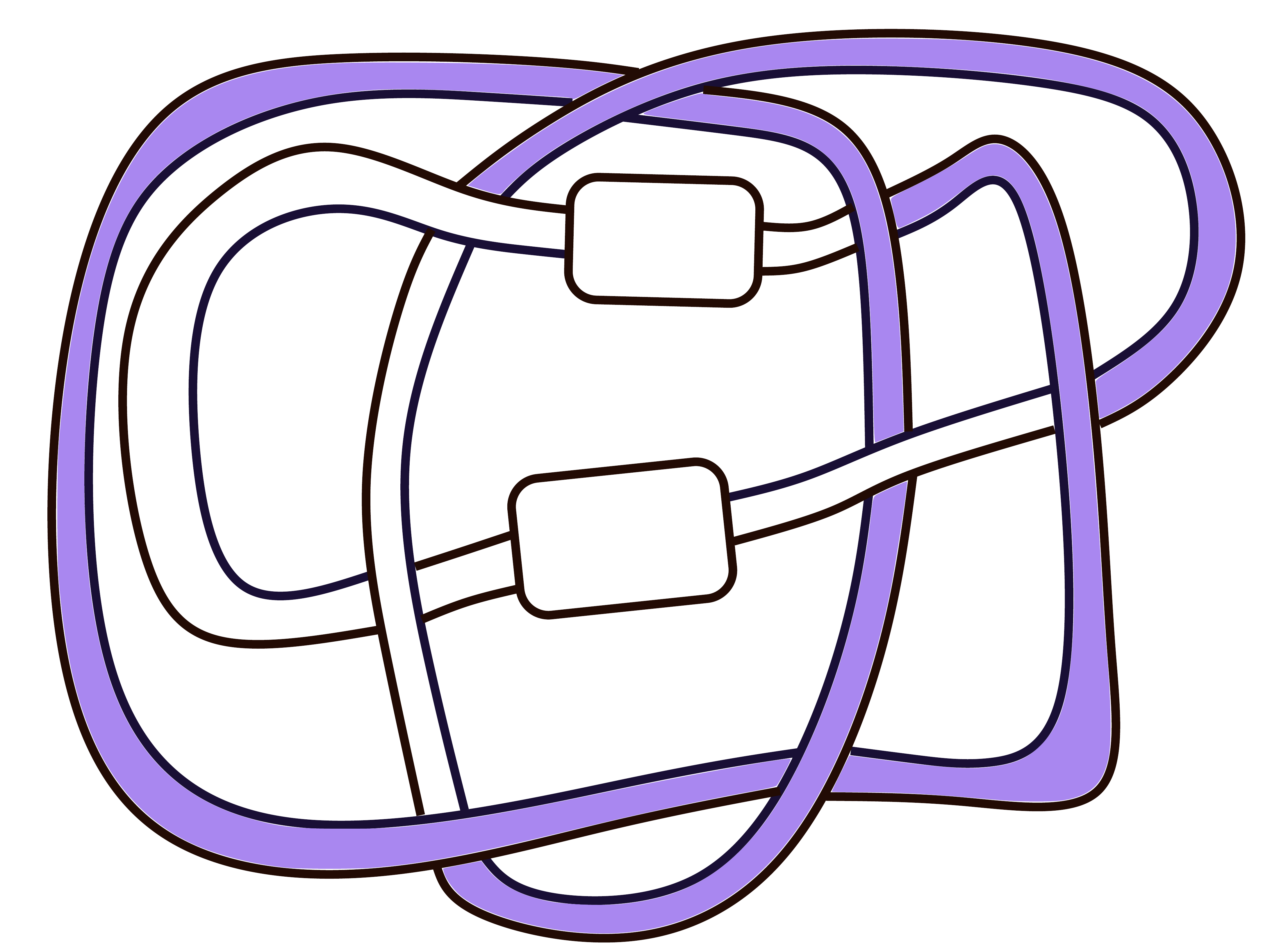
\caption{$p$-annulus $\mathcal{A}$ associated to $(8_{16},\tau)$.}
\label{fig:p_annulus_8_16_tau}
\end{subfigure} 
\medskip

%
\begin{subfigure}{0.8\textwidth}
\centering
\def\svgwidth{.6\columnwidth}
\begingroup%
  \makeatletter%
  \providecommand\color[2][]{%
    \errmessage{(Inkscape) Color is used for the text in Inkscape, but the package 'color.sty' is not loaded}%
    \renewcommand\color[2][]{}%
  }%
  \providecommand\transparent[1]{%
    \errmessage{(Inkscape) Transparency is used (non-zero) for the text in Inkscape, but the package 'transparent.sty' is not loaded}%
    \renewcommand\transparent[1]{}%
  }%
  \providecommand\rotatebox[2]{#2}%
  \newcommand*\fsize{\dimexpr\f@size pt\relax}%
  \newcommand*\lineheight[1]{\fontsize{\fsize}{#1\fsize}\selectfont}%
  \ifx\svgwidth\undefined%
    \setlength{\unitlength}{1559.05511811bp}%
    \ifx\svgscale\undefined%
      \relax%
    \else%
      \setlength{\unitlength}{\unitlength * \real{\svgscale}}%
    \fi%
  \else%
    \setlength{\unitlength}{\svgwidth}%
  \fi%
  \global\let\svgwidth\undefined%
  \global\let\svgscale\undefined%
  \makeatother%
  \begin{picture}(1,0.63636364)%
    \lineheight{1}%
    \setlength\tabcolsep{0pt}%
    \put(0,0){\includegraphics[width=\unitlength,page=1]{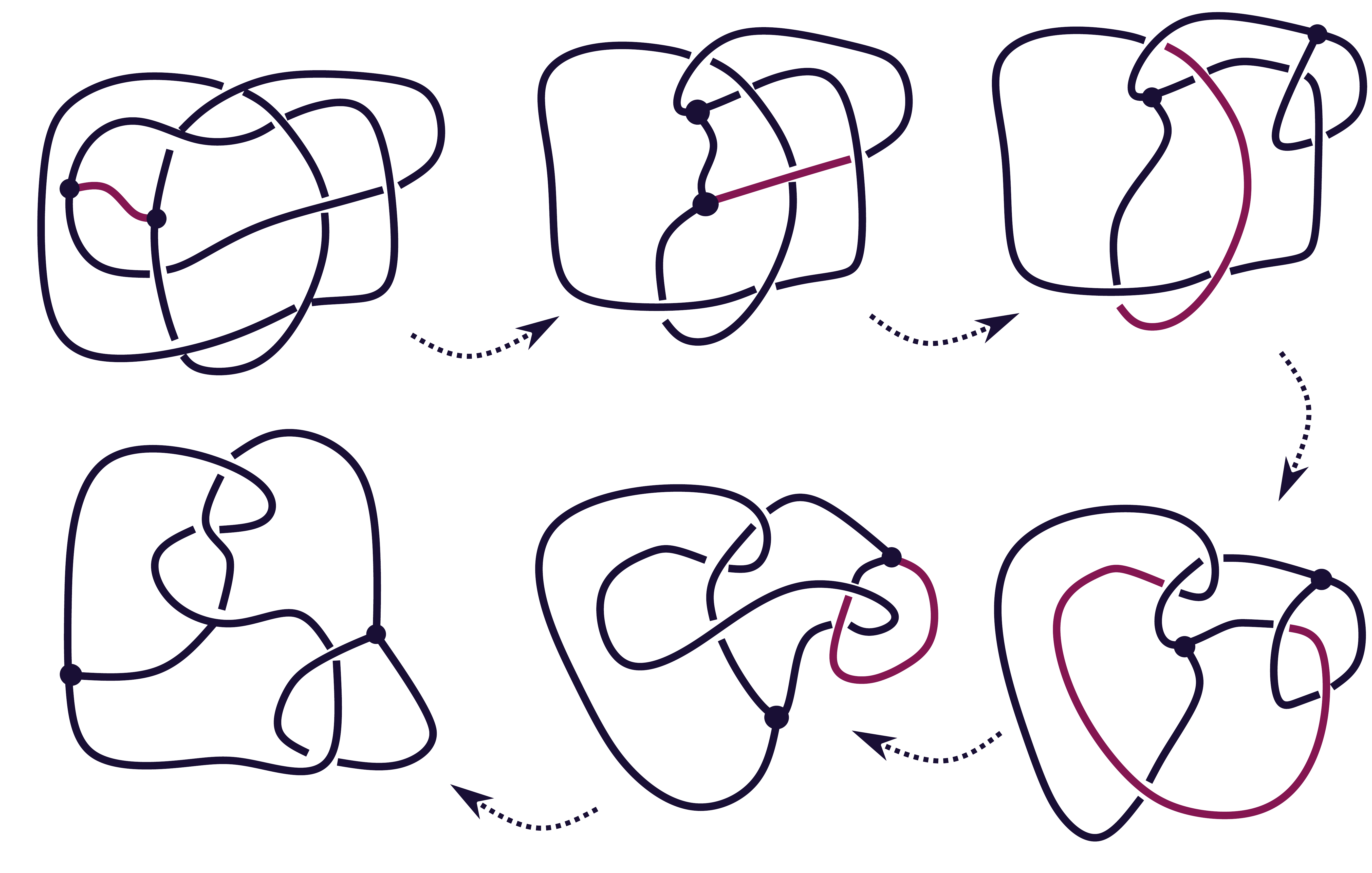}}%
    \put(0.00224398,0.01865995){\color[rgb]{0,0,0}\makebox(0,0)[lt]{\lineheight{1.25}\smash{\begin{tabular}[t]{l}{\footnotesize $(\sphere,\mirror 5_1)$}\end{tabular}}}}%
    \put(0.00524367,0.60700627){\color[rgb]{0,0,0}\makebox(0,0)[lt]{\lineheight{1.25}\smash{\begin{tabular}[t]{l}{\footnotesize $(\sphere, \rnbhd{8_{16}\cup\tau})$}\end{tabular}}}}%
    \put(-0.15690161,-0.57813665){\color[rgb]{0,0,0}\makebox(0,0)[lt]{\begin{minipage}{1.39515201\unitlength}\raggedright \end{minipage}}}%
  \end{picture}%
\endgroup%

\caption{Equivalence: $(\sphere,\rnbhd{8_{16}\cup\tau})$ and $(\sphere,\mirror 5_1)$.}
\label{fig:equivalence_to_m5_1}
\end{subfigure} 
\caption{Construction of $\pairUmn$.}
\end{figure}

For example, consider the knot $8_{16}$ 
in the Rolfsen knot table with 
a tunnel $\tau$ depicted in Fig.\ \ref{fig:8_16_w_tunnel}.
Let $\mathcal{A}$ be the $p$-annulus associated to $(8_{16},\tau)$
given in Fig.\ \ref{fig:p_annulus_8_16_tau}
where $\mu+\nu=p$. Denote by $\pairUmn$ and $A$ the 
handlebody-kont $\pairAtau$ and the type $3$-$3$
annulus obtained by
the construction in Section \ref{subsec:examples_irreatoro}, respectively.

The isotopy in Fig.\ \ref{fig:equivalence_to_m5_1} shows that 
a regular neighborhood of $8_{16}\cup \tau$ in $\sphere$
is equivalent to the mirror image $(\sphere, \mirror 5_1)$ of $(\sphere, 5_1)$ 
in the handlebody-knot table \cite{IshKisMorSuz:12}.
Thus by Corollary \ref{cor:irreducibility:irreducible_HKA} 
and Lemma \ref{lm:atoroidality_general},
$\pairUmn$
is irreducible and atoroidal, for every $\mu,\nu\in\mathbb{Z}$.

Now, in terms of the meridional basis 
of $H_1(\ComplUmn)$ 
given by $x_1,x_2$ in Fig.\ \ref{fig:8_16_w_tunnel}
and with the orientation of $A$ in Fig.\ \ref{fig:p_annulus_8_16_tau},
we have $[l_+]=(\mu,\nu)$,
$[l_-]=(\mu-1,\nu+1)$.
Hence for every member $\pair$ in the handlebody-knot family 
\[\mathcal{U}:=\{\pairUmn\mid \mu>\nu+1>1 \text{ or } 0>\mu>\nu+1\}.\]  
the meridional basis is normalized, 
the slope type of $A$ is not $(\frac{p+1}{2},\frac{p-1}{2})$, and
none of $l_+,l_-$ represents the $\absp$-th multiple of some generator
of $H_1(\ComplHKA)$.
Thus up to isotopy, $A\subset\ComplHK$ is the unique type $3$-$3$ annulus 
by Theorem \ref{teo:uniqueness_HKA_irreducible}, and 
$\Sym\HK=1$ by Corollary \ref{cor:non_invertiable_no_symmetries}.
In addition, $\mathcal{U}$ is an infinite family.



\end{document}